\documentclass[12pt]{article}

\usepackage[applemac]{inputenc}
\usepackage{amsmath,amssymb,fullpage}

\newcommand{\tsum}{\sum}
\newcommand{\dint}{\int}
\newcommand{\tint}{\int}

\newtheorem{example}{Example}[section]
\newtheorem{theorem}[example]{Theorem}

\newtheorem{definition}[example]{Definition}
\newtheorem{lemma}[example]{Lemma}

\newtheorem{problem}[example]{Problem}
\newtheorem{remark}[example]{Remark}

\usepackage{amsfonts}
\usepackage{makeidx}
\usepackage{graphicx}
\usepackage{multicol}
\usepackage[bottom]{footmisc}

\numberwithin{equation}{section}

\def\1B{\text{1\!\!I}}

\begin{document}

\title{Introduction to White Noise, Hida-Malliavin Calculus and Applications}
\author{Nacira AGRAM$^{1}$ and Bernt ØKSENDAL$^{2}$}
\date{10 April 2019}
\maketitle

\footnotetext[1]{%
Department of Mathematics, Linnaeus University (LNU), Sweden.\newline
Email: \texttt{nacira.agram@lnu.se}}

\footnotetext[2]{%
Department of Mathematics, University of Oslo, Norway. \newline
Email: \texttt{oksendal@math.uio.no.}
\par
This research was carried out with support of the Norwegian Research
Council, within the research project Challenges in Stochastic Control,
Information and Applications (STOCONINF), project number 250768/F20.}

\begin{abstract}
This paper is based on lectures given by one of us (N.Agram) at the CIMPA Research School on Stochastic Analysis and Applications at Saida, Algeria, 1-9 March 2019. The purpose of the paper is threefold:

\begin{itemize}
\item We first give a short and simple survey of the Hida white noise calculus, and in
this context we introduce the Hida-Malliavin derivative as a stochastic
gradient with values in the Hida stochastic distribution space $(\mathcal{S}%
)^*$. We show that this Hida-Malliavin derivative defined on $L^2(\mathcal{F}%
_T,P)$ is a natural extension of the classical Malliavin derivative defined
on the subspace $\mathbb{D}_{1,2}$ of $L^2(P)$.

\item Second, the Hida-Malliavin calculus allows us to prove new results under
weaker assumptions than could be obtained by the classical theory. In
particular, we prove the following:\newline
(i) A general integration by parts formula and duality theorem for Skorohod
integrals,\newline
(ii) a generalised fundamental theorem of stochastic calculus, and\newline
(iii) a general Clark-Ocone theorem, valid for all $F \in L^2(\mathcal{F}%
_T,P).$
\end{itemize}

Thirdly, we present applications of the above theory. For example, we discuss the following:

\begin{itemize}
\item A general representation theorem for backward stochastic differential
equations with jumps, in terms of Hida-Malliavin derivatives,

\item a general stochastic maximum principle for optimal control,

\item backward stochastic Volterra integral equations,

\item optimal control of stochastic Volterra integral equations and other
stochastic systems.
\end{itemize}
\end{abstract}

\paragraph{MSC(2010):}

60H05, 60H20, 60J75, 93E20, 91G80,91B70.

\paragraph{Keywords:}

Brownian motion; chaos expansion; classical Malliavin derivative; white
noise probability space, Hida space of stochastic distributions,
Hida-Malliavin calculus, optimal control of stochastic Volterra integral
equations.

\section{White noise theory and Hida-Malliavin calculus \newline
with applications - and without tears ...}

\label{intro} The purpose of these lectures is to give a short and easy
introduction to the Hida white noise theory and the associated
Hida-Malliavin calculus. This theory is important for many applications and
we believe it deserves to be better known. The problem has been that most of
the literature in this area has been too general and formidable for the
reader who just wants to know enough of the basic features of the theory in
order to be able to apply it to his or her area of research. These lectures
aim to fill that gap in the literature. Moreover, they intend to convince
the reader that there are indeed a lot of applications of this theory, and
to explain where and how.

This is basically a survey paper. Most of the text in these lectures are
taken from other published sources referred to in the reference list.
However, some of the proofs are new and shorter than the originals.\newline

The stochastic calculus of variations, now also know as Malliavin calculus,
was introduced by P. Malliavin \cite{M} as a tool for studying the
smoothness of densities of solutions of stochastic differential equations.
Subsequently other applications of this theory was found. In \cite{O} Ocone
used Malliavin calculus to prove an explicit representation theorem for
Brownian motion functionals and in a subsequent paper \cite{OK} Karatzas and
Ocone applied this to study portfolio problems in finance. \newline
The original presentation of Malliavin was quite complicated, but
subsequently simpler constructions of this theory have been found. See e.g. 
\cite{DOP} and the references therein. In particular, we think that the use
of white noise theory makes the theory of Malliavin calculus (in this
context also known as Hida-Malliavin calculus) quite natural within the
context of directional derivatives and Fr\' echet derivatives on the space $%
\mathcal{S}^{^{\prime }}$ of tempered distributions, both in the Brownian
motion case and in the case of Poisson random measure. See definitions below.

A major advantage of presenting the Malliavin calculus in the context of
white noise theory is that the corresponding Hida-Malliavin derivative can
be extended from the subspace $\mathbb{D}_{1.2}$ to all of $L^2(P)$. This
enables us to prove stronger results compared with what would be possible in
the classical setting. In particular, we will prove\newline
(i) A general integration by parts formula and duality theorem for Skorohod
integrals,\newline
(ii) a generalised fundamental theorem of stochastic calculus,\newline
(iii) a general Clark-Ocone theorem, valid for all $F \in L^2(\mathcal{F}%
_T,P),$\newline
(iv) a general representation for solutions of backward stochastic
differential equations (BSDEs) with jumps, in terms of Hida-Malliavin
derivatives, and (v) applications to stochastic control.

\section{White noise theory for Brownian motion}

\label{sec2} In this section we give a short introduction to the Hida white
noise calculus. A general reference for this section is \cite{HOUZ}. See
also \cite{HKPS}.

\subsection{ The white noise probability space}

We start with the construction of the \emph{white noise probability space}.
Let $\mathcal{S}=\mathcal{S}(\mathbb{R}^d)$\label{simb-028} be \textit{the
Schwartz space of rapidly decreasing smooth $C^{\infty}(\mathbb{R}^d)$ real
functions} 
on $\mathbb{R}^d$. The space $\mathcal{S}=\mathcal{S}(\mathbb{R}^d)$ is a
Fr\'echet space with respect to the family of seminorms:\label{simb-029} 
\begin{equation*}
\Vert f \Vert_{K,\alpha} := \sup_{x \in \mathbb{R}^d}\big\{ (1+|x|^K) \vert
\partial^\alpha f(x)\vert \big\},
\end{equation*}
where $K = 0,1,...$, $\alpha=(\alpha_1,...,\alpha_d)$ is a multi-index with $%
\alpha_j= 0,1,...$ $(j=1,...,d)$ and\label{simb-030} 
\begin{equation*}
\partial^\alpha f := \frac{\partial^{|\alpha|}}{\partial
x_1^{\alpha_1}\cdots \partial x_d^{\alpha_d}}f
\end{equation*}
for $|\alpha|=\alpha_1+ ... +\alpha_d$. 

Let $\mathcal{S}^{\prime }=\mathcal{S}^{\prime }(\mathbb{R}^{d})$\label%
{simb-031} be its dual, called the space of \emph{tempered distributions}. 
\index{tempered distributions} Let $\mathcal{B}$ denote the family of all
Borel subsets of $\mathcal{S}^{\prime }(\mathbb{R}^{d})$ equipped with the
weak* topology. If $\omega \in \mathcal{S}^{\prime }$ and $\phi \in \mathcal{%
S}$ we let \label{simb-033} 
\begin{equation}
\omega (\phi )=\langle \omega ,\phi \rangle  \label{3.1}
\end{equation}%
denote the action of $\omega $ on $\phi $. For example, if $\omega =m$ is a
measure on $\mathbb{R}^{d}$ then 
\begin{equation*}
\langle \omega ,\phi \rangle =\int\limits_{\mathbb{R}^{d}}\phi (x)dm(x),
\end{equation*}%
and, in particular, if this measure $m$ is concentrated on $x_{0}\in \mathbb{%
R}^{d}$, then 
\begin{equation*}
\langle \omega ,\phi \rangle =\phi (x_{0})
\end{equation*}%
is the evaluation of $\phi $ at $x_{0}\in \mathbb{R}^{d}$.\newline
Other examples include 
\begin{equation*}
\left\langle \omega ,\phi \right\rangle =\phi ^{\prime }(x_{1}).
\end{equation*}%
i.e. $\omega $ takes the derivative of $\phi $ at a point $x_{1}$.%
\newline
Or, more generally, 
\begin{equation*}
\left\langle \omega ,\phi \right\rangle =\phi ^{(k)}(x_{k}),
\end{equation*}%
i.e. $\omega $ takes the $k$'th derivative at the point $x_{k}$,
or linear combinations of the above.\newline

From now on we consider only the $1$-dimensional case,
i.e. $d=1$. For a multidimensional presentation see \cite{HOUZ}. We fix the
sample space to be $\Omega =\mathcal{S}^{\prime }(\mathbb{R})=\mathcal{S}%
^{\prime }$ and $\mathcal{F}=\mathcal{B}$. In the following we will use the
Bochner-Minlos-Sazonov theorem%
\index{Bochner-Minlos-Sazonov theorem} (see e.g. \cite{HOUZ}), which in our
setting states the following:

\begin{theorem}
{(Bochner-Minlos-Sazonov)}. \label{bms} Let $g: \mathcal{S} \mapsto 
\mathbb{R}$ be given. Then there exists a probability measure $\mu $ on $%
\Omega =\mathcal{S}^{\prime }(\mathbb{R})$ such that 
\begin{equation}
\mathbb{E}_{\mu }[e^{i\left\langle \omega ,\phi \right\rangle
}]:=\int_{\Omega }e^{i\left\langle \omega ,\phi \right\rangle }d\mu (\omega
)=g(\phi );%
\text{ for all }\phi \in \mathcal{S}
\end{equation}%
if and only if the function $g$ satisfies the following 3 conditions:\newline
(i) $g(0)=1$\newline
(ii) $g$ is continuous in the Fr\' echet topology on $\mathcal{S}$\newline
(iii) $g$ is positive definite, i.e. 
\begin{equation} \label{posdef}
\sum_{j,\ell }^{n}z_{j}\bar{z}_{\ell }g(\phi _{j}-\phi _{\ell })\geq 0\quad 
\text{ for all }z_{j}\in \mathbb{C},\phi _{j}\in \mathcal{S}; \quad j=1,2,...n,
\end{equation}%
where $\mathbb{C}$ denotes the set of complex numbers.
\end{theorem}

In particular, if we choose 
\begin{equation}
g(\phi )=e^{-\frac{1}{2}||\phi ||^{2}};\quad \phi \in \mathcal{S},
\end{equation}%
we can check that $g$ satisfies the conditions (i) - (iii) in the above
theorem, and hence we get that there exists a probability measure $P$ on $%
\Omega $ such that \label{simb-501} 
\begin{equation}
\mathbb{E}[e^{i\langle \omega ,\phi \rangle }]:=\int\limits_{\Omega
}e^{i\langle \omega ,\phi \rangle }P(d\omega )=e^{-\frac{1}{2}\Vert \phi
\Vert ^{2}},\text{ }\phi \in \mathcal{S},  \label{boch-min}
\end{equation}%
where 
\begin{equation*}
\Vert \phi \Vert ^{2}=\Vert \phi \Vert _{L^{2}(\mathbb{R})}^{2}=\int\limits_{%
\mathbb{R}}|\phi (x)|^{2}dx.
\end{equation*}%
The measure $P$ is called the \emph{white noise probability measure}%
\index{white noise!probability measure} and $(\Omega ,\mathcal{F},P)=(%
\mathcal{S}^{\prime },\mathcal{B},P)$ is called \emph{the white noise
probability space}.%
\index{white noise!probability space}\label{simb-502}

\bigskip

\begin{definition}
\label{Def. 3.1} \textit{The (smoothed) white noise process\/} 
\index{white noise!smoothed}is the measurable map 
\begin{equation*}
w:\mathcal{S}\times \mathcal{S}^{\prime }\rightarrow \mathbb{R}
\end{equation*}%
given by \label{simb-032} 
\begin{equation}
w(\phi ,\omega )=w_{\phi }(\omega )=\langle \omega ,\phi \rangle , \quad
\phi \in \mathcal{S},%
\text{ }\omega \in \mathcal{S}^{\prime }.  \label{3.4}
\end{equation}
\end{definition}

From $w_{\phi }$ we can construct a Brownian motion process $B({t})$, $t \in 
\mathbb{R}$, as follows:

\begin{enumerate}
\item[\textbf{Step 1}] \bigskip First we verify that the isometry 
\begin{equation}
\mathbb{E}[w_{\phi }^{2}]=\Vert \phi \Vert ^{2},\;\phi \in \mathcal{S},
\label{3.5}
\end{equation}%
holds true where, according to our notation, the left-hand side is 
\begin{equation*}
\mathbb{E}[w_{\phi }^{2}]=\int\limits_{\mathcal{S}^{\prime }}\langle \omega
,\phi \rangle ^{2}P(d\omega ).
\end{equation*}

\item[\textbf{Step 2}] Next we use Step 1 to define the value $\langle
\omega ,\psi \rangle $ for arbitrary $\psi \in L^{2}(\mathbb{R})$, as $%
\langle \omega ,\psi \rangle :=\lim $ $\langle \omega ,\phi _{n}\rangle $,
where $\phi _{n}\in \mathcal{S}$, $n\in \mathbb{N}=\{1,2,...\}$, and $%
\phi _{n}\rightarrow \psi \;$in$\;L^{2}(\mathbb{R}).$\newline
By Step 1 it follows that this definition does not depend on the choice of
the approximating sequence $\{\phi _{n}\}_{n\in \mathbb{N}}$.

\item[\textbf{Step 3}] Using Step 2 we can define 
\begin{equation*}
\widetilde{B}(t,\omega ):=\langle \omega ,\chi _{\lbrack 0,t]}\rangle
,\qquad t\in \mathbb{R},
\end{equation*}%
by choosing 
\begin{equation*}
\psi (s)=\chi _{\lbrack 0,t]}(s)=\left\{ 
\begin{array}{cc}
1 & \text{if }s\in \lbrack 0,t)\text{ (or }s\in \lbrack t,0),\text{ if }t<0)
\\ 
0 & \text{ otherwise }%
\end{array}%
\right.
\end{equation*}%
which belongs to $L^{2}(\mathbb{R})$ for all $t\in \mathbb{R}$.\newline

\item[\textbf{Step 4}] By the Kolmogorov continuity theorem (see \cite{Loeve}%
) we obtain that $\widetilde{B}({t})$, $t \in \mathbb{R}$, has a \emph{%
continuous version} $B({t})$, $t\in \mathbb{R}$, i.e. for all $t$ we have $P%
\big\{\widetilde{B}({t}) = B({t})\big\}=1$. This continuous process $B(t)$, $%
t\in \mathbb{R}$, is a Brownian motion (Wiener process).
\end{enumerate}

\begin{remark}
Note that for each $t$ we are defining $\widetilde{B}(t)= \widetilde{B}%
(t,\omega)=\langle \omega,\chi_{[0,t]}\rangle$ by using a sequence of functions $%
\phi_n^{(t)} \in \mathcal{S}(\mathbb{R})$ converging to $\chi_{[0,t]}$ in $L^2(%
\mathbb{R})$. Hence $\widetilde{B}(t)$ is only defined almost everywhere on $%
\Omega$, where the exceptional set of measure zero depends on $t$. Since
there are uncountably many $t \in [0,\infty)$ there is no common set $%
\Omega_0$ of measure 0 in $\Omega$ such that $\widetilde{B}%
(t,\omega)=\langle \omega,\chi_{[0,t]} \rangle $ is defined for all $t \in [0,\infty)$ and
for all $\omega \in \Omega \setminus \Omega_0$. Therefore we cannot prove
the continuity of $t \mapsto \widetilde{B}(t,\omega) $ by arguing $\omega$%
-wise. But we can use the Kolmogorov continuity theorem to conclude that $%
\widetilde{B}(t,\omega)$ has a continuous version.
\end{remark}

Note that when the Brownian motion process $B(t,\omega )$, $t \in \mathbb{R}$%
, $\omega\in\Omega$ with $\Omega :=\mathcal{S}^{\prime }(\mathbb{R})$ is
constructed this way, then each $\omega \in \Omega =\mathcal{S}^{\prime }(%
\mathbb{R}) $ is a tempered distribution. Hence $\Omega$ is a Fr\' echet
space, i.e. a topological vector space with a topology given by a family of
seminorms. This gives us a topological structure on $\Omega$ which we will
use frequently in the following.\newline

From the above Step 2 it follows that the smoothed white noise $w_{\phi }$
can be extended to all (deterministic) $\phi \in L^{2}(\mathbb{R)}$ and that
the relation between smoothed white noise $w_{\phi }$ and the Brownian
motion process $B({t})$, $t\in\mathbb{R}$, is 
\begin{equation}
w_{\phi }(\omega )=\int\limits_{\mathbb{R}}\phi (t)dB(t,\omega ), \quad
\omega\in \Omega, \qquad \phi \in L^{2}(\mathbb{R)},  \label{3.8}
\end{equation}
where the integral on the right-hand side is the Wiener-Itô integral. Note
that the isometry \eqref{3.5} is then the classical \emph{It\^o isometry}.

\subsection{The Wiener-Itô chaos expansion}

We now present an orthogonal expansion of the space $L^2(\mathcal{F},P)$. 
This presentation will be useful for the extended Hida-Malliavin calculus we
come to later.\newline

The \emph{Hermite polynomials} $h_{n}(x)$ are defined by 
\begin{equation*}
h_{n}(x)=(-1)^{n}e^{\frac{1}{2}x^{2}}\frac{d^{n}}{dx^{n}}(e^{-\frac{1}{2}%
x^{2}})\;,\;n=0,1,2,\ldots  \label{3.9}
\end{equation*}
The first Hermite polynomials are 
\begin{eqnarray*}
&&h_{0}(x)=1,h_{1}(x)=x,h_{2}(x)=x^{2}-1,h_{3}(x)=x^{3}-3x \\
&&h_{4}(x)=x^{4}-6x^{2}+3,h_{5}(x)=x^{5}-10x^{3}+15x,\ldots
\end{eqnarray*}
Some useful properties of the Hermite polynomials are

\begin{itemize}
\item $h^{\prime }_{n} (x)=n h_{n-1}(x); \quad n= 1,2, ...$

\item $h_{n+1}(x) - 2x h_{n}(x) +2n h_{n-1}(x) = 0; \quad n=1,2, ....$
\end{itemize}

Let $e_{k}$\label{simb-034} be the $k$'th \emph{Hermite function}%
\index{Hermite function} defined by 
\begin{equation}
e_{k}(x):=\pi ^{-%
\frac{1}{4}}((k-1)!)^{-\frac{1}{2}}\,e^{-\frac{1}{2}x^{2}}h_{k-1}(\sqrt{2}%
x),\quad k=1,2,\ldots  \label{3.10}
\end{equation}%
Then $\{e_{k}\}_{k\geq 1}$ constitutes an orthonormal basis for $L^{2}(%
\mathbb{R)}$ and $e_{k}\in \mathcal{S}(\mathbb{R})$ for all $k$. 
\newline
Define 
\begin{equation}
\theta _{k}(\omega ):=\langle \omega ,e_{k}\rangle =w_{e_{k}}(\omega
)=\int\limits_{\mathbb{R}}e_{k}(x)dB(x,\omega ),\quad \omega \in \Omega .
\label{3.11}
\end{equation}

\begin{definition}
Let $\mathcal{J}$\label{simb-035} denote the set of all finite multi-indices 
$\alpha =(\alpha _{1},\alpha _{2},\ldots ,\alpha _{m})$, $m=1,2,\ldots$, of
non-negative integers $\alpha _{i}$. If $\alpha =(\alpha _{1},\cdots ,\alpha
_{m})\in \mathcal{J}$ , $\alpha \neq 0,$\ we put \label{simb-036} 
\begin{equation}  \label{3.12a}
H_{\alpha }(\omega ):=\prod\limits_{j=1}^{m}h_{\alpha _{j}}(\theta
_{j}(\omega)) = h_{\alpha_1}(\theta_1) h_{\alpha_2}(\theta_2) ...
h_{\alpha_m}(\theta_m), \quad \omega\in\Omega.
\end{equation}
\end{definition}

We set $H_{0}:=1$. 
Hereafter we put 
\begin{equation}
\epsilon ^{(k)}=(0,0,...,1,0,...,0)  \label{3.13bis}
\end{equation}%
with $1$ on $k$'th position. For example, we have 
\begin{equation*}
H_{\epsilon ^{(k)}}(\omega )=h_{1}(\theta _{k}(\omega ))=\theta _{k}=\langle
\omega ,e_{k}\rangle ,
\end{equation*}%
and if $\alpha =(3,0,2)$, then 
\begin{equation*}
H_{(3,0,2)}=h_{3}(\theta _{1})h_{0}(\theta _{2})h_{2}(\theta _{3})=(\theta
_{1}^{3}-3\theta _{1})(\theta _{3}^{2}-1).
\end{equation*}%
We have the following fundamental result: 

\begin{theorem}
\label{Th. 3.2} \textbf{The Wiener-Itô chaos expansion theorem. }%
\index{chaos expansion} The family\textbf{\ }$\{H_{\alpha }\}_{\alpha \in 
\mathcal{J}%
\text{ }}$ constitutes an orthogonal basis of $L^{2}(P)$. More precisely,
for all $\mathcal{F}$-measurable $X\in L^{2}(P )$ there exist (uniquely
determined) numbers $c_{\alpha }\in \mathbb{R}$ such that 
\begin{equation}
X=\sum\limits_{_{\alpha \in \mathcal{J}\text{ }}}c_{\alpha }H_{\alpha }
\quad \in L^{2}(P).  \label{3.14}
\end{equation}
Moreover, we have the isometry 
\begin{equation}
\Vert X\Vert _{L^{2}(P)}^{2}=\sum_{\alpha \in \mathcal{J}} \alpha !
c_{\alpha }^{2},  \label{3.15}
\end{equation}
where $\alpha !=\alpha _{1}!\alpha _{2}!\cdots \alpha _{m}!$ for $%
\alpha=(\alpha _{1},\alpha _{2},\ldots, \alpha _{m})$.
\end{theorem}

\begin{example}
\label{Example 6.3} To find the chaos expansion of the Brownian motion
process $B(t)$ at time $t$, we proceed as follows: \label{extra3.14} 
\begin{align}
&B(t) = \int_\mathbb{R}\chi_{[0,t]} (s) dB(s) = \int_\mathbb{R}%
\sum_{k=1}^{\infty} \big( \chi_{[0,t]}, e_k \big)_{L^2(\mathbb{R})} e_k(s)
dB(s)  \notag \\
&= \sum_{k=1}^{\infty} \Big( \int_0^t e_k(y)dy\Big) \int_\mathbb{R} e_k (s)
dB(s) = \sum_{k=1}^{\infty} \Big( \int_0^t e_k(y) dy \Big) %
H_{\epsilon^{(k)}}.
\end{align}
\end{example}

\subsection{The Hida stochastic test and distribution spaces}

Analogous to the test functions $\mathcal{S}(\mathbb{R})$ and the tempered
distributions $\mathcal{S}^{\prime }(\mathbb{R})$ on the real line $\mathbb{R%
}$, there is a useful space of \emph{(Hida) stochastic test functions} $(%
\mathcal{S})$ and a space of \emph{(Hida) stochastic distributions} \textsl{%
\ $(\mathcal{S})^{* }$} on the white noise probability space. In the
following we will use the notation\label{simb-037} 
\begin{equation}  \label{3.23bis}
\big( 2\mathbb{N} \big)^{\alpha} = \prod\limits_{j=1}^{m }(2j)^{\alpha
_{j}}=(2\cdot1)^{\alpha_1}(2\cdot 2)^{\alpha _2} (2\cdot 3)^{\alpha _3} ...
(2m)^{\alpha_m} , \quad \text{ for }\quad \alpha =
(\alpha_{1},...,\alpha_{m}) \in \mathcal{J}.
\end{equation}

\begin{definition}
{(Hida stochastic test function spaces $(\mathcal{S})_k, (\mathcal{S})$)} %
\label{Def. 3.4} Let $k\in\mathbb{R}$. We say that $f=\sum\limits_{\alpha
\in \mathcal{J}}a_{\alpha } H_{\alpha }\in L^{2}(P )$ belongs to the \emph{%
Hida test function Hilbert space} $(\mathcal{S})_{k}$ \label{simb-038}if %
\label{simb-039} 
\begin{equation}
\Vert f \Vert^2_k := \sum\limits_{\alpha \in \mathcal{J}}\alpha !a_{\alpha
}^{2}(2\mathbb{N})^{\alpha k}<\infty.  \label{3.23}
\end{equation}
We define the \emph{Hida stochastic test function space}%
\index{Hida test function space} $(\mathcal{S})$\label{simb-040} as the
space 
\begin{equation*}
(\mathcal{S}) = \bigcap_{k\in\mathbb{R}} (\mathcal{S})_{k}
\end{equation*}
equipped with the projective topology, 
\index{projective topology}, i.e.\newline
$f_{n} \longrightarrow f$, $n\to \infty$, in $(\mathcal{S})$ if and only if $%
\Vert f_{n} - f \Vert_{k} \longrightarrow 0$, $n\to \infty$, for all $k$. 
\newline
\end{definition}

We illustrate this concept with an example: 
\vspace{3mm}

\begin{example}
\label{Example 3.6} The smoothed white noise $w_{\phi }$\ belongs to $(%
\mathcal{S})$ if $\phi \in \mathcal{S}(\mathbb{R})$.\newline
In fact, if $\phi =\sum\limits_{j=1}^\infty c_{j}e_{j}$\ we have 
\begin{equation}
w_{\phi }=\sum\limits_{j=1}^\infty c_{j}H_{\epsilon ^{(j)}}.  \label{3.25}
\end{equation}
Therefore, using (\ref{3.23}) we can see that $w_{\phi }\in (\mathcal{S})$\
if and only if 
\begin{equation*}
\sum\limits_{j=1}^\infty c_{j}^{2}(2j)^{k}<\infty
\end{equation*}
for all $k$, which holds because $\phi \in \mathcal{S}(\mathbb{R})$. See
e.g. \cite{RS}.\newline
\end{example}

\begin{definition}
{(Hida stochastic distribution spaces $(\mathcal{S})_q, (\mathcal{S})^*$) }

\begin{itemize}
\item Let $q\in \mathbb{R}$. We say that the formal sum $F=\sum\limits_{%
\alpha \in \mathcal{J}} b_{\alpha }H_{\alpha }$ belongs to the \emph{Hida
distribution Hilbert space} $(\mathcal{S})_{-q}$\label{simb-041} if\label%
{simb-042} 
\begin{equation}
\Vert F \Vert^2_{-q} :=\sum\limits_{\alpha \in \mathcal{J}}\alpha
!c_{\alpha}^{2} (2\mathbb{N})^{-\alpha q}<\infty .  \label{3.24}
\end{equation}
We define the \emph{Hida stochastic distribution space}%
\index{Hida distribution space} $(\mathcal{S})^{*}$\label{simb-043} as the
space 
\begin{equation*}
(\mathcal{S})^{*} = \bigcup_{q\in\mathbb{R}} (\mathcal{S})_{-q}
\end{equation*}
equipped with the inductive topology,%
\index{inductive topology} i.e. $F_{n} \longrightarrow F$, $n\to \infty$, in 
$(\mathcal{S})^{*}$ if and only if there exists $q$ such that $\Vert F_{n} -
F \Vert_{-q} \longrightarrow 0$, $n\to \infty$.

\item If $F=\sum_{\alpha \in \mathcal{J}}b_{\alpha }H_{\alpha }\in (\mathcal{%
S})^{\ast }$, we define the \emph{generalized expectation }$\mathbb{E}$\emph{%
$[F]$%
\index{generalised expectation} of $F$} by 
\begin{equation}
\mathbb{E}\big[F\big]=b_{0}.  \label{Wick6.27}
\end{equation}%
(Note that if $F\in L^{2}(P)$ then the generalised expectation coincides
with the usual expectation, since $\mathbb{E}[H_{\alpha }]=0$ for all $\alpha \neq 0$%
).
\end{itemize}
\end{definition}

Note that $(\mathcal{S})^{* }$ can be regarded as the dual of $(\mathcal{S})$%
: Namely, the action of $F=\sum\limits_{\alpha }b_{\alpha }H_{\alpha }\in (%
\mathcal{S})^{* }$ on $f=\sum\limits_{\alpha }a_{\alpha }H_{\alpha }\in (%
\mathcal{S})$, where $b_{\alpha}, a_{\alpha} \in \mathbb{R}$, is given by 
\begin{equation*}
\langle F,f\rangle =\sum\limits_{\alpha }\alpha !a_{\alpha }b_{\alpha }.
\end{equation*}
We have the inclusions 
\begin{equation*}
(\mathcal{S}) \subset (\mathcal{S})_{k} \subset L^{2}(P)\subset(\mathcal{S}%
)_{-q} \subset (\mathcal{S})^{* }, \quad 
\text{ for all }\quad k,q.
\end{equation*}

\begin{example}
The \emph{singular} (also called \emph{pointwise}) \emph{white noise} 
\index{white noise!singular} $\overset{\bullet }{B} (t)$, $t\in\mathbb{R}$,%
\label{simb-044} is defined as follows: 
\begin{equation}
\overset{\bullet }{B}({t}):=\sum\limits_{k=1}^{\infty}e_{k}(t)H_{\epsilon
^{(k)}}.  \label{3.26}
\end{equation}
We can verify that $\overset{\bullet }{B}({t})\in (\mathcal{S})^{* }$ for
all $t$, as follows:\newline
\begin{equation*}
\Vert\overset{\bullet }{B}({t})\Vert^2_{-q} = \sum_{k=0}^\infty e^2_k(t)
\epsilon^{(k)}! \big( (2\mathbb{N})^{\epsilon^{(k)}}\big)^{-q}
=\sum_{k=0}^\infty e^2_k(t) \big( 2k\big)^{-q} < \infty, \qquad q\geq 2,
\end{equation*}
because 
\begin{equation*}  \label{extra3.25}
\sup_{t\in\mathbb{R}} | e_k(t)| = \mathcal{O} \big( k^{-1/12} \big).
\end{equation*}
Similarly we see that by \eqref{extra3.14} we get 
\begin{equation}  \label{extraextra3.25}
\frac{d}{dt} B(t) = \frac{d}{dt} \sum_{k=1}^{\infty} \Big( \int_0^t e_k(y)dy %
\Big) H_{\epsilon^{(k)}} = \overset{\bullet }{B}(t),
\end{equation}
where the derivative is taken in $(\mathcal{S})^*$.\newline
Thus we see that although the time derivative $\frac{d}{dt} B(t)$ does not
exist in the classical sense, it does exist as an element of $(\mathcal{S}%
)^{*}$, and this derivative is the singular white noise $\overset{\bullet }{B%
}(t)$.
\end{example}

The following useful characterisation of the Hida spaces is due to Zhang 
\cite{Z1}:

\begin{theorem}
(i) The Hida stochastic test function space $(\mathcal{S})$ consists of
those $F=\sum\limits_{\alpha \in \mathcal{J}}c_{\alpha}H_{\alpha} \in L^2(P)$
such that 
\begin{equation}  \label{zang1}
\sup_{\alpha} \{ c_{\alpha}^2 \alpha ! (2\mathbb{N})^{k \alpha} \}< \infty 
\text{ for \emph{all } } k \in \mathbb{N}.
\end{equation}
\medskip (ii) The Hida stochastic distribution space $(\mathcal{S})^*$
consists of those formal expansions $F=\sum\limits_{\beta \in \mathcal{J}%
}b_{\beta}H_{\beta} $ such that 
\begin{equation}  \label{zang2}
\sup_{\beta} \{ c_{\beta}^2 \beta ! (2\mathbb{N})^{-q \beta} \}< \infty 
\text{ for \emph{some} } q \in \mathbb{N}.
\end{equation}
\end{theorem}

We will also need the following result:

\begin{theorem}
\begin{equation}  \label{zang3}
\sum\limits_{\alpha \in \mathcal{J}} (2\mathbb{N})^{-\alpha q}<\infty \text{
if and only if } q > 1.
\end{equation}
\end{theorem}

\section{The Wick product}

\label{Wick-Hermite} In addition to a canonical vector space structure, the
spaces $(\mathcal{S})$ and $(\mathcal{S})^{* }$ also have a natural
multiplication given by the \emph{Wick product}.\\

The main reference for this
section is \cite{HOUZ}.

\begin{definition}
\label{Def. 3.6} If $X=\sum\limits_{\alpha }a_{\alpha }H_{\alpha }\in (%
\mathcal{S})^{* },Y=\sum\limits_{\beta }b_{\beta }H_{\beta }\in (\mathcal{S}%
)^{* }$ then the \textsl{Wick product} 
\index{Wick product}$X\diamond Y$ \label{simb-045}of $X$ and $Y$ is defined
by 
\begin{equation}
X\diamond Y :=\sum\limits_{\alpha ,\beta }a_{\alpha }b_{\beta }H_{\alpha
+\beta }=\sum\limits_{\gamma }(\sum\limits_{\alpha +\beta =\gamma }a_{\alpha
}b_{\beta })H_{\gamma }.  \label{3.27}
\end{equation}
\end{definition}

Using (\ref{3.24}) and (\ref{3.23}) one can now verify the following: 
\begin{equation}
X,Y\in (\mathcal{S})^{* }\Rightarrow X\diamond Y\in (\mathcal{S})^{* }.
\label{3.28}
\end{equation}%
\begin{equation}
X,Y\in (\mathcal{S})\Rightarrow X\diamond Y\in (\mathcal{S}).  \label{3.29}
\end{equation}%
See \cite[Lemma 2.4.4]{HOUZ}. Note, however, that $X,Y\in
L^{2}(P)\not\Rightarrow X\diamond Y\in L^{2}(P)$ in general. See \cite[%
Example 2.4.8]{HOUZ}.

\begin{example}
\label{Ex. 3.7}

\begin{itemize}
\item[(i)] \bigskip The Wick square of the singular white noise is 
\begin{equation*}
(\overset{\bullet }{W})^{\diamond 2}({t}) =\sum\limits_{k,m=1}^\infty
e_{k}(t)e_{m}(t)H_{\epsilon^{(k)}+\epsilon ^{(m)}}.
\end{equation*}
One can show that 
\begin{equation}  \label{3.29bis}
(\overset{\bullet }{W})^{\diamond 2}({t}) \in (\mathcal{S})^{*}, \quad t\in%
\mathbb{R}.
\end{equation}

\item[(ii)] \bigskip The Wick square of the smoothed white noise is 
\begin{equation*}
(w_{\phi })^{\diamond 2}=\sum\limits_{k,m=1}^{\infty }c_{k}c_{m}H_{\epsilon
^{(k)}+\epsilon ^{(m)}}\quad 
\text{\emph{if}}\quad \phi =\sum_{k=1}^{\infty }c_{k}e_{k}\in L^{2}(\mathbb{%
R)}
\end{equation*}%
Since 
\begin{equation*}
H_{_{\epsilon ^{(k)}+\epsilon ^{(m)}}}=\left\{ 
\begin{array}{cc}
H_{\epsilon ^{(k)}}\cdot H_{\epsilon ^{(m)}} & \text{\emph{if} }k\neq m \\ 
H_{\epsilon ^{(k)}}^{2}-1 & \text{\emph{if} }k=m%
\end{array}%
\right.
\end{equation*}%
we see that 
\begin{equation}
(w_{\phi })^{\diamond 2}=w_{\phi }^{2}-\sum_{k=1}^{\infty }c_{k}^{2}=w_{\phi
}^{2}-\Vert \phi \Vert ^{2}.  \label{3.29bisbis}
\end{equation}%
Note, in particular, that $(w_{\phi })^{\diamond 2}$ is not positive. In
fact, $\mathbb{E}[(w_{\phi })^{\diamond 2}]=0$\ by (\ref{3.26}) and the fact
that $\mathbb{E}[H_{\alpha }]=0$ for $\alpha \neq 0$ (see Theorem \ref{Th.
3.2}).
\end{itemize}
\end{example}

Before proceeding further, we list some reasons that the Wick product is
natural to use in stochastic calculus:

\begin{description}
\item[a)] First, note that if (at least) one of the factors $X,Y$ is
deterministic, then 
\begin{equation*}
X\diamond Y=X\cdot Y
\end{equation*}%
Therefore the two types of products, the Wick product and the ordinary ($%
\omega $-pointwise) product, coincide in the deterministic calculus. So when
one extends a deterministic model to a stochastic model by introducing
noise, it is not obvious which interpretation to choose for the products
involved. The choice should be based on additional modelling and
mathematical considerations.

\item[b)] The Wick product is the only product which is defined for singular
white noise $\overset{\bullet }{B}.$ Pointwise product $X\cdot Y$ does not
make sense in $(\mathcal{S})^{* }$!

\item[c)] The Wick product has been used for 50 years already in quantum
physics as a renormalization procedure.

\item[d)] There is a fundamental relation between Itô/Skorohod integrals and
Wick products, given by 
\begin{equation}
\int_\mathbb{R} Y({t})\delta B(t)= \int_\mathbb{R} Y({t})\diamond \overset{%
\bullet}{B}({t})dt  \label{3.32}
\end{equation}%
Here the integral on the right is interpreted as a Bochner integral with
values in $(\mathcal{S})^{* }$. See Theorem \ref{Th5.12} below.

\item[e)] A big class of strong solutions to stochastic differential
equations can be explicitly solved by using the Wick product. See \cite{LaP}.
\end{description}

\subsection{Some basic properties of the Wick product}

We list below some useful properties of the Wick product. Some are easy to
prove, others harder. For complete proofs see \cite{HOUZ}.

\medskip The Wick product is a binary operation on $(\mathcal{S})^{\ast }$,
i.e. for all $X,Y\in (\mathcal{S})^{\ast }$ we have $X\diamond Y\in (%
\mathcal{S})^{\ast }$. Moreover, for arbitrary $X,Y,Z\in (\mathcal{S})^{\ast
}$ we have 
\begin{eqnarray}
&&X\diamond Y=Y\diamond X\qquad \text{(commutative law)},  \label{3.46} \\
&&X\diamond (Y\diamond Z)=(X\diamond Y)\diamond Z\qquad \text{(associative
law)},  \label{3.47} \\
&&X\diamond (Y+Z)=(X\diamond Y)+(X\diamond Z)\qquad \text{(distributive law)}.
\label{3.48}
\end{eqnarray}%
In view of the above we can define the \emph{Wick powers}%
\index{Wick power} 
\begin{equation*}
X^{\diamond n}=X\diamond X\diamond \cdots \diamond X\quad (n%
\text{ times})\quad \text{for}\quad X\in (\mathcal{S})^{\ast },\quad
n=1,2,...\,.
\end{equation*}%
We put $X^{\diamond 0}=1$. Similarly, the \emph{Wick exponential}%
\index{Wick exponential} of $X\in (\mathcal{S})^{\ast }$ is defined by 
\begin{equation}
\exp ^{\diamond }X=\sum_{n=0}^{\infty }%
\frac{1}{n!}X^{\diamond n},  \label{wick-exponential-brownian}
\end{equation}%
if convergent in $(\mathcal{S})^{\ast }$. Thus the Wick algebra obeys the
same rules as the ordinary algebra. For example, 
\begin{equation}
(X+Y)^{\diamond 2}=X^{\diamond 2}+2X\diamond Y+Y^{\diamond 2}  \label{3.49}
\end{equation}%
(no Itô formula!) and 
\begin{equation}
\exp ^{\diamond }(X+Y)=\exp ^{\diamond }(X)\diamond \exp ^{\diamond }(Y).
\label{3.50}
\end{equation}%
Note, however, that combinations of ordinary products and Wick products
require caution. For example, in general we have 
\begin{equation*}
X\cdot (Y\diamond Z)\not=(X\cdot Y)\diamond Z\,.
\end{equation*}%
Note that since $E[H_{\alpha }]=0$ for all $\alpha \neq 0$, we have that if $%
X=\sum_{\alpha \in \mathcal{J}}c_{\alpha }H_{\alpha }\in L^{2}(P)$, then 
\begin{equation*}
\mathbb{E}[X]=c_{0}.
\end{equation*}%
From this we deduce the remarkable property of the Wick product: 
\begin{equation}
\mathbb{E}[X\diamond Y]=\mathbb{E}[X]\cdot \mathbb{E}[Y],  \label{3.51}
\end{equation}%
whenever $X,Y$ and $X\diamond Y$ are $P$-integrable. Note that it is \textit{%
not} required that $X$ and $Y$ are independent!

By induction it follows from \eqref{3.51} that 
\begin{equation}
\mathbb{E}[\exp ^{\diamond }X]=\exp \mathbb{E}[X].  \label{3.51bis}
\end{equation}%
From Example \ref{Ex. 3.7} (ii) we deduce that 
\begin{equation*}
w_{\phi }\diamond w_{\psi }=w_{\phi }\cdot w_{\psi }-{\frac{1}{2}}%
\int\limits_{\mathbb{R}}\phi (t)\psi (t)dt,\text{ }\phi ,\psi \in L^{2}(%
\mathbb{R})\,.
\end{equation*}%
In particular, 
\begin{equation}
B^{\diamond 2}(t)=B^{2}({t})-t,\qquad t\geq 0.  \label{3.52}
\end{equation}%
Moreover, if $\mathrm{supp\,}\phi \cap \mathrm{supp\,}\psi =\emptyset $,
then 
\begin{equation}
w_{\phi }\diamond w_{\psi }=w_{\phi }\cdot w_{\psi }.  \label{3.53}
\end{equation}%
Hence if $0\leq t_{1}\leq t_{2}\leq t_{3}\leq t_{4}$, then 
\begin{equation}
(B({t_{4}})-B({t_{3}}))\diamond (B({t_{2}})-B({t_{1}}))=(B({t_{4}})-B({t_{3}}%
))\cdot (B({t_{2}})-B({t_{1}})).  \label{3.54}
\end{equation}%
More generally, it can be proved that if $F$ is $\mathcal{F}_{t}$-measurable
and $h>0$, then 
\begin{equation}
F\diamond (B({t+h})-B({t}))=F\cdot (B({t+h})-B({t})).  \label{3.55}
\end{equation}%
For a proof see e.g. \cite[Exercise 2.22]{HOUZ}.

\subsection{Wick product and Hermite polynomials}

There is a striking connection between Wick powers and Hermite polynomials $%
h_n; n=0,1,2, ...$, as follows:

\begin{theorem}
\label{th3.3} (a) Choose $\varphi \in L^2(\mathbb{R})$.Then 
\begin{equation}
h_n(w_{\varphi}) = ||\varphi||^{-n} w_{\varphi} ^{\diamond n}.
\end{equation}
where $||\varphi|| = ||\varphi||_{L^2(\mathbb{R}) }$.\newline

(b) In particular, let $\theta_k = \int_{\mathbb{R}} e_k(s)dB(s); k=1,2, ...$
Then 
\begin{equation}
h_n (\theta_k) = \theta_k ^{\diamond n}: \quad n=0,1, ...
\end{equation}
\end{theorem}

\noindent {Proof.} \quad (a). The generating function of $h_n$ is given by 
\begin{equation}  \label{eq4.20}
\exp(tx -\frac{1}{2}t^2 )= \sum_{n=0}^{\infty} h_n(x) \frac{t^n}{n!}.
\end{equation}
On the other hand we know that 
\begin{equation}  \label{eq4.21}
\exp\Big(t\frac{w_{\varphi}}{||\varphi||} -\frac{1}{2}t^2 \Big)=\exp^{\diamond} \Big(t 
\frac{w_{\varphi} }{||\varphi||}\Big) = \sum_{n=0}^{\infty} \frac{t^n}{n!} \Big(
\frac{w_{\varphi}}{||\varphi||}\Big)^{\diamond n}.
\end{equation}
Substituting $x = \frac{w_{\varphi}}{||\varphi||}$ in \eqref{eq4.20} and
comparing the terms with equal power of $t$ in \eqref{eq4.21} we get that 
\begin{equation*}
h_n\Big(\frac{w_{\varphi}}{||\varphi||}\Big)=\Big(\frac{w_{\varphi}}{||\varphi||}\Big)^{\diamond n};
\quad \text{ for all } n.
\end{equation*}%
\newline
(b) This follows from (a) by using $\varphi=e_k$, using that $||e_k|| = 1$.
\hfill $\square$ \bigskip

Applying this result to the basis elements $H_{\alpha}$ defined in %
\eqref{3.12a} and using that Wick products and ordinary products are the
same for independent variables, we get

\begin{theorem}
\begin{align}
H_{\alpha} = \theta_1^{\diamond \alpha_1} \theta_2^{\diamond \alpha 2} .
...= \theta_1^{\diamond \alpha_1} \diamond \theta_2^{\diamond \alpha_3}
\diamond ...
\end{align}
\end{theorem}

\subsection{Wick products and Skorohod integration}

We now prove the fundamental relation \eqref{3.32} between Wick products and
Skorohod integration.

\begin{definition}
\label{definition3.13-added} A function $Y:\mathbb{R}\longrightarrow (%
\mathcal{S})^{\ast }$ (also called an $(\mathcal{S})^{\ast }$-valued
process) is \emph{$(\mathcal{S})^{\ast }$-integrable} if 
\begin{equation*}
\left\langle Y(t),f\right\rangle \in L^{1}(\mathbb{R}),\quad \text{ for all}%
\quad f\in (\mathcal{S}).
\end{equation*}%
Then the \emph{$(\mathcal{S})^{\ast }$-integral of $Y$}%
\index{$(\mathcal{S})^{\ast }$-integral}, denoted by $\int_{\mathbb{R}%
}Y(t)dt $, is the (unique) element in $(\mathcal{S})^{\ast }$ such that 
\begin{equation}
\left\langle \int_{\mathbb{R}}Y(t)dt,f\right\rangle =\int_{\mathbb{R}%
}\left\langle Y(t),f\right\rangle dt,\quad 
\text{ for all}\quad f\in (\mathcal{S}).  \label{3.68-added}
\end{equation}
\end{definition}

\begin{remark}
\label{remark3.14-added} The fact that \eqref{3.68-added} does indeed define 
$\int_{\mathbb{R}}Y(t)dt$ as an element of $(\mathcal{S})^{*}$ is a
consequence of \cite[Proposition 8.1]{HKPS}.
\end{remark}

Note that if $Y$ is $(\mathcal{S})^{*}$-integrable, then so is $%
Y\chi_{(a,b]} $, for all $a,b\in\mathbb{R}$, and we put 
\begin{equation*}
\int_{a}^{b} Y(t)dt := \int_{\mathbb{R}}Y(t)\chi_{(a,b]}(t)dt.
\end{equation*}

\begin{theorem}
{(\cite{DOP},Theorem 5.20)} \label{Th5.12}\newline
Assume that $\varphi (t)$ is $\mathbb{F}$-adapted and $\mathbb{E}[\int_{%
\mathbb{R}}\varphi (t)^{2}dt]<\infty $. Then $\varphi (t)\diamond \overset{%
\bullet }{B}(t)$ is integrable in $(\mathcal{S})^{\ast }$ and 
\begin{equation}
\int_{\mathbb{R}}\varphi (t)dB(t)=\int_{\mathbb{R}}\varphi (t)\diamond 
\overset{\bullet }{B}(t)dt.  \label{3.66}
\end{equation}
\end{theorem}

\label{rieman-skorohod} In particular, note that if $Y(t)= \sum_{i=1}^n c_i
\chi_{(t_i,t_{i+1}]}(t)$, $t \in \mathbb{R}$, with $c_i \in (\mathcal{S}%
)^{*} $ for $i=1,...,n$ and $t_1 < ... < t_n$. Then we have 
\begin{equation*}
\int_{\mathbb{R}} Y(t) \diamond \overset{\bullet}{B}(t)dt = \sum_{i=1}^n c_i
\diamond \big( B(t_{i+1}) - B(t_{i}) \big).
\end{equation*}

\vspace{2mm} In view of the above theorem, the following terminology is
natural.

\begin{definition}
\label{Wick-definition6.16} Suppose $Y$ is an $(\mathcal{S})^{*}$-valued
process such that 
\begin{equation*}
\int_{\mathbb{R}} Y(t)\diamond \overset{\bullet}{B}(t)dt \: \in (\mathcal{S}%
)^{*},
\end{equation*}
then we call this integral \emph{the generalised Skorohod integral of $Y$}%
\index{Skorohod integral!generalized}.
\end{definition}

Combining the properties above with the fundamental relation (\ref{3.32})
for Skorohod integration, we get a powerful calculation technique for
stochastic integration. First of all, note that, by (\ref{3.32}), 
\begin{equation}
\int\limits_{0}^{t}\overset{\bullet }{B}(s)ds=B({t}); \quad t \in [0,T].
\label{3.63}
\end{equation}
From this we deduce that 
\begin{equation}  \label{3.63a}
\frac{d}{dt}B(t) \text{ exists in } (\mathcal{S})^{*} \text{ and } \frac{d}{dt}%
B(t)=\overset{\bullet }{B}(t); \quad t \in [0,T].
\end{equation}

Moreover, using (\ref{3.47}) we get
\begin{equation}
\int\limits_{0}^{T}X\diamond Y({t})\diamond \overset{\bullet }{B}({t})dt
=X\diamond \int\limits_{0}^{T}Y({t})\diamond \overset{\bullet }{B}({t})dt,
\label{3.64}
\end{equation}%
if $X$ does not depend on $t$. Compare this with the fact that for Skorohod
integrals we generally have 
\begin{equation}
\int\limits_{0}^{T}X\cdot Y({t})\delta B(t)\not=X\cdot \int\limits_{0}^{T}Y({%
t})\delta B(t)\,,  \label{3.65}
\end{equation}%
even if $X$ does not depend on $t$.

\begin{example}
To illustrate the use of Wick calculus, let us consider the following: 
\begin{eqnarray*}
\int\limits_{0}^{T}B(t)\, \lbrack B({T})-B(t)]\delta B(t)
&=&\int\limits_{0}^{T}B({t})\diamond (B({T})-B({t}))\diamond \overset{%
\bullet }{B}({t})dt \\
&=&\int\limits_{0}^{T}B(t)\diamond B({T})\diamond \overset{\bullet }{B}({t}%
)dt-\int\limits_{0}^{T}B^{\diamond 2}({t})\diamond \overset{\bullet }{B}({t}%
)dt \\
&=&B({T})\diamond \int\limits_{0}^{T}B(t)\diamond \overset{\bullet }{B}({t}%
)dt -{\frac{1}{3}}B^{\diamond 3}({T}) \\
&=&{\frac{1}{6}}B^{\diamond 3}({T})=\frac{1}{6}[B^{3}({T})-3TB({T})],
\end{eqnarray*}%
where we have correspondingly used (\ref{3.54}), (\ref{3.48}), (\ref{3.64})
and Theorem \ref{th3.3}, keeping in mind that $|| \chi_{[0,T]}(\cdot)||_{L^2(\mathbb{R})} = T^{%
\frac{1}{2}}$.
\end{example}

We proceed to establish some useful properties of generalised Skorohod
integrals.

\begin{lemma}
\label{Wick-lemma6.16} Suppose $f\in (\mathcal{S})$ and $G(t)\in (\mathcal{S}%
)_{-q}$ for all $t\in{\mathbb{R}}$, for some $q\in \mathbb{N}$. Put 
\begin{equation*}
\hat q = q +\frac{1}{\log 2}.
\end{equation*}
Then 
\begin{equation*}
\int_{{\mathbb{R}}} \vert \langle G(t)\diamond \overset{\bullet}{B}(t),f \rangle \vert
dt \leq \Vert f \Vert_{\hat q}\Big( \int_{{\mathbb{R}}} \Vert G(t)
\Vert^{2}_{-q}dt \Big)^{1/2}.
\end{equation*}
\end{lemma}

\noindent {Proof.} \quad Suppose $G(t)=\sum_{\alpha \in \mathcal{J}%
}a_{\alpha }(t)H_{\alpha }$, $f=\sum_{\beta \in \mathcal{J}}b_{\beta
}H_{\beta }$. Then 
\begin{equation*}
\begin{split}
\left\langle G(t)\diamond \overset{\bullet }{B}(t),f\right\rangle =&
\langle\sum_{\alpha ,k}a_{\alpha }(t)e_{k}(t)H_{\alpha +\epsilon
^{(k)}},\sum_{\beta \in \mathcal{J}}b_{\beta }H_{\beta }\rangle \\
=& \sum_{\alpha ,k}a_{\alpha }(t)e_{k}(t)b_{\alpha +\epsilon ^{(k)}}(\alpha
+\epsilon ^{(k)})!\,.
\end{split}%
\end{equation*}%
Hence 
\begin{equation*}
\begin{split}
\left\vert \int_{{\mathbb{R}}}\left\langle G(t)\diamond \overset{\bullet }{B}%
(t),f\right\rangle dt \right\vert \leq & \sum_{\alpha ,k}|b_{\alpha +\epsilon
^{(k)}}|\alpha !(\alpha _{k}+1)\int_{{\mathbb{R}}}|a_{\alpha }(t)e_{k}(t)|dt
\\
\leq & \sum_{\alpha ,k}|b_{\alpha +\epsilon ^{(k)}}|\alpha !(\alpha _{k}+1)%
\Big(\int_{{\mathbb{R}}}a_{\alpha }^{2}(t)dt\Big)^{1/2} \\
\leq & \Big(\sum_{\alpha ,k}b_{\alpha +\epsilon ^{(k)}}^{2}(\alpha +\epsilon
^{(k)})!(2\mathbb{N})^{\hat{q}(\alpha +\epsilon ^{(k)})}\Big)^{1/2} \\
& \cdot \Big(\sum_{\alpha ,k}\Big(\int_{{\mathbb{R}}}a_{\alpha }^{2}(t)dt%
\Big)\alpha !(\alpha _{k}+1)(2\mathbb{N})^{-\hat{q}(\alpha +\epsilon ^{(k)})}%
\Big)^{1/2} \\
\leq & \Vert f\Vert _{\hat{q}}\Big(\sum_{\alpha ,k}\Big(\int_{{\mathbb{R}}%
}a_{\alpha }^{2}(t)dt\Big)\alpha !(\alpha _{k}+1)(2k)^{-\frac{\alpha _{k}}{%
\log 2}}(2\mathbb{N})^{-q\alpha }\Big)^{1/2} \\
\leq & \Vert f\Vert _{\hat{q}}\Big(\int_{{\mathbb{R}}}\Vert G(t)\Vert
_{-q}^{2}dt\Big)^{1/2}.
\end{split}%
\end{equation*}%
\hfill $\square $ \bigskip

\vspace{2mm} \noindent Using this result we obtain the following:

\begin{theorem}
\label{Wick-theorem6.17}

\begin{enumerate}
\item[(i)] Suppose $G:{\mathbb{R}} \mapsto (\mathcal{S})_{-q}$ satisfies 
\begin{equation*}
\int_{{\mathbb{R}}} \Vert G(t) \Vert^{2}_{-q} dt < \infty, \quad\text{ for
some } q\in\mathbb{N}.
\end{equation*}
Then 
\begin{equation*}
\int_{{\mathbb{R}}} G(t)\diamond \overset{\bullet}{B}(t )dt \quad\text{
exists in } (\mathcal{S})^{*}.
\end{equation*}

\item[(ii)] Suppose $F(t)$, $F_{n}(t)$, $n=1,2,...$, are elements of $(%
\mathcal{S})_{-q}$ for all $t\in{\mathbb{R}}$ and 
\begin{equation*}
\int_{{\mathbb{R}}} \Vert F_{n}(t) - F(t) \Vert^{2}_{-q}dt \longrightarrow
0, \qquad n\to\infty.
\end{equation*}
Then 
\begin{equation*}
\int_{{\mathbb{R}}} F_{n}(t) \diamond \overset{\bullet}{B}(t)dt
\longrightarrow \int_{{\mathbb{R}}} F (t) \diamond \overset{\bullet}{B}%
(t)dt, \quad n \to \infty,
\end{equation*}
in the $\text{weak}^{*}$-topology on $(\mathcal{S})^{*}$.
\end{enumerate}
\end{theorem}

\noindent {Proof.} \quad (i) The proof follows from Lemma \ref%
{Wick-lemma6.16} and Definition \ref{definition3.13-added}.\newline
(ii)\thinspace\ By Lemma \ref{Wick-lemma6.16} we have 
\begin{equation*}
\begin{split}
\left\vert \left\langle \int_{{\mathbb{R}}}\big(F_{n}(t)-F(t)\big)\diamond 
\overset{\bullet }{B}(t)dt,f\right\rangle \right\vert \leq & \int_{{\mathbb{R%
}}}\left\vert \left\langle \big(F_{n}(t)-F(t)\big)\diamond \overset{\bullet }%
{B}(t),f\right\rangle \right\vert dt \\
\leq & \Vert f\Vert _{\hat{q}}\int_{{\mathbb{R}}}\Vert F_{n}(t)-F(t)\Vert
_{-q}^{2}dt\longrightarrow 0,\quad n\rightarrow \infty .
\end{split}%
\end{equation*}%
\hfill $\square $ \bigskip

\section{The Hida-Malliavin calculus}

As in previous sections we assume that the Brownian motion $B(t)$, $t\in 
\mathbb{R}$, is constructed on the space $(\Omega ,\mathcal{B},P)$ with $%
\Omega =\mathcal{S}^{\prime }(\mathbb{R})$. Note that any $\gamma \in L^{2}(%
\mathbb{R})$ can be regarded as an element of $\Omega =\mathcal{S}^{\prime }(%
\mathbb{R})$ by the action 
\begin{equation*}
\left\langle \gamma ,\phi \right\rangle =\int_{\mathbb{R}}\gamma (t)\phi
(t)dt;\quad \phi \in \mathcal{S}(\mathbb{R}).
\end{equation*}

\subsection{The Hida-Malliavin derivative}

We are now ready to define the Hida-Malliavin derivative. A general
reference for this section is \cite{DOP}.

\begin{definition}
\label{direct.derivative}

\begin{description}
\item[(i)] Let $F\in L^{2}(P)$ and let $\gamma \in L^2(\mathbb{R})$ be
deterministic. Then the \emph{directional derivative of $F$} in $(\mathcal{S}%
)^*$ (respectively, in $L^{2}(P)$) in the direction $\gamma$ is defined by 
\begin{equation}  \label{6.1}
D_\gamma F(\omega) = \lim_{\varepsilon \to 0} \frac{1}{\varepsilon} \big[ %
F(\omega + \varepsilon \gamma) - F(\omega)\big]
\end{equation}
whenever the limit exists in $(\mathcal{S})^*$ (respectively, in $L^{2}(P)$).

\item[(ii)] Suppose there exists a function $\psi: \mathbb{R} \mapsto (%
\mathcal{S})^*$ (respectively, $\psi: \mathbb{R} \mapsto L^{2}(P)$) such
that 
\begin{equation}  \label{6.2}
\begin{split}
& \int_{{\mathbb{R}}} \psi(t) \gamma(t)dt \quad \text{ exists in } (%
\mathcal{S})^{*} \text{ (respectively, in } L^{2}(P)) \text{ and} \\
&D_\gamma F = \int_\mathbb{R} \psi(t) \gamma(t) dt, \quad \text{ for all }
\gamma \in L^2(\mathbb{R}).
\end{split}%
\end{equation}
Then we say that $F$ is \emph{Hida-Malliavin differentiable} in $(\mathcal{S}%
)^*$ (respectively, in $L^{2}(P)$) and we write 
\begin{equation*}
\psi (t) = D_t F, \quad t \in \mathbb{R}.
\end{equation*}
We call $D_tF$ \label{simb-053}the \emph{Hida-Malliavin derivative at $t$ in 
$(\mathcal{S})^*$ 
\index{Hida-Malliavin derivative}%
\index{Malliavin derivative}(respectively, in $L^{2}(P)$)} or the \emph{%
stochastic gradient} of $F$ at $t$.
\end{description}
\end{definition}

\vspace{3mm}

\begin{example}
\begin{enumerate}

\item[(i)] Suppose $F(\omega )=\left\langle \omega ,f\right\rangle =\int_{%
\mathbb{R}}f(t)dB(t)$, $f\in L^{2}(\mathbb{R})$. Then 
\begin{equation*}
D_{\gamma }F=\frac{1}{\varepsilon }\big[\left\langle \omega +\varepsilon
\gamma ,f\right\rangle -\left\langle \omega ,f\right\rangle \big]%
=\left\langle \gamma ,f\right\rangle =\int_{\mathbb{R}}f(t)\gamma (t)dt.
\end{equation*}%
Therefore $F$ is Hida-Malliavin differentiable and 
\begin{equation*}
D_{t}\Big(\int_{\mathbb{R}}f(t)dB(t)\Big)=f(t),\qquad t-a.a.
\end{equation*}

\item[(ii)] Let $F\in L^{2}(P)$ be Hida-Malliavin differentiable in $%
L^{2}(P) $ for a.a. $t$. Suppose that $\varphi \in C^1(\mathbb{R})$ and $%
\varphi^{\prime }(F) D_tF \in L^{2}(P\times \lambda)$. Then if $\gamma \in
L^2(\mathbb{R})$ we have 
\begin{equation*}
\begin{split}
D_\gamma \big( \varphi(F) \big) = &\lim_{\varepsilon \to 0} \frac{1}{%
\varepsilon} \big[ \varphi(F(\omega + \varepsilon \gamma)) -
\varphi(F(\omega))\big] \\
= & \lim_{\varepsilon \to 0} \frac{1}{\varepsilon} \big[ \varphi(F(\omega) +
\varepsilon D_\gamma F) - \varphi(F(\omega))\big] \\
= & \frac{1}{\varepsilon} \varphi^{\prime }(F(\omega)) \varepsilon D_\gamma
F = \varphi^{\prime }(F) D_\gamma F \\
= &\int_\mathbb{R} \varphi^{\prime }(F) D_t F \gamma(t) dt.
\end{split}%
\end{equation*}
This proves that $\varphi(F)$ is also Hida-Malliavin differentiable and we
have the \emph{chain rule} 
\begin{equation}  \label{chain.rule}
D_t\big(\varphi(F)\big) = \varphi^{\prime }(F) D_tF.
\end{equation}
\end{enumerate}
\end{example}

\medskip More generally, the same proof gives the following extension: 

\begin{theorem}
\textbf{(Chain rule)} \label{Theorem 7.2bis} Let $F_{1},...,F_{m}\in L^{2}(P)$
be Hida-Malliavin differentiable in $L^{2}(P)$. Suppose that $\varphi \in
C^{1}({\mathbb{R}}^{m})$, $D_{t}F_{i}\in L^{2}(P)$, for all $t\in {\mathbb{R}%
}$, and $\frac{\partial \varphi }{\partial x_{i}}(F)D_{\cdot }F_{i}\in
L^{2}(\lambda \times P)$ for $i=1,...,m$, where $F=(F_{1},...,F_{m})$. Then $%
\varphi (F)$ is Hida-Malliavin differentiable and 
\begin{equation}
D_{t}\varphi (F)=\sum_{i=1}^{m}\frac{\partial \varphi }{\partial x_{i}}%
(F)D_{t}F_{i}.  \label{7.4bis}
\end{equation}
\end{theorem}

\subsection{The general Hida-Malliavin derivative}

It is useful to note how the Hida-Malliavin derivative can be expressed in
terms of the Wiener-It\^o chaos expansion (see Theorem \ref{Th. 3.2}). To
this aim observe that from 
the chain rule \eqref{chain.rule} we have 
\begin{equation}  \label{6.6}
D_t H_\alpha = \sum_{k=1}^m \prod_{j\ne k} h_{\alpha_j}(\theta_j) \alpha_k
h_{\alpha_k -1}(\theta_k) e_k (t) = \sum_{k=1}^m \alpha_k e_k(t) H_{\alpha -
\epsilon^{(k)}}.
\end{equation}
In view of this, the following definition is natural:

\begin{definition}
\textbf{The general Hida-Malliavin derivative.} \label{HM}\newline
If $F=\sum_{\alpha \in \mathcal{J}}c_{\alpha }H_{\alpha }\in (\mathcal{S}%
)^{\ast }$ we define the \emph{Hida-Malliavin derivative} $D_{t}F$ of $F$ at 
$t$ in $(\mathcal{S})^{\ast }$ by the following expansion: 
\begin{equation}
D_{t}F=\sum_{\alpha \in \mathcal{J}}\sum_{k=1}^{\infty }c_{\alpha }\alpha
_{k}e_{k}(t)H_{\alpha -\epsilon ^{(k)}},  \label{hida7.7}
\end{equation}%
whenever this sum converges in $(\mathcal{S})^{\ast }$. We shall denote $%
Dom(D_{t})$\label{simb-054} the set of all $F\in (\mathcal{S})^{\ast }$ for
which the above series converges in $(\mathcal{S})^{\ast }$ for all $t$.
\end{definition}

The following result gives that in fact $Dom(D_t) = (\mathcal{S})^*$:

\begin{theorem}

(i) If $F \in (\mathcal{S})$ then $D_tF \in (\mathcal{S})$ for all $t$.%
\newline

(ii) If $F \in (\mathcal{S})^*$ then $D_tF \in (\mathcal{S})^*$ for all $t$.
\end{theorem}

\noindent {Proof.} \quad (ii) We prove only the second part; the proof of
the first part being similar:\newline
Suppose $F=\sum_{\alpha \in \mathcal{J}}c_{\alpha }H_{\alpha }\in (\mathcal{S%
})^{\ast }$. Then we know by \eqref{zang2} that there exists $q_{0}\in 
\mathbb{N}$ such that 
\begin{equation}
\alpha !c_{\alpha }^{2}(2\mathbb{N})^{-q_{0}\alpha }\leq 1\text{ for all }%
\alpha .  \label{zhang1}
\end{equation}%
We have to prove that 
\begin{equation*}
D_{t}F=\sum_{\alpha \in \mathcal{J}}\sum_{k=1}^{\infty }c_{\alpha }\alpha
_{k}e_{k}(t)H_{\alpha -\epsilon ^{(k)}}\in (\mathcal{S})^{\ast }.
\end{equation*}%
To this end, it suffices to prove that there exists $q_{1}\in \mathbb{N}$
such that 
\begin{equation}
(\alpha -\epsilon ^{(k)})!c_{\alpha }^{2}\alpha _{k}^{2}e_{k}^{2}(t)(2%
\mathbb{N})^{-(q_{0}+q_{1})(\alpha -\epsilon ^{(k)})}\leq 1.
\end{equation}%
Since $\{e_{k}(t)\}$ is a bounded family we get by \eqref{zhang1} that 
\begin{align*}
& (\alpha -\epsilon ^{(k)})!c_{\alpha }^{2}\alpha _{k}^{2}e_{k}^{2}(t)(2%
\mathbb{N})^{-(q_{0}+q_{1})(\alpha -\epsilon ^{(k)})}\leq C_{1}(\alpha
-\epsilon ^{(k)})!c_{\alpha }^{2}\alpha _{k}^{2}(2\mathbb{N}%
)^{-(q_{0}+q_{1})(\alpha -\epsilon ^{(k)})} \\
& =C_{1}\alpha !\frac{\alpha _{k}-1}{\alpha _{k}}c_{\alpha }^{2}\alpha
_{k}^{2}(2\mathbb{N})^{-q_{0}(\alpha -\epsilon ^{(k)})}(2\mathbb{N}%
)^{-q_{1}(\alpha -\epsilon ^{(k)})} \\
& =C_{1}\alpha !c_{\alpha }^{2}(2\mathbb{N})^{-q_{0}(\alpha -\epsilon
^{(k)})}(\alpha _{k}-1)(\alpha _{k})(2\mathbb{N})^{-q_{1}(\alpha -\epsilon
^{(k)})} \\
& \leq C_{1}\alpha !c_{\alpha }^{2}(2\mathbb{N})^{-q_{0}\alpha }(2k)^{q_{0}}%
\Big[(\alpha _{k}-1)(2(k-1))^{-q_{1}\alpha _{k-1}}\Big]\Big[\alpha
_{k}(2k)^{-q_{1}(\alpha _{k}-1)}\Big] \\
& \leq C_{1}\alpha !c_{\alpha }^{2}\Big[(\alpha _{k}-1)2^{-q_{1}\alpha
_{k-1}}\Big]\Big[\alpha _{k}(2k)^{-q_{1}(\alpha _{k}-1)+q_{0}}\Big]\mathbf{1}%
_{\alpha _{k}>1} \\
& \leq 1,
\end{align*}%
if $k>1$ and $q_{1}$ is large enough.

\hfill $\square$ \bigskip

We end this section by stating some crucial properties of the Hida-Malliavin
derivative:

\subsection{The fundamental theorem of stochastic calculus for $B(\cdot)$}

\begin{theorem}
\label{th4.6}\textbf{(Fundamental theorem)} \label{fundamental} Suppose that $%
\varphi \in L^{2}(\lambda \times P)$ is $\mathbb{F}$-adapted. Then 
\begin{equation}
\int_{\mathbb{R}}\varphi (s)\diamond \overset{\bullet}{B}(s)ds\in L^{2}(P),
\label{eq2.77}
\end{equation}%
and for all $t>0$ we have

\begin{align}  \label{eq2.68}
&D_t\Big( \int_{\mathbb{R}} \varphi(s) \diamond \overset{\bullet}{B}(s) ds \Big) = \int_{%
\mathbb{R}} D_t\varphi(s) \diamond \overset{\bullet}{B}(s) ds + \varphi(t) \\
&=\int_t^{\infty} D_t\varphi(s) \diamond \overset{\bullet}{B}(s) ds + \varphi(t).
\label{eq2.70}
\end{align}
\end{theorem}

\noindent {Proof.} \quad The first statement \eqref{eq2.77} follows from
Theorem \ref{Th5.12}.\newline
Recall that 
\begin{equation}
\overset{\bullet}{B}(s)= \sum_j e_j(s) H_{\epsilon^{(j)}}.
\end{equation}
Hence, if we assume that $\varphi$ has the expansion 
\begin{equation}
\varphi(s)= \sum_{\beta} c_{\beta}(s)H_{\beta},
\end{equation}
we get 
\begin{align}
&D_t\Big( \int_{\mathbb{R}} \varphi(s) \diamond \overset{\bullet}{B}(s) ds \Big)=D_t\Big( %
\int_{\mathbb{R}}(\sum_{\beta} c_{\beta}(s) H_{\beta}) \diamond (\sum_j
e_j(s) H_{\epsilon^{(j)}} )ds \Big)  \notag \\
&=D_t\Big( \int_{\mathbb{R}} \sum_{\beta,j} c_{\beta}(s) e_j(s) H_{\beta +
\epsilon^{(j)}} ds \Big)  \notag \\
&=\int_{\mathbb{R}} \sum_{\beta,j,k} c_\beta(s) e_j(s) e_k(t) (\beta_k +
\epsilon^{(j)}) H_{\beta+\epsilon^{(j)}-\epsilon^{(k)} }ds  \notag \\
&=\int_{\mathbb{R}} \sum_{\beta,j,k} c_\beta(s) e_j(s) e_k(t) \beta_k
H_{\beta+\epsilon^{(j)}-\epsilon^{(k)} }ds +\int_{\mathbb{R}}
\sum_{\beta,j,k} c_\beta(s) e_j(s) e_k(t) \epsilon^{(j)}
H_{\beta+\epsilon^{(j)}-\epsilon^{(k)} }ds  \notag \\
&=\int_{\mathbb{R}} \Big( \sum_{\beta,k} c_\beta(s) e_k(t) \beta_k
H_{\beta-\epsilon^{(k)} }\Big) \diamond \Big( \sum_j
e_j(s)H_{\epsilon^{(j)}} \Big) ds+\sum_{\beta,k} (c_{\beta},e_k)_{L^2(%
\mathbb{R})} e_k(t) H_{\beta}  \notag \\
&=\int_{\mathbb{R}} D_t \Big(\sum_{\beta} c_\beta(s) H_{\beta }\Big) %
\diamond \overset{\bullet}{B}(s) ds +\sum_{\beta} c_{\beta}(t)H_{\beta}  \notag \\
&=\int_{\mathbb{R}} D_t\varphi(s) \diamond \overset{\bullet}{B}(s) ds+ \varphi(t).
\end{align}
This proves \eqref{eq2.68}. 
\hfill $\square$ \bigskip

\subsection{A generalised Clark-Ocone theorem}

We can apply Theorem \ref{th4.6} to prove the following, which was first
obtained in \cite{AaOPU} (with a different proof). See also \cite{DOP}:

\begin{theorem}
\textbf{(Generalised Clark-Ocone theorem)} Let $F\in L^{2}(\mathcal{F}_{T},P)$%
. Then $D_{t}F\in (\mathcal{S})^{\ast }$ for all $t$, $\mathbb{E}[D_{t}F|%
\mathcal{F}_{t}]\in L^{2}(\lambda \times P)$ and 
\begin{equation}
F=\mathbb{E}[F]+\int_{0}^{T}\mathbb{E}[D_{t}F|\mathcal{F}_{t}]dB(t)
\end{equation}
\end{theorem}

\noindent {Proof.} \quad By the Itô representation theorem there exists a unique $%
\mathbb{F}$-adapted process $\varphi \in L^{2}(\lambda \times P)$ such that 
\begin{equation}
F=\mathbb{E}[F]+\int_{0}^{T}\varphi (s)dB(s).
\end{equation}%
Taking the Hida-Malliavin derivative and conditional expectation of both
sides and applying the fundamental theorem (Theorem \ref{th4.6}) we get 
\begin{align}
\mathbb{E}[D_{t}F|\mathcal{F}_{t}]& =\mathbb{E}[D_{t}\big(%
\int_{0}^{T}\varphi (s)dB(s)\big)|\mathcal{F}_{t}]=\mathbb{E}%
[\int_{0}^{T}D_{t}\varphi (s)dB(s)+\varphi (t)|\mathcal{F}_{t}]  \notag \\
& =\mathbb{E}[\int_{t}^{T}D_{t}\varphi (s)dB(s)+\varphi (t)|\mathcal{F}%
_{t}]=\varphi (t),
\end{align}%
since $D_{t}\varphi (s)=0$ for all $s<t$. \hfill $\square $ \bigskip

\subsection{Integration by parts}

\begin{lemma}
\label{lemma 6.3} Suppose $g\in L^2(\mathbb{R})$ and $F \in \mathbb{D}_{1,2}$%
. Then 
\begin{equation}  \label{6.7}
F \diamond \int_\mathbb{R} g(t) dB(t) = F \int_\mathbb{R} g(t) dB(t) - \int_%
\mathbb{R} g(t) D_t F dt .
\end{equation}
\end{lemma}

\noindent {Proof.} \quad To ease the notation, let $\Vert \cdot \Vert =
\Vert \cdot \Vert_{L^2(\mathbb{R})}$ and $(\cdot, \cdot) = (\cdot,
\cdot)_{L^2(\mathbb{R})}$. For $y \in \mathbb{R}$ we define 
\begin{equation*}
G_y := \exp^\diamond \Big\{ y \int_\mathbb{R} g(t) dB(t) \Big\} = \exp\Big\{ %
y \int_\mathbb{R} g(t) dB(t) -\frac{1}{2} y^2 \Vert g \Vert^2 \Big\}.
\end{equation*}
Choose $F = \exp^\diamond \big\{ \int_\mathbb{R} f(t) dB(t) \big\}$ $=\exp%
\big\{\int_\mathbb{R} f(t) dB(t) -\frac{1}{2} \Vert f \Vert^2 \big\}$, where 
$f\in L^2(\mathbb{R})$. Then 
\begin{equation*}
\begin{split}
F \diamond G_y & = \exp^\diamond \Big\{ \int_\mathbb{R} f(t) dB(t) \Big\} %
\diamond \exp^\diamond \Big\{ y \int_\mathbb{R} g(t) dB(t) \Big\} \\
&= \exp^\diamond \Big\{ \int_\mathbb{R}\big[ f(t) + y g(t)\big] dB(t) \Big\}
\\
&= \exp \Big\{ \int_\mathbb{R} \big[f(t) + y g(t)\big] dB(t) - \frac{1}{2}
\Vert f+ y g \Vert^2 \Big\} \\
& = \exp^\diamond \Big\{ \int_\mathbb{R} f(t) dB(t) \Big\} \, \exp^\diamond %
\Big\{ \int_\mathbb{R} y g(t) dB(t) \Big\} \,\exp \Big\{- y (f,g) \Big\} \\
&= F\, G_y \, \exp \Big\{- y (f,g) \Big\}.
\end{split}%
\end{equation*}
Now differentiating with respect to $y$, we get 
\begin{equation}  \label{6.8}
\frac{d}{dy} \big( F \diamond G_{y}\big) = F \diamond \big( G_{y} \diamond
\int_\mathbb{R} g(t) dB(t) \big)
\end{equation}
and 
\begin{equation}  \label{6.9}
\begin{split}
\frac{d}{dy} \big( F \, G_{y} \, &\exp\big\{ -y (f,g)\big\} \big) \\
&= F \, G_{y} \Big[ \int_\mathbb{R} g(t) dB(t)\, \exp\big\{ -y (f,g)\big\} -
(f,g) \exp\big\{ -y (f,g) \big\} \Big].
\end{split}
\end{equation}
Comparing \eqref{6.8} and \eqref{6.9} we get 
\begin{equation*}
F \diamond \big( G_{y} \diamond \int_\mathbb{R} g(t) dB(t) \big) = F\,
G_{y}\exp\big\{ -y (f,g)\big\} \big[ \int_\mathbb{R} g(t) dB(t) - (f,g) \big].
\end{equation*}

In particular, putting $y=0$ we get 
\begin{equation*}
\begin{split}
F \diamond \int_\mathbb{R} g(t) dB(t) &= F \, \int_\mathbb{R} g(t) dB(t) - F
\int_\mathbb{R} f(t) g(t) dt \\
& = F \, \int_\mathbb{R} g(t) dB(t) - \int_\mathbb{R} g(t) D_{t}F dt.
\end{split}%
\end{equation*}
This proves the result if $F = \exp^{\diamond} \big\{ \int_\mathbb{R} f(t)
dB(t)\big\}$ for some $f\in L^{2}(\mathbb{R})$. Since linear combinations of
such $F$'s are dense in $\mathbb{D}_{1,2}$, the result follows by an
approximation argument. 
\hfill $\square$ \bigskip

\vspace{3mm}

\begin{example}
Choose $F = \int_{\mathbb{R}}f(t) dB(t)$ with $f \in L^{2}(\mathbb{R})$.
Then \eqref{6.7} gives 
\begin{equation}  \label{6.10}
\Big( \int_{\mathbb{R}}f(t) dB(t)\Big) \diamond \Big( \int_{\mathbb{R}} g(t)
dB(t)\Big) = \Big( \int_{\mathbb{R}}f(t) dB(t)\Big) \Big( \int_{\mathbb{R}}
g(t) dB(t)\Big) - \int_{\mathbb{R}}f(t) g(t) dt,
\end{equation}
which is in agreement with Theorem \ref{th3.3}.
\end{example}

\begin{remark}
A general formula for the relation between Wick products and ordinary
products can be found in \cite{HO}.
\end{remark}

\begin{theorem}
\label{theorem 6.6} \textbf{(Integration by parts) \cite{DOP}.}\newline
Suppose $F u(t)$, $0\leq t \leq T$, is Skorohod integrable, with $F\in L^2(%
\mathcal{F}_T,P)$.\newline
Then $F \, u(t)$, $0\leq t \leq T$, is Skorohod integrable and 
\begin{equation}  \label{6.12}
\int_{0}^{T} F \, u(t) \delta B(t) = F \int_{0}^{T} u(t) \delta B(t) -
\int_{0}^{T} u(t) D_{t } F dt.
\end{equation}
\end{theorem}

\noindent {Proof.} \quad First assume that $u(t)$, $0\leq t \leq T$, is a
simple function, i.e. it can be written as a finite linear combination of the
form $u(t) = \sum_{i} a_{i} \chi_{(t_{i}, t_{i+1}]}(t)$, $0\leq t \leq T$,
where $a_{i} \in \mathbb{D}_{1,2}$ for all $i$. Then by applying Lemma \ref
{lemma 6.3} twice we get 
\begin{equation*}
\begin{split}
\int_{0}^{T}F\, u(t) \delta B(t) &= \sum_{i}\big( F a_{i}\big) \diamond
\Delta B(t_{i}) \\
&=\sum_{i} F a_{i} \Delta B(t_{i}) - \sum_{i} \int_{t_{i}}^{t_{i}+1} D_{t} %
\big( F a_{i}\big) dt \\
&= F \sum_{i} a_{i} \Delta B(t_{i}) - \sum_{i} \int_{t_{i}}^{t_{i}+1} D_{t} %
\big( F a_{i}\big) dt \\
&= F \Big( \sum_{i} a_{i}\diamond \Delta B(t_{i}) + \sum_{i}
\int_{t_{i}}^{t_{i}+1} D_{t} a_{i} dt \Big) - \sum_{i}
\int_{t_{i}}^{t_{i}+1} D_{t} \big( F a_{i}\big) dt \\
&= F \int_{0}^{T} u(t) \delta B(t) - \int_{0}^{T} u(t) D_{t }F dt.
\end{split}%
\end{equation*}
Now approximate the general $u$ by a sequence $u_{m}$ of simple functions in 
$Dom(\delta) \subseteq L^{2}( P \times\lambda)$ converging to $u$ in $L^{2}(
P \times\lambda)$. We omit the details. \quad 
\hfill $\square$ \bigskip

\subsection{The duality formula}

The following useful result is a consequence of the generalised Clark-Ocone
theorem:

\begin{theorem}
\textbf{{(Generalised duality formula)}.}\newline
Let $F\in L^{2}(\mathcal{F}_{T},P)$ and let $u(t)=u(t,\omega )\in
L^{2}(\lambda \times P)$ be $\mathbb{F}$-adapted. Then $\mathbb{E}[D_{t}F|%
\mathcal{F}_{t}]\in L^{2}(\lambda \times P)$ and 
\begin{equation}
\mathbb{E}\Big[F\int_{0}^{T}u(t)dB(t)\Big]=\mathbb{E}\Big[\int_{0}^{T}u(t)\mathbb{E}[D_{t}F|%
\mathcal{F}_{t}]dt\Big].
\end{equation}
\end{theorem}

\noindent {Proof.} \quad By the generalised Clark-Ocone theorem and the It\^ o isometry we have 
\begin{align*}
\mathbb{E}\Big[F\int_{0}^{T}u(t)dB(t)\Big]& =\mathbb{E}\Big[\Big(\mathbb{E}%
[F]+\int_{0}^{T}\mathbb{E}[D_{t}F|\mathcal{F}_{t}]dB(t)\Big)%
\int_{0}^{T}u(t)dB(t)\Big] \\
& =\mathbb{E}\Big[\int_{0}^{T}u(t)\mathbb{E}[D_{t}F|\mathcal{F}_{t}]dt\Big],
\end{align*}%
\hfill $\square $ \bigskip

\subsection{Connection to the classical Malliavin derivative}


In this section, we recall the basic definition and properties of the
classical Malliavin calculus for Brownian motion. A general reference for
this presentation is the book \cite{DOP}. See also \cite{I}, 
\cite{N} and \cite{S}. \vskip 0.3cm 

A natural starting point is the classical \emph{Wiener-Itô chaos expansion
theorem}, which states that any $F\in L^2(\mathcal{F}_T,P)$ can be written 
\begin{eqnarray}
F=\sum_{n=0}^{\infty}I_n(f_n)
\end{eqnarray}
for a unique sequence of symmetric deterministic functions $f_n\in
L^2(\lambda^n)$, where $\lambda$ is Lebesgue measure on $[0,T]$ and 
\begin{eqnarray}
I_n(f_n)=n!\int^T_0\int^{t_n}_0\cdots\int^{t_2}_0f_n(t_1,%
\cdots,t_n)dB(t_1)dB(t_2)\cdots dB(t_n)
\end{eqnarray}
(the $n$-times iterated integral of $f_n$ with respect to $B(\cdot)$) for $%
n=1,2,\ldots$ and $I_0(f_0)=f_0$ when $f_0$ is a constant.

Moreover, we have the isometry 
\begin{equation}
\mathbb{E}[F^2]=||F||^2_{L^2(P)}=\sum^\infty_{n=0}n!||f_n||^2_{L^2(%
\lambda^n)}.
\end{equation}

\begin{definition}[Classical Malliavin derivative $\widetilde{D}_t$ with
respect to $B(\cdot)$]
\hfill\break \textrm{Let $\mathbb{D}^{(B)}_{1,2}=\mathbb{D}_{1,2}$ be the
space of all $F\in L^2({\mathcal{F}}_T,P)$ such that its chaos expansion
(2.1) satisfies 
\begin{eqnarray}
||F||^2_{\mathbb{D}^{(B)}_{1,2}}:=\sum^\infty_{n=1}n
n!||f_n||^2_{L^2(\lambda^n)}<\infty.
\end{eqnarray}
}

\textrm{For $F\in \mathbb{D}^{(B)}_{1,2}$ and $t\in [0,T]$, we define the
classical \emph{Malliavin derivative} $\widetilde{D}_tF$ of $F$ at $t$ (with
respect to $B(\cdot)$), by 
\begin{eqnarray}  \label{eq2.5a}
\widetilde{D}_tF=\sum^\infty_{n=1}nI_{n-1}(f_n(\cdot,t)),
\end{eqnarray}
where the notation $I_{n-1}(f_n(\cdot,t))$ means that we apply the $(n-1)$%
-times iterated integral to the first $n-1$ variables $t_1,\cdots, t_{n-1}$
of $f_n(t_1,t_2,\cdots,t_n)$ and keep the last variable $t_n=t$ as a
parameter.}
\end{definition}

Using the classical It\^o isometry repeatedly, we can prove the following
important isometry: 
\begin{eqnarray}  \label{isometry}
\mathbb{E}\Big[\int^T_0(\widetilde{D}_tF)^2dt\Big]=\sum^\infty_{n=1}n
n!||f_n||^2_{L^2(\lambda^n)}=:||F||^2_{\mathbb{D}^{(B)}_{1,2}}.
\end{eqnarray}
In particular, this proves that the function $(t,\omega)\rightarrow
D_tF(\omega)$ belongs to $L^2(\lambda \times P)$ if $F \in \mathbb{D}_{1,2}$%
. 

\begin{example}
If $F=\int^T_0f(t)dB(t)$ with $f\in L^2(\lambda)$ deterministic, then 
\begin{equation*}
\widetilde{D}_t F=f(t) \mbox{ for } a.a. \,t\in[0,T].
\end{equation*}
\end{example}

We now proceed to compare the classical Malliavin derivative $\widetilde{D}%
_t $ with the Hida-Malliavin derivative $D_t$:\newline

The following crucial connection between Hermite polynomials $h_n(x)$ and
iterated Wiener-It\^o integrals $I_n$ was proved by It\^o \cite{I}:

\begin{theorem}
Let $g\in L^{2}(\lambda )$ be deterministic. Then 
\begin{equation}
I_{n}(g^{\otimes n})=||g||^{n}h_{n}\left( \frac{\textstyle{%
\int_{0}^{T}g(t)dB(t)}}{||g||}\right) ,
\end{equation}%
where $||g||=||g||_{L^{2}(\lambda )}$ and $\otimes $ denotes the tensor
product, i.e. 
\begin{equation}
g^{\otimes n}(t_{1},t_{2},...,t_{n}):=g(t_{1})g(t_{2})...g(t_{n}).
\end{equation}
\end{theorem}

Combining this result with the chain rule for the Hida-Malliavin derivative
and the properties of the Hermite polynomials, we obtain that if $f_n =
e_k^{\otimes n} \in \widetilde L^2(\mathbb{R}^n)$ then 
\begin{align*}  \label{6.5}
D_t \big( I_n(f_n)\big) &= D_t(h_n(\theta_k)) = h^{\prime }_n (\theta_k)
e_k(t) =n h_{n-1}(\theta_k) e_k(t) \\
&= n I_{n-1}\big( f_n(\cdot,t)\big).= \widetilde{D}_t \big(I_n(f_n)) \big).
\end{align*}

Since any symmetric function $f$ on $[0,T]^n$ can be written as a linear
combination of tensor products of $e_k ^{\prime }s$ we have proved the
following:

\begin{theorem}
Let $F \in \mathbb{D}_{1,2}$. Then the classical Malliavin derivative $%
\widetilde{D}_tF$ of $F$ coincides with the Hida-Malliavin derivative $D_tF$
of $F$. 
\end{theorem}

We conclude that the Hida-Malliavin derivative defined above on $L^2(P)$ is
an extension of the Malliavin derivative defined on the space $\mathbb{D}%
_{1,2}$. Therefore we can from now on without ambiguity use the notation $D_t$
both for the Hida-Malliavin derivative and the classical Malliavin
derivative. \medskip

\section{ White noise theory for L\' evy processes and Poisson random
measures}

The construction we did in Section 2 of the white noise probability space
for Brownian motion can be modified to apply to other processes. For
example, we obtain a white noise theory for L\' evy processes if we proceed
as follows (see \cite{DOP} for details):

\begin{definition}
Let $\nu $ be a measure on $\mathbb{R}_{0}$ such that 
\begin{equation}
\int_{\mathbb{R}}\zeta ^{2}\nu (d\zeta )<\infty .
\end{equation}%
Define 
\begin{equation}
h(\varphi )=\exp (\int_{\mathbb{R}}\Psi (\varphi (x))dx);\quad \varphi \in (%
\mathcal{S}),
\end{equation}%
where 
\begin{equation}
\Psi (w)=\int_{\mathbb{R}}(e^{iw\zeta }-1-iw\zeta )\nu (d\zeta );\quad w\in 
\mathbb{R},\quad i=\sqrt{-1}.
\end{equation}%
Then $h$ satisfies the conditions (i) - (iii) of the Bochner - Minlos-
Sazonov theorem of Section 2. Therefore there exists a probability measure $%
Q $ on $\Omega =\mathcal{S}^{\prime }(\mathbb{R})$ such that 
\begin{equation}
\mathbb{E}_{Q}[e^{i\left\langle \omega ,\varphi \right\rangle
}]:=\int_{\Omega }e^{i\langle \omega ,\varphi \rangle}dQ(\omega )=h(\varphi );\quad
\varphi \in (\mathcal{S}).
\end{equation}%
The triple $(\Omega ,\mathcal{F},Q)$ is called the (pure jump) L\' evy
white noise probability space.
\end{definition}

One can now easily verify the following

\begin{itemize}
\item $\mathbb{E}_{Q}[\left\langle \cdot ,\varphi \right\rangle ]=0;\quad
\varphi \in (\mathcal{S})$

\item $\mathbb{E}_{Q}[\left\langle \cdot ,\varphi \right\rangle ^{2}]=K\int_{%
\mathbb{R}}\varphi ^{2}(y)dy;\quad \varphi \in (\mathcal{S})$,\newline
where $K=\int_{\mathbb{R}}\zeta ^{2}\nu (d\zeta ).$
\end{itemize}

As in Section 2 we use an approximation argument to define 
\begin{equation}
\widetilde{\eta }(t)=\widetilde{\eta }(t,\omega )=\left\langle \omega ,\chi
_{\lbrack 0,t]}\right\rangle ;\quad a.a.(t,\omega )\in \lbrack 0,\infty
)\times \Omega .
\end{equation}%
Then the following holds:

\begin{theorem} \label{th5.2}
The stochastic process $\widetilde{\eta}(t) $ has a c\`adl\`ag version. This
version $\eta(t)$ is a pure jump L\' evy process with L\' evy measure $\nu$.
\end{theorem}

We can now proceed as in Section 2 and develop a white noise theory and
Hida-Malliavin calculus with respect to the process $\eta$. We refer to \cite%
{DOP} for more information.

\section{Exercises}

\begin{enumerate}
\item Prove that 
\begin{equation*}
\tilde{B}(t):=\left\langle \omega ,\chi _{\lbrack 0,t]}\right\rangle ,\text{ 
}t\in 
\mathbb{R}
,
\end{equation*}%
where%
\begin{equation*}
\chi _{\lbrack 0,t]}=\left\{ 
\begin{array}{ll}
1 & \text{if }s\in \lbrack 0,t)\text{ (or }s\in \lbrack t,0)\text{, if }t<0%
\text{),} \\ 
0 & \text{otherwise},
\end{array}
\right.
\end{equation*}
\ has the following properties:\newline
1a) $\mathbb{E}[\tilde{B}(t)]=0;\quad t\geq 0$,\newline
1b) $\mathbb{E}[\tilde{B}(t)^{2}]=t;\quad t\geq 0$,\newline
1c) $\mathbb{E}[\Tilde{B}(t)^{4}]=3t^{2};\quad t\geq 0.$

\item Prove that $L^{2}(\mathcal{F}_{T},P)\subseteq (\mathcal{S})^{\ast }.$

\item Prove that the functional 
\begin{equation}
g(\phi )=e^{-\frac{1}{2}||\phi ||^{2}};\quad \phi \in L^{2}(\mathbb{R}),
\end{equation}%
is positive definite. (See \eqref{posdef}.)


\item Define $F(t)=\mathbb{E}[\exp (\tint_{0}^{T}f(s)dB(s))|\mathcal{F}%
_{t}], $ with $f\in L^{2}([0,T])$ deterministic. Prove that 
\begin{equation*}
D_{t}F=Ff(t);\quad t\in \lbrack 0,T]
\end{equation*}

\item Let $G=\int_{0}^{T}\int_{\mathbb{R}_{0}}\gamma (s,\zeta )\tilde{N}%
(ds,d\zeta )$ and $\varphi (G)=\exp (\int_{0}^{T}\int_{\mathbb{R}_{0}}\gamma
(s,\zeta )\tilde{N}(ds,d\zeta ))$. Prove that 
\begin{equation*}
D_{t,\zeta }\varphi (G)=\exp (G)[\exp (\gamma (t,\zeta ))-1].
\end{equation*}

\item Put $\varphi (G)=\exp (\tint_{0}^{T}\tint_{\mathbb{R}_{0}}\ln
(1+\gamma (s,\zeta ))\tilde{N}(ds,d\zeta ))$. Prove that 
\begin{equation}
D_{t,z}\varphi (G)=\varphi (G)\gamma (t,z).  \label{iv}
\end{equation}

\item Write down the expansion%
\begin{equation*}
F=\underset{\alpha \in J}{\sum }c_{\alpha }H_{\alpha }
\end{equation*}%
for the following random variables and use Definition \ref{HM} to find $%
D_{t}F$:

\begin{description}
\item[7a.] $F=B(t_{0})$ for some $t_{0}\in \lbrack 0,T].$

\item[7b.] $F=\tint_{0}^{T}f(s)dB(s)$ for some deterministic $f\in
L^{2}(\lambda ).$

\item[7c.] $F=\overset{\bullet }{B}({t_{0}})$ for some $t_{0}\in \lbrack
0,T].$
\end{description}
\end{enumerate}


\section{Applications}

\subsection{Backward stochastic differential equations}

\noindent Backward sde's (bsde's) were first introduced in their linear form
by Bismut \cite{b} in connection with a stochastic version of the Pontryagin
maximum principle. Subsequently, this theory was extended by Pardoux and
Peng \cite{PP} to the nonlinear case. The first work applying bsde to
finance was the paper by El Karoui \textit{et al} \cite{EPQ} where they
studied several applications to option pricing and recursive utilities. All
the above mentioned works are in the Brownian motion framework (continuous
case). The discontinuous case is more involved. 
Tang and Li \cite{TL} proved an existence and uniqueness result in the case
of a natural filtration associated with a Brownian motion and a Poisson
random measure. Barles \textit{et al} \cite{BBP} proved a comparison theorem
for such equations and later Royer \cite{R} extended comparison theorem
under weaker assumptions. We define the following spaces for the solution to
live:

\begin{itemize}
\item $S^{2}$ consists of the $\mathbb{F}$-adapted c\`adl\`ag processes $%
Y:\Omega\times\lbrack0,T]\rightarrow\mathbb{R},$ equipped with the norm 
\begin{equation*}
\parallel Y\parallel_{S^{2}}^{2}:=\mathbb{E[}\sup_{t\in%
\lbrack0,T]}|Y(t)|^{2}]<\infty.
\end{equation*}

\item $L^{2}$ consists of the $\mathbb{F}$-predictable processes $%
Z:\Omega\times\lbrack0,T]\rightarrow\mathbb{R},$ with%
\begin{equation*}
\parallel Z\parallel_{L^{2}}^{2}:=\mathbb{E}\Big[\tint _{0}^{T}\left\vert
Z(t)\right\vert ^{2}dt\Big]<\infty.
\end{equation*}
{}

\item $L_{\nu}^{2}$ consists of Borel functions $K:%
\mathbb{R}
_{0}\rightarrow\mathbb{R},$ such that 
\begin{equation*}
\parallel R\parallel_{L_{\nu}^{2}}^{2}:=\tint _{\mathbb{R}%
_{0}}|R(\zeta)|^{2}\nu(d\zeta)<\infty.
\end{equation*}

\item $H_{\nu}^{2}$ consists of $\mathbb{F}$-predictable processes$\
R:\Omega \times\lbrack0,T]\times%
\mathbb{R}
_{0}\rightarrow\mathbb{R},$ such that for any fixed $t\in\lbrack0,T]$, $%
R(t,\zeta)$ is any element in $L_{\nu}^{2}$ and 
\begin{equation*}
\parallel R\parallel_{H_{\nu}^{2}}^{2}:=\mathbb{E}\Big[\tint _{0}^{T}\tint _{%
\mathbb{R}_{0}}R(t,\zeta)^{2}\nu(d\zeta)dt\Big]<\infty.
\end{equation*}
\end{itemize}

\noindent Let us develop the basic idea which is behind BSDEs. For this let
us consider a SDE driven by the Brownian motion $B$, which is of the
following form 
\begin{equation*}
dY(t)=f(Y(t))dt+\sigma (Y(t))dB(t),\,t\in \lbrack 0,T].
\end{equation*}%
Here we suppose for simplicity that the coefficients $f,\sigma :%
\mathbb{R}
\rightarrow 
\mathbb{R}
$ are Lipschitz. Then it is a classical result that, if $Y(0)=\xi \in 
\mathbb{R}
$ is an imposed initial condition, the SDE, also called forward SDE, has a
unique continuous $\mathbb{F}$-adapted solution $Y$. But what does happen,
if we replace now the initial condition by a terminal one? As long as the
condition $\xi $ remains still a deterministic real value, we make a time
inversion and write the equation with terminal condition for $V(t):=Y(T-t),\,%
\tilde{B}(t):=B(T)-B(T-t),\,t\in \lbrack 0,T],$ as an SDE with initial
condition: 
\begin{equation*}
\left\{ 
\begin{array}{ll}
dV(t) & =-f(V(t))dt-\sigma (V(t))d\tilde{B}(t),\,t\in \lbrack 0,T], \\ 
V(0) & =\xi .%
\end{array}%
\,\right.
\end{equation*}%
We observe that $\,\tilde{B}=(\,\tilde{B}(t))_{t\in \lbrack 0,T]}$ is a
Brownian motion, and due to the Lipschitz condition on the coefficients,
there is a unique continuous solution $V$ adapted with respect to the
filtration generated by $\,\tilde{B}$. But this means, that $Y(t)=V(T-t)$ is 
$\tilde{\mathcal{F}}_{t}=\sigma \{B(T)-B(s),\,s\in \lbrack t,T]\}$%
-measurable, for all $t\in \lbrack 0,T]$ (The $\sigma $-fields are
considered as completed), and, for $t_{i}^{n}=T-t+t\frac{i}{n},\ 0\leq i\leq
n,$ the stochastic integral of the SDE for the process $Y$ can be described
by%
\begin{equation*}
\begin{array}{ll}
\int_{0}^{t}\sigma (Y(s))dB(s) & :=\int_{T-t}^{T}\sigma (V(s))d\tilde{B}(s)
\\ 
& =L^{2}-\lim_{n\rightarrow +\infty }\sum_{i=0}^{n-1}\sigma \left(
V(t_{i}^{n})\right) (\tilde{B}(t_{i+1}^{n})-\tilde{B}(t_{i}^{n})) \\ 
& =L^{2}-\lim_{n\rightarrow +\infty }\sum_{i=0}^{n-1}\sigma \left(
Y(T-t_{i}^{n})\right) (B(T-t_{i}^{n})-B(T-t_{i+1}^{n})) \\ 
& =L^{2}-\lim_{n\rightarrow +\infty }\sum_{i=0}^{n-1}\sigma \left( Y(t-t%
\frac{i}{n})\right) (B(t-t\frac{i}{n})-B(t-t\frac{i+1}{n})) \\ 
& =L^{2}-\lim_{n\rightarrow +\infty }\sum_{i=1}^{n}\sigma \left( Y(t\frac{i+1%
}{n})\right) (B(t\frac{i+1}{n})-B(t\frac{i}{n})),%
\end{array}%
\end{equation*}

\noindent i.e., we have to do with the so-called Itô backward integral.

\noindent But how about a terminal condition $Y(T)=\xi$ with $\xi\in L^{2}(%
\mathcal{F}_{T},P)$? We see that in this case a time inversion of our SDE
for $Y$ leads to a forward equation for $V(t)=Y(T-t)$ with an anticipating
initial condition $V(0)=\xi.$ But we are not interested in studying SDEs
with anticipation.

\noindent In order to understand better what to do, let us first consider
the special case where $\sigma \equiv 0$ and $f(y)=ay,$ with $a\in 
\mathbb{R}
$ is a real constant, i.e., we have the SDE 
\begin{equation*}
\left\{ 
\begin{array}{ll}
dY(t) & =aY(t)dt,\,t\in \lbrack 0,T], \\ 
Y(T) & =\xi \in L^{2}(\mathcal{F}_{T},P).%
\end{array}%
\right.
\end{equation*}%
The unique solution of this equation is the process $Y(t)=\xi \exp
\{-a(T-t)\},\,t\in \lbrack 0,T].$ This process is, obviously, not adapted to
the Brownian filtration $\mathbb{F}$. Being interested in adapted solutions
we replace this process $Y(t),\,t\in \lbrack 0,T]$, by its best
approximation in $L^{2}$ by a process which is adapted with respect to the
filtration $\mathbb{F}$, i.e., we consider the process $U(t)=\exp \{-a(T-t)\}%
\mathbb{E}[\xi |\mathcal{F}_{t}],\ t\in \lbrack 0,T],$ which is just the
optional projection of the process $Y$.

\noindent In the special case, where $a=0$, we see that $U(t)=\mathbb{E}[\xi
|\mathcal{F}_{t}],\ t\in \lbrack 0,T],$ is just the martingale generated by $%
U(T)=Y(T)=\xi $, and from the martingale representation property in a
Brownian setting we have the existence and the uniqueness of a square
integrable $\mathbb{F}$-adapted process $Z$ such that 
\begin{equation*}
\xi =\mathbb{E}[\xi ]+\int_{0}^{T}Z(t)dB(t),\,P\text{-a.s.},
\end{equation*}%
i.e., 
\begin{equation*}
U(t)=\mathbb{E}[\xi |\mathcal{F}_{t}]=\mathbb{E}[\xi
]+\int_{0}^{t}Z(s)dB(s),\,t\in \lbrack 0,T],
\end{equation*}%
which shows that $dU(t)=Z(t)dB(t),\,t\in \lbrack 0,T],\,U(T)=\xi .$ This
indicates that we have to reinterpret our equation for $Y$ in the following
sense:\newline
\begin{equation*}
\left\{ 
\begin{array}{ll}
dY(t) & =(f(Y(t))dt+\sigma (Y(t))dB(t))+V(t)dB(t)\,t\in \lbrack 0,T], \\ 
Y(T) & =\xi (\in L^{2}(\mathcal{F}_{T},P)),%
\end{array}%
\right. \,
\end{equation*}%
where the solution we have to look for is a couple of square integrable $%
\mathbb{F}$-adapted processes $(Y,V),$ where $V$ has its origin from the
martingale representation property. However, having $(Y,V)$, we can define
the process $Z=(Z(t))_{t\in \lbrack 0,T]}$ by putting $Z(t):=V(t)+\sigma
(Y(t)),\,t\in \lbrack 0,T]$ (observe that, knowing $(Y,Z)$ we can compute $%
V(t)=Z(t)-\sigma (Y(t))$), and this leads to the SDE, 
\begin{equation*}
\left\{ 
\begin{array}{ll}
dY(t) & =f(Y(t))dt+Z(t)dB(t),t\in \lbrack 0,T],\  \\ 
\,Y(T) & =\xi .%
\end{array}%
\right.
\end{equation*}%
Finally, in order to have the backward SDE in the general form studied by
Pardoux and Peng in 1990, we replace the Lipschitz function $f:\mathbb{R}%
\rightarrow \mathbb{R}$ by a more general, adapted coefficient $f:\Omega
\times \lbrack 0,T]\times \mathbb{R}\times \mathbb{R}\rightarrow \mathbb{R}$
which is now allowed to depend also on $(\omega ,t,y)\in \Omega \times
\lbrack 0,T]\times \mathbb{R}$, but also on the solution component $Z(t)$,
and in order to be coherent with the notation which is usually used, we
endow the function $f$ with the negative sign. This leads to the so-called
backward stochastic differential equation (BSDE) introduced and studied by
Pardoux and Peng in their pioneering work of 1990 (see \cite{PP}): 
\begin{equation*}
\left\{ 
\begin{array}{ll}
dY(t) & =-f(t,Y(t),Z(t))dt+Z(t)dB(t),\,t\in \lbrack 0,T], \\ 
Y(T) & =\xi \in L^{2}(\mathcal{F}_{T},P).%
\end{array}%
\right. \,
\end{equation*}

\begin{theorem}
{\cite{PP}} Let $\xi\in L^{2}(\mathcal{F}_{T},P)$ and $f:\Omega\times\lbrack
0,T]\times\mathbb{R}\times\mathbb{R}\rightarrow\mathbb{R}$ a jointly
measurable mapping satisfying the following assumptions:

i) $f(\cdot,\cdot,0,0)\in L^{2},$

ii) $f(\omega,t,\cdot,\cdot)$ is uniformly Lipschitz, $dtP(d\omega)$-a.e.,
i.e., there is some constant $C\in\mathbb{R}$ such that, $dtP(d\omega)$%
-a.e., for all $y,y^{\prime},z,z^{\prime}\in\mathbb{R},$

\centerline{$|f(t,\omega,y,z)-f(t,\omega,y',z')|\le C(|y-y'|+|z-z'|).$}%
\smallskip

\noindent Then the BSDE 
\begin{equation*}
\left\{ 
\begin{array}{ll}
dY(t) & =-f(t,Y(t),Z(t))dt+Z(t)dB(t),\,t\in \lbrack 0,T], \\ 
Y(T) & =\xi .%
\end{array}%
\right. \,
\end{equation*}%
possesses a unique solution $(Y,Z)\in S^{2}\times L^{2}.$
\end{theorem}

\begin{remark}
So far we have been dealing with BSDEs driven by Brownian motion $B(\cdot)$ only. We now turn to the more general case, with BSDEs driven by both Brownian motion and an independent  compensated Poissson random measure $\widetilde{N}$. As explained in Theorem \ref{th5.2} we can construct both $B$ and $\widetilde{N}$ on the same space $\Omega=\mathcal{S}'(\mathbb{R})$. We refer to \cite{OS1} for  more information on SDEs and BSDEs driven by Brownian motion and Poisson random measure and optimal control of such equations.
\end{remark}

\subsubsection{Representation of solutions of BSDE}

The following result is new:

\begin{theorem} \label{th7.3}
\label{Th2.9} Suppose that $f,p,q$ and $r$ are given càdlàg adapted
processes in $L^{2}(\lambda\times P),L^{2}(\lambda\times
P),L^{2}(\lambda\times P)$ and $L^{2}(\lambda\times\nu\times P)$
respectively, and they satisfy a BSDE of the form 
\begin{equation}
\begin{cases}
dp(t) & =f(t)dt+q(t)dB(t)+\tint _{\mathbb{R}_{0}}r(t,\zeta)\tilde{N}%
(dt,d\zeta);0\leq t\leq T, \\ 
p(T) & =F\in L^{2}(\mathcal{F}_{T},P).%
\end{cases}
\label{BSDE}
\end{equation}

Then for a.a. $t$ and $\zeta$ the following holds: 
\begin{equation}
q(t)=D_{t}p(t^{+}):=\underset{\varepsilon\rightarrow0^{+}}{\lim}%
D_{t}p(t+\varepsilon)\text{ (limit in }(\mathcal{S})^{\ast}),  \label{eq2.24}
\end{equation}%
\begin{equation}
q(t)=\mathbb{E}[D_{t}p(t^{+})|\mathcal{F}_{t}]:=\underset{\varepsilon
\rightarrow0^{+}}{\lim}\mathbb{E}[D_{t}p(t+\varepsilon)|\mathcal{F}_{t}]%
\text{ (limit in }L^{2}(P)),  \label{eq2.25}
\end{equation}

and 
\begin{equation}
r(t,\zeta)=D_{t^{,}\zeta}p(t^{+}):=\underset{\varepsilon\rightarrow0^{+}}{%
\lim}D_{t,\zeta}p(t+\varepsilon)\text{ (limit in }(\mathcal{S})^{\ast}),
\label{eq2.26}
\end{equation}%
\begin{equation}
r(t,\zeta)=\mathbb{E}[D_{t^{,}\zeta}p(t^{+})|\mathcal{F}_{t}]:=\underset{%
\varepsilon\rightarrow0^{+}}{\lim}\mathbb{E}[D_{t,\zeta}p(t+\varepsilon )|%
\mathcal{F}_{t}]\text{ (limit in }L^{2}(P)).  \label{eq2.27}
\end{equation}
\end{theorem}

\noindent {Proof.} \quad Using white noise calculus and the Wick product
representation of the stochastic integrals, we can write the BSDE (\ref{BSDE}%
) as a forward SDE 
\begin{equation}
p(t)=p_{0}+\tint _{0}^{t}f(u)du+\tint _{0}^{t}q(u)\diamond\dot{B}(u)du+\tint
_{0}^{t}\tint _{\mathbb{R}_{0}}r(u,\zeta)\diamond\dot{\tilde{N}}(u,d\zeta)du.
\end{equation}
for some initial value $p(0)=p_{0}$ (constant).\newline
This implies that for all $s<t$ we have, as equations in $(\mathcal{S}%
)^{\ast}$, 
\begin{equation*}
D_{t}p(t+\varepsilon)=\tint _{t}^{t+\varepsilon}D_{t}f(u)du+\tint
_{t}^{t+\varepsilon}D_{t}q(u)\diamond\dot{B}(u)du+q(t)
\end{equation*}
and 
\begin{equation*}
D_{t,\zeta}p(t+\varepsilon)=\tint
_{t}^{t+\varepsilon}D_{t,\zeta}f(u)du+\tint _{t}^{t+\varepsilon}\tint _{%
\mathbb{R}_{0}}D_{t,\zeta}r(u,\zeta)\diamond\dot{\tilde{N}}%
(u,\zeta)du+r(t,\zeta)
\end{equation*}

\noindent Taking the limit in $(\mathcal{S})^{\ast}$ as $\varepsilon
\rightarrow0^{+}$ we get \eqref{eq2.24} and \eqref{eq2.26}.\newline
Taking the conditional expectation and then the limit as $\varepsilon$ goes
to $0^{+}$, we get \eqref{eq2.25} and \eqref{eq2.27}. \hfill $\square$
\bigskip 

\subsubsection{Closed formula for mean-field BSDE}

We shall find the closed formula corresponding to the linear mean-field BSDE
of the form%
\begin{equation}
\left\{ 
\begin{array}{ll}
dY(t) & =-[\alpha _{1}(t)Y(t)+\beta _{1}(t)Z(t)+\tint_{\mathbb{R}_{0}}\eta
_{1}(t,\zeta )K(t,\zeta )\nu (d\zeta )+\alpha _{2}(t)\mathbb{E[}Y(t)] \\ 
& +\beta _{2}(t)\mathbb{E[}Z(t)]+\tint_{\mathbb{R}_{0}}\eta _{2}(t,\zeta )%
\mathbb{E[}K(t,\zeta )]\nu (d\zeta )+\gamma (t)]dt \\ 
& +Z(t)dB(t)+\tint_{\mathbb{R}_{0}}K(t,\zeta )\tilde{N}(dt,d\zeta ),t\in %
\left[ 0,T\right] , \\ 
Y(T) & =\xi ,%
\end{array}%
\right.  \label{lbsde}
\end{equation}%
where the coefficients $\alpha _{1}(t),\alpha _{2}(t),\beta _{1}(t),\beta
_{2}(t),\eta _{1}(t,\cdot ),\eta _{2}(t,\cdot )$ are given deterministic
functions; $\gamma (t)$ is a given $\mathbb{F}$-adapted process and $\xi \in
L^{2}\left( \Omega ,\mathcal{F}_{T}\right) $ is a given $\mathcal{F}_{T}$
measurable random variable. Applying the closed formula for linear BSDE with
jums, the above linear mean-field bsde $(\ref{lbsde})$ can be written as
follows. 
\begin{equation}
\begin{array}{c}
Y(t)=\mathbb{E}[(\xi \Gamma (t,T)+\tint_{t}^{T}\Gamma (t,s)\{\alpha _{2}(s)%
\mathbb{E[}Y(s)]+\beta _{2}(s)\mathbb{E[}Z(s)] \\ 
+\tint_{\mathbb{R}_{0}}\eta _{2}(s,\zeta )\mathbb{E[}K(s,\zeta )]\nu (d\zeta
)+\gamma (s)\}ds)|\mathcal{F}_{t}],\quad t\in \left[ 0,T\right] \,,%
\end{array}
\label{cfbsde}
\end{equation}%
where $\Gamma (t,s)$ is the solution of the following linear sde 
\begin{equation}
\left\{ 
\begin{array}{ll}
d\Gamma (t,s) & =\Gamma (t,s^{-})[\alpha _{1}(t)dt+\beta _{1}(t)dB(t)+\tint_{%
\mathbb{R}_{0}}\eta _{1}(t,\zeta )\widetilde{N}(dt,d\zeta )],\quad s\in %
\left[ t,T\right] , \\ 
\Gamma (t,t) & =1\,.%
\end{array}%
\right.  \label{gam}
\end{equation}%
%
Since we are in one dimension, Equation $(\ref{gam})$ can be solved
explicitly and the solution is given by 
\begin{equation}
\begin{array}{c}
\Gamma (t,s)=\exp \{\tint_{t}^{s}\beta _{1}(r)dB(r)+\tint_{t}^{s}(\alpha
_{1}(r)-\tfrac{1}{2}(\beta _{1}(r))^{2})dr \\ 
\text{ \ \ \ \ \ \ \ \ \ \ \ \ \ \ \ \ \ }+\tint_{t}^{s}\tint_{\mathbb{R}%
_{0}}(\ln (1+\eta _{1}(r,\zeta ))-\eta _{1}(r,\zeta ))\nu (d\zeta )dr \\ 
\text{ \ \ \ \ }+\tint_{t}^{s}\tint_{\mathbb{R}_{0}}(\ln (1+\eta
_{1}(r,\zeta ))\widetilde{N}(dr,d\zeta )\}.%
\end{array}
\label{gexp}
\end{equation}%
%
Notice that 
\begin{equation}
\mathbb{E}\Gamma (t,s)=\exp \{\tint_{t}^{s}\alpha _{1}(r)dr\}\,.
\end{equation}%
To solve \eqref{cfbsde} we take the expectation on both sides of $(\ref%
{cfbsde})$. Denoting $\overline{Y}(t):=\mathbb{E[}Y(t)],$ $\overline{Z}(t):=%
\mathbb{E[}Z(t)]$, and $\overline{K}(t,\zeta ):=\mathbb{E[}K(t,\zeta )]$, we
obtain%
\begin{equation}
\begin{array}{c}
\overline{Y}(t)=\mathbb{E}[\xi \Gamma (t,T)+\tint_{t}^{T}\Gamma
(t,s)\{\alpha _{2}(s)\overline{Y}(s)+\beta _{2}(s)\overline{Z}(s) \\ 
+\tint_{\mathbb{R}_{0}}\eta _{2}(s,\zeta )\overline{K}(s,\zeta )\nu (d\zeta
)+\gamma (s)\}ds],t\in \left[ 0,T\right] .%
\end{array}
\label{e.meany}
\end{equation}%
To find equations for $\overline{Z}(t)$ and $\overline{K}(t,\zeta )$ we
write the original equation \eqref{lbsde} as a forward one: 
\begin{equation*}
\begin{array}{c}
Y(t)=Y(0)+\tint_{0}^{t}[\alpha _{1}(s)Y(s)+\alpha _{2}(s)\overline{Y}%
(s)+\beta _{1}(s)Z(s)+\beta _{2}(s)\overline{Z}(s) \\ 
+\tint_{\mathbb{R}_{0}}(\eta _{1}(s,\zeta )K(s,\zeta )+\eta _{2}(s,\zeta )%
\overline{K}(s,\zeta ))\nu (d\zeta )+\gamma (s)]ds \\ 
+\tint_{0}^{t}Z(s)dB(s)+\tint_{0}^{t}\tint_{\mathbb{R}_{0}}K(s,\zeta )\tilde{%
N}(ds,d\zeta ),\quad t\in \left[ 0,T\right] ,%
\end{array}%
\end{equation*}%
for some deterministic initial value $Y(0)$. Then using the fundamental
theorem of stochastic calculus (Theorem \ref{fundamental}), we compute the
Hida-Malliavin derivative of $Y(t)$ for all $r<t$ as follows: 
\begin{align*}
D_{r}Y(t)& =\tint_{r}^{t}D_{r}[\alpha _{1}(s)Y(s)+\alpha _{2}(s)\overline{Y}%
(s)+\beta _{1}(s)Z(s)+\beta _{2}(s)\overline{Z}(s) \\
& +\tint_{\mathbb{R}_{0}}(\eta _{1}(s,\zeta )K(s,\zeta )+\eta _{2}(s,\zeta )%
\overline{K}(s,\zeta ))\nu (d\zeta )+\gamma (s)]ds \\
& +{\int_{r}^{t}}D_{r}Z(s)dB(s)+Z(r).
\end{align*}%
Letting $r\rightarrow t-$, we get that $Z(t)=D_{t}Y(t).$ Thus, to find $Z(t)$
we only need to compute $D_{t}Y(t)$. We shall use the expression %
\eqref{cfbsde} for $Y(t)$ and the identity 
\begin{equation*}
D_{t}\mathbb{E}[F|\mathcal{F}_{t}]=\mathbb{E}[D_{t}F|\mathcal{F}_{t}]\,.
\end{equation*}%
We also notice that $D_{t}\Gamma (t,T)=\Gamma (t,T)\beta _{1}(t)$. Then 
\begin{align*}
Z(t)& =\mathbb{E}[(D_{t}\xi \Gamma (t,T)+\xi \Gamma (t,T)\beta
_{1}(t)+\tint_{t}^{T}\Gamma (t,s)\beta _{1}(t)\{\alpha _{2}(s)\overline{Y}(s)
\\
& \qquad +\beta _{2}(s)\overline{Z}(s)+\tint_{\mathbb{R}_{0}}\eta
_{2}(s,\zeta )\overline{K}(s,\zeta )\nu (d\zeta )+\gamma (s)\}ds]\,.
\end{align*}%
Taking the expectation, we have 
\begin{align}
\overline{Z}(t)& =\mathbb{E}[D_{t}\xi \Gamma (t,T)+\beta _{1}(t)\mathbb{E(}%
\xi \Gamma (t,T))+\tint_{t}^{T}\mathbb{E(}\Gamma (t,s))\beta _{1}(t)\{\alpha
_{2}(s)\overline{Y}(s)  \notag \\
& \qquad +\beta _{2}(s)\overline{Z}(s)+\tint_{\mathbb{R}_{0}}\eta
_{2}(s,\zeta )\overline{K}(s,\zeta )\nu (d\zeta )+\gamma (s)\}ds].
\label{e.meanz}
\end{align}%
Similarly, we have $K(t,\zeta )=D_{t,\zeta }Y(t)$ which yields 
\begin{align*}
K(t,\zeta )& =\mathbb{E}[(D_{t,\zeta }\xi \Gamma (t,T)+\xi \Gamma (t,T)\eta
_{1}(t,\zeta )+\tint_{t}^{T}\Gamma (t,s)\eta _{1}(t,\zeta )\{\alpha _{2}(s)%
\overline{Y}(s) \\
& \qquad +\beta _{2}(s)\overline{Z}(s)+\tint_{\mathbb{R}_{0}}\eta
_{2}(s,\zeta )\overline{K}(s,\zeta )\nu (d\zeta )+\gamma (s)\}ds)|\mathcal{F}%
_{t}]\,.
\end{align*}%
Taking the expectation yields 
\begin{eqnarray}
\overline{K}(t,\zeta ) &=&\mathbb{E}[D_{t,\zeta }\xi \Gamma (t,T)+\xi \Gamma
(t,T)\eta _{1}(t,\zeta )+\tint_{t}^{T}\Gamma (t,s)\eta _{1}(t,\zeta
)\{\alpha _{2}(s)\overline{Y}(s)\,  \label{e.meank} \\
&&\qquad +\beta _{2}(s)\overline{Z}(s)+\tint_{\mathbb{R}_{0}}\eta
_{2}(s,\zeta )\overline{K}(s,\zeta )\nu (d\zeta )+\gamma (s)\}ds]\,.  \notag
\end{eqnarray}%
Equations \eqref{e.meany}, \eqref{e.meanz} and \eqref{e.meank} can be used
to obtain $\bar{Y},\bar{Z},\bar{K}$. In fact, we let 
\begin{equation*}
V(t)=\left( 
\begin{array}{c}
V_{1}(t) \\ 
V_{2}(t) \\ 
V_{3}(t,\zeta )%
\end{array}%
\right) =\left( 
\begin{array}{c}
\overline{Y}(t) \\ 
\overline{Z}(t) \\ 
\overline{K}(t,\zeta )%
\end{array}%
\right) \in L^{2}\times L^{2}\times H_{\nu }^{2},
\end{equation*}%
and 
\begin{align}
& 
\begin{array}{c}
A(t,s,\zeta )=\left( A_{ij}(t,s,\zeta )\right) _{1\leq i,j\leq 3}%
\end{array}
\label{A} \\
& 
\begin{array}{c}
=\left( 
\begin{matrix}
\exp \{\int_{t}^{s}\alpha _{1}(r)dr\}\alpha _{2}(s) & \exp
\{\int_{t}^{s}\alpha _{1}(r)dr\}\beta _{2}(s) & \exp \{\int_{t}^{s}\alpha
_{1}(r)dr\}\eta _{2}(s,\zeta ) \\ 
\exp \{\int_{t}^{s}\alpha _{1}(r)dr\}\beta _{1}(t)\alpha _{2}(s) & \exp
\{\int_{t}^{s}\alpha _{1}(r)dr\}\beta _{1}(t)\beta _{2}(s) & \exp
\{\int_{t}^{s}\alpha _{1}(r)dr\}\beta _{1}(t)\eta _{2}(s,\zeta ) \\ 
\exp \{\int_{t}^{s}\alpha _{1}(r)dr\}\eta _{1}(t,\zeta )\alpha _{2}(s) & 
\exp \{\int_{t}^{s}\alpha _{1}(r)dr\}\eta _{1}(t,\zeta )\beta _{2}(s) & \exp
\{\int_{t}^{s}\alpha _{1}(r)dr\}\eta _{1}(t,\zeta )\eta _{2}(s,\zeta )\ 
\end{matrix}%
\right) \,.%
\end{array}
\notag
\end{align}%
Define a mapping $A=A^{T}$ from $V=(V_{1},V_{2},V_{3})^{T}\in L^{2}\times
L^{2}\times H_{\nu }^{2}$ to itself by 
\begin{equation}
(AV)_{i}(t,\zeta
)=\tsum_{j=1}^{2}\tint_{t}^{T}A_{ij}(t,s)V_{j}(s)ds+\tint_{t}^{T}\tint_{%
\mathbb{R}_{0}}A_{i3}(t,s,\zeta )V_{3}(s,\zeta )\nu (d\zeta )\,ds.
\end{equation}%
Then \eqref{e.meany}, \eqref{e.meanz} and \eqref{e.meank} can be written as 
\begin{equation}
V=F+AV\,,  \label{e.v}
\end{equation}%
where 
\begin{equation}
F(t,\zeta )=\left( 
\begin{matrix}
\mathbb{E}(\xi \Gamma (t,T))+\int_{t}^{T}\gamma (s)ds \\ 
\mathbb{E[}D_{t}\xi \Gamma (t,T)+\beta _{1}(t)\xi \Gamma
(t,T)]+\int_{t}^{T}\gamma (s)ds \\ 
\mathbb{E[}D_{t,\zeta }\xi \Gamma (t,T)+\xi \Gamma (t,T)\eta _{1}(t,\zeta
)]+\int_{t}^{T}\gamma (s)ds%
\end{matrix}%
\right) \,.  \label{F}
\end{equation}%
Note that the operator norm of $A$, $||A||$, is less than $1$ if $t$ is
close enough to $T$. Therefore there exists $\delta >0$ such that $||A||<1$
if we restrict the operator to the interval $[T-\delta ,T]$ for some $\delta
>0$ small enough. In this case the linear equation equation \eqref{e.v} can
now be solved easily as follows: 
\begin{equation*}
(I-A)V=F\,,
\end{equation*}%
or 
\begin{equation}
V=(I-A)^{-1}F=\tsum_{n=0}^{\infty }A^{n}F\,;\quad t\in \lbrack T-\delta ,T].
\end{equation}%
Next, using $V(T-\delta )$ as the terminal value of the corresponding BSDE
in the interval $[T-2\delta ,T-\delta ]$ and repeating the argument above,
we find that there exists a solution $V$ of the BSDE in this interval, given
by the equation 
\begin{equation}
V(t,\zeta )=V(T-\delta ,\zeta )+A^{T-\delta }(t,\cdot ,\zeta )V(\cdot
);\quad T-2\delta \leq t\leq T-\delta .
\end{equation}%
Proceeding by induction we end up with a solution on the whole interval $%
[0,T]$. We summarise this as follows:


\begin{theorem}[Closed formula \protect\cite{AHO}]
Assume that $\alpha _{1}(t),\alpha _{2}(t),\beta _{1}(t),\beta _{2}(t),\eta
_{1}(t,\cdot ),\eta _{2}(t,\cdot )$ are given bounded deterministic
functions and that $\gamma (t)$ is $\mathbb{F}$-adapted and $\xi \in
L^{2}\left( \Omega ,\mathcal{F}_{T}\right) $. Then the component $Y(t)$ of
the solution of the linear mean-field BSDE (\ref{lbsde}) can be written on
its closed formula as follows%
\begin{equation}
Y(t)=\mathbb{E}[(\xi \Gamma (t,T)+\tint_{t}^{T}\Gamma (t,s)\{(\alpha
_{2}(s),\beta _{2}(s),\eta _{2}(s,\zeta ))V(s)+\gamma (s)\}ds)|\mathcal{F}%
_{t}],t\in \left[ 0,T\right] ,\text{ }\mathbb{P}\text{-a.s.,}  \label{cf_y}
\end{equation}%
where 
\begin{equation*}
\begin{array}{c}
\Gamma (t,s)=\exp \{\tint_{t}^{s}\beta _{1}(r)dB(r)+\tint_{t}^{s}(\alpha
_{1}(r)-\tfrac{1}{2}(\beta _{1}(r))^{2})dr \\ 
\text{ \ \ \ \ \ \ \ \ \ \ \ \ \ \ \ \ \ }+\tint_{t}^{s}\tint_{\mathbb{R}%
_{0}}(\ln (1+\eta _{1}(r,\zeta ))-\eta _{1}(r,\zeta ))\nu (d\zeta )dr \\ 
\text{ \ \ \ \ }+\tint_{t}^{s}\tint_{\mathbb{R}_{0}}(\ln (1+\eta
_{1}(r,\zeta ))\widetilde{N}(dr,d\zeta )\}.%
\end{array}%
\end{equation*}%
and, inductively, 
\begin{equation}
V(t,\zeta )=V(T-k\delta ,\zeta )+A^{T-k\delta }(t,\cdot ,\zeta )V(\cdot
);\quad T-(k+1)\delta \leq t\leq T-k\delta ;\quad k=0,1,2,...
\end{equation}%
Or, equivalently, 
\begin{equation*}
V(t,\zeta )=(A^{T-k\delta }(t,\cdot ,\zeta ))^{n}V(T-k\delta ,\cdot );\quad
T-(k+1)\delta \leq t\leq T-k\delta ;\quad k=0,1,2,...
\end{equation*}%
where $A^{S};S>0$ is given by (\ref{A}) and $V(T,\zeta )=F$.
\end{theorem}

\subsection{Stochastic maximum principles via Hida-Malliavin calculus}

Recall the two main methods of optimal control of systems described by Itô- L%
évy processes:

\begin{itemize}
\item \emph{Dynamic programming and the Hamilton-Jacobi-Bellman (HJB)
equation}\newline
This method was introduced by R. Bellman in the 1950's, first in the  deterministic case.

\item \emph{The stochastic maximum principle} \newline
This method was established at around the same time by L. Pontryagin and his group in the deterministic case.\newline
The maximum principle was extended to the stochastic case by J.-M. Bismut  for the linear-quadratic case and for Brownian motion driven SDE's (1976)  and subsequently further developed by 
A. Bensoussan, S. Peng, E. Pardoux and others (still for Brownian motion driven SDE's, 1980 -1990), 
and then by 
S. Rong and N.-C. Framstad, A. Sulem \& B. \O ksendal (for jump diffusions, 1990 -).
\end{itemize}

Dynamic programming is efficient when applicable, but it requires that the
system is Markovian. The maximum principle has the advantage that it also
applies to non-Markovian SDE's, but the drawback is the corresponding
complicated BSDE for the adjoint processes.\bigskip

\noindent The state of our system $X^{u}(t)=X(t)$ satisfies the following SDE%
\begin{equation}
\left\{ 
\begin{array}{ll}
dX(t) & =b(t,X(t),u(t))dt+\sigma (t,X(t),u(t))dB(t) \\ 
& +\int_{\mathbb{R}_{0}}\gamma (t,X(t),u(t),\zeta )\tilde{N}(dt,d\zeta
);0\leq t\leq T, \\ 
X(0) & =x_{0}\in \mathbb{R}\ (\text{constant}),%
\end{array}%
\right.  \label{sde-}
\end{equation}%
where $b(t,x,u)=b(t,x,u,\omega ):\left[ 0,T\right] \times \mathbb{R}\times
U\times \Omega \rightarrow \mathbb{R}$, $\sigma (t,x,u)=\sigma (t,x,u,\omega
):\left[ 0,T\right] \times \mathbb{R}\times U\times \Omega \rightarrow 
\mathbb{R}
$ and $\gamma (t,x,u,\zeta )=:\left[ 0,T\right] \times \mathbb{R}\times
U\times 
\mathbb{R}
_{0}\times \Omega \rightarrow 
\mathbb{R}
$.\newline
From now on we fix an open convex set $U$ such that $V\subset U$ and we
assume that $b$, $\sigma $ and $\gamma $ are continuously differentiable and
admits uniformly bounded partial derivatives in $U$ with respect to $x$ and $%
u$.\newline
Moreover, we assume that the coefficients $b$, $\sigma $ and $\gamma $ are $%
\mathbb{F}$-adapted, and uniformly Lipschitz continuous with respect to $x$,
in the sense that there is a constant $C$ such that, for all $t\in \lbrack
0,T],u\in V,\zeta \in \mathbb{R}_{0},\,x,x^{\prime }\in \mathbb{R}$ we have 
\begin{equation*}
\begin{array}{c}
\left\vert b\left( t,x,u\right) -b\left( t,x^{\prime },u\right) \right\vert
^{2}+\left\vert \sigma \left( t,x,u\right) -\sigma \left( t,x^{\prime
},u\right) \right\vert ^{2} \\ 
+\tint_{\mathbb{R}_{0}}\left\vert \gamma \left( t,x,u,\zeta \right) -\gamma
\left( t,x^{\prime },u,\zeta \right) \right\vert ^{2}\nu (d\zeta )\leq
C\left\vert x-x^{\prime }\right\vert ^{2},\text{a.s.}%
\end{array}%
\end{equation*}%
Under this assumption, there is a unique solution $X\in \mathcal{S}^{2}$ to
the equation $\left( \ref{sde-}\right) $, such that%
\begin{equation*}
\begin{array}{c}
X(t)=x_{0}+{\tint_{0}^{t}}b(s,X(s),u(s))ds+{\tint_{0}^{t}}\sigma
(s,X(s),u(s))dB(s) \\ 
+{\tint_{0}^{t}}\tint_{\mathbb{R}_{0}}\gamma (s,X(s),u(s),\zeta )\tilde{N}%
(ds,d\zeta );0\leq t\leq T.%
\end{array}%
\end{equation*}%
\newline
For a given set $\mathcal{A}$ of admissible contyrols, the \emph{performance functional} has the form%
\begin{equation}
\begin{array}{ll}
J(u) & =\mathbb{E[}{\tint_{0}^{T}}\text{ }f(t,X(t),u(t))dt+g(X(T))],\quad
u\in \mathcal{A},%
\end{array}
\label{perf-}
\end{equation}%
with given functions $f:\left[ 0,T\right] \times \mathbb{R}\times U\times
\Omega \rightarrow \mathbb{R}$ and $g:\Omega \times \mathbb{R}\rightarrow 
\mathbb{R},$ assumed to be $\mathbb{F}$-adapted and $\mathcal{F}_{T}$%
-measurable, respectively, and continuously differentiable with respect to $%
x $ and $u$ with bounded partial derivatives in $U$.\newline
Suppose that $\hat{u}$ is an optimal control. Fix $\tau \in \lbrack
0,T),0<\epsilon <T-\tau $ and a bounded $\mathcal{F}_{\tau }$-measurable $v$
and define the spike perturbed $u^{\epsilon }$ of the optimal control $\hat{u%
}$ by

\begin{equation*}
u^{\epsilon}(t)=%
\begin{cases}
\hat{u}(t); & t\in\lbrack0,\tau)\cup(\tau+\epsilon,T], \\ 
v; & t\in\lbrack\tau,\tau+\epsilon].%
\end{cases}%
\end{equation*}
Let $X^{\epsilon}(t):=X^{u^{\epsilon}}(t)$ and $\hat{X}(t):=X^{\hat{u}}(t)$
be the solutions of $\left( \ref{sde-}\right) $\ corresponding to $%
u=u^{\epsilon}$ and $u=\hat{u}$, respectively.\newline
Define

\begin{equation}
Z^{\epsilon}(t):=X^{\epsilon}(t)-\hat{X}(t);\text{ }t\in\lbrack0,T].
\label{eq3.4}
\end{equation}
Then by the \emph{mean value theorem} \footnote{%
Recall that if a function $f$ is continuously differentiable on an open
convex set $U\subset\mathbb{R}^{n}$ and continuous on the closure $\bar{U}$,
then for all $x,y\in\bar{U}$ there exists a point $\tilde{x}$ on the
straight line connecting $x$ and $y$ such that 
\begin{equation}
f(y)-f(x)=f^{\prime}(\tilde{x})(y-x):=\sum_{i=1}^{n}\frac{\partial f}{%
\partial x_{i}}(\tilde{x})(y_{i}-x_{i})
\end{equation}%
}, we can write%
\begin{equation*}
\begin{array}{ll}
b^{\epsilon}(t)-\hat{b}(t) & =\frac{\partial\tilde{b}}{\partial x}%
(t)Z^{\epsilon}(t)+\frac{\partial\tilde{b}}{\partial u}(t)(u^{\epsilon }(t)-%
\hat{u}(t)),%
\end{array}%
\end{equation*}
where 
\begin{equation*}
b^{\epsilon}(t)=b(t,X^{\epsilon}(t),u^{\epsilon}(t)),\hat{b}(t)=b(t,\hat {X}%
(t),\hat{u}(t)),
\end{equation*}
and 
\begin{equation*}
\begin{array}{ll}
\frac{\partial\tilde{b}}{\partial x}(t) & =\frac{\partial b}{\partial x}%
(t,x,u)_{x=\tilde{X}(t),u=\tilde{u}(t)},%
\end{array}%
\end{equation*}
and%
\begin{equation*}
\begin{array}{ll}
\frac{\partial\tilde{b}}{\partial u}(t) & =\frac{\partial b}{\partial u}%
(t,x,u)_{x=\tilde{X}(t),u=\tilde{u}(t)}.%
\end{array}%
\end{equation*}
Here $(\tilde{u}(t),\tilde{X}(t))$ is \emph{a point on the straight line}
between $(\hat{u}(t),\hat{X}(t))$ and $(u^{\epsilon}(t),X^{\epsilon}(t))$.
With a similar notation for $\sigma$ and $\gamma$, we get

\begin{align}
Z^{\epsilon}(t) & ={\tint _{\tau}^{t}}\{\tfrac{\partial\tilde{b}}{\partial x}%
(s)Z^{\epsilon}(s)+\tfrac {\partial\tilde{b}}{\partial u}(s)(u^{\epsilon}(s)-%
\hat{u}(s))\}ds+{\tint _{\tau}^{t}}\{\tfrac{\partial\tilde{\sigma}}{\partial
x}(s)Z^{\epsilon}(s)+\tfrac {\partial\tilde{\sigma}}{\partial u}%
(s)(u^{\epsilon}(s)-\hat{u}(s))\}dB(s)  \notag  \label{eq3.5} \\
& +{\tint _{\tau}^{t}}\tint _{\mathbb{R}_{0}}\{\tfrac{\partial\tilde{\gamma}%
}{\partial x}(s,\zeta)Z^{\epsilon}(s)+\tfrac{\partial\tilde{\gamma}}{%
\partial u}(s,\zeta)(u^{\epsilon}(s)-\hat{u}(s))\}\tilde{N}%
(ds,d\zeta);\tau\leq t\leq\tau+\epsilon,
\end{align}
and

\begin{equation*}
\begin{array}{l}
Z^{\epsilon}(t)={\tint _{\tau+\epsilon}^{t}}\tfrac{\partial\tilde{b}}{%
\partial x}(s)Z^{\epsilon}(s)ds+\tint _{\tau+\epsilon}^{t}\tfrac{\partial%
\tilde{\sigma}}{\partial x}(s)Z^{\epsilon}(s)dB(s) \\ 
\text{ \ \ \ \ \ \ \ \ \ \ \ }+\tint _{\tau+\epsilon}^{t}\tint _{\mathbb{R}%
_{0}}\tfrac{\partial\tilde{\gamma}}{\partial x}(s,\zeta)(s)Z^{\epsilon}(s)%
\tilde {N}(ds,d\zeta);\tau+\epsilon\leq t\leq T.%
\end{array}%
\end{equation*}
On other words,%
\begin{equation}
\left\{ 
\begin{array}{ll}
dZ^{\epsilon}(t) & =\{\tfrac{\partial\tilde{b}}{\partial x}(t)Z^{\epsilon
}(t)+\tfrac{\partial\tilde{b}}{\partial u}(t)(v-\hat{u}(t))\}dt+\{\tfrac {%
\partial\tilde{\sigma}}{\partial x}(t)Z^{\epsilon}(t)+\tfrac{\partial \tilde{%
\sigma}}{\partial u}(t)(v-\hat{u}(t))\}dB(t) \\ 
& +\tint _{\mathbb{R}_{0}}\{\tfrac{\partial\tilde{\gamma}}{\partial x}%
(t,\zeta)Z^{\epsilon}(t)+\tfrac{\partial\tilde{\gamma}}{\partial u}%
(t,\zeta)(v-\hat{u}(t))\}\tilde{N}(dt,d\zeta);\tau\leq t\leq\tau+\epsilon,%
\end{array}
\right.  \label{var1}
\end{equation}
and 
\begin{equation}
\begin{array}{ll}
dZ^{\epsilon}(t) & =\tfrac{\partial\tilde{b}}{\partial x}(t)Z^{\epsilon
}(t)dt+\tfrac{\partial\tilde{\sigma}}{\partial x}(t)Z^{\epsilon}(t)dB(t)+%
\tint _{\mathbb{R}_{0}}\tfrac{\partial\tilde{\gamma}}{\partial x}%
(t,\zeta)Z^{\epsilon}(t)\tilde {N}(dt,d\zeta);\tau+\epsilon\leq t\leq T.%
\end{array}
\label{var2}
\end{equation}

\begin{remark}
\-

\begin{enumerate}
\item Note that since the process 
\begin{equation*}
\eta(t):=\tint _{0}^{t}\tint _{\mathbb{R}_{0}}\zeta\tilde{N}%
(ds,d\zeta);t\geq0
\end{equation*}
is a Lévy process, we know that for every given (deterministic) time $t\geq0$
the probability that $\eta$ jumps at $t$ is $0$. Hence, for each $t$, the
probability that $X$ makes jump at $t$ is also $0$. Therefore we have 
\begin{equation*}
Z^{\epsilon}(\tau)=0\text{ a.s. }
\end{equation*}

\item We remark that the equations $\left( \ref{var1}\right) -\left( \ref%
{var2}\right) $ are linear SDE and then by our assumptions on the
coefficients, they admit a unique solution.\newline
\end{enumerate}
\end{remark}

\noindent Let $\mathcal{R}$ denote the set of (Borel) measurable functions $%
r:\mathbb{R}_{0}\rightarrow\mathbb{R}$ and define the Hamiltonian $H:\left[
0,T\right] \times\mathbb{R}\times U\times\mathbb{R}\times\mathbb{R}\times%
\mathcal{R}\times\Omega\rightarrow\mathbb{R}$, to be

\begin{align}
H(t,x,u,p,q,r) & :=H(t,x,u,p,q,\omega)=f(t,x,u)+b(t,x,u)p  \notag \\
& \quad\quad+\sigma(t,x,u)q+\tint _{\mathbb{R}_{0}}\gamma(t,x,u,\zeta)r(%
\zeta)\nu(d\zeta).  \label{h}
\end{align}
Let $(p^{\epsilon},q^{\epsilon},r^{\epsilon})\in\mathcal{S}^{2}\times
L^{2}\times L_{\nu}^{2}$ be the solution of the following associated adjoint
BSDE: 
\begin{equation}
\begin{cases}
dp^{\epsilon}(t) & =-\tfrac{\partial\tilde{H}}{\partial x}(t)dt+q^{\epsilon
}(t)dB(t)+\tint _{\mathbb{R}_{0}}r^{\epsilon}(t,\zeta)\tilde{N}%
(dt,d\zeta);t\in\lbrack0,T], \\ 
p^{\epsilon}(T) & =\tfrac{\partial\tilde{g}}{\partial x}(\tilde{X}(T)),%
\end{cases}
\label{p-}
\end{equation}
where 
\begin{equation*}
\begin{array}{ll}
\tfrac{\partial\tilde{H}}{\partial x}(t) & =\tfrac{\partial\tilde{f}}{%
\partial x}(t)+\tfrac{\partial\tilde{b}}{\partial x}(t)p^{\epsilon}(t)+%
\tfrac {\partial\tilde{\sigma}}{\partial x}(t)q^{\epsilon}(t)+\tint _{%
\mathbb{R}_{0}}\tfrac{\partial\tilde{\gamma}}{\partial x}(t,\zeta)r^{%
\epsilon}(t,\zeta )\nu(d\zeta).%
\end{array}%
\end{equation*}

\begin{lemma}
{\cite{AO-}} \label{est} The following holds, 
\begin{equation}
Z^{\epsilon}(t)\rightarrow0\text{ as }\epsilon\rightarrow0^{+}\text{;}\quad%
\text{ for all }t\in\lbrack\tau,T].  \label{esz}
\end{equation}%
\begin{equation}
(p^{\epsilon},q^{\epsilon},r^{\epsilon})\rightarrow(\hat{p},\hat{q},\hat {r})%
\text{ when }\epsilon\rightarrow0^{+},  \label{esa}
\end{equation}
where $(\hat{p},\hat{q},\hat{r})$ is the solution of the BSDE 
\begin{equation*}
\begin{cases}
d\hat{p}(t) & =-\tfrac{\partial\hat{H}}{\partial x}(t)dt+\hat{q}%
(t)dB(t)+\tint _{\mathbb{R}_{0}}\hat{r}(t,\zeta)\tilde{N}(dt,d\zeta);t\in%
\lbrack0,T], \\ 
\hat{p}(T) & =\frac{\partial g}{\partial x}(\hat{X}(T)).%
\end{cases}%
\end{equation*}
\end{lemma}

\noindent {Proof.} \quad By the Itô formula, we see that the solutions of
the equations $\left( \ref{var1}\right) -\left( \ref{var2}\right) ,$ are 
\begin{align}
Z^{\epsilon }(t)& =Z^{\epsilon }(\tau +\epsilon )\exp (\tint_{\tau +\epsilon
}^{t}\{\tfrac{\partial \tilde{b}}{\partial x}(s)-\tfrac{1}{2}(\tfrac{%
\partial \tilde{\sigma}}{\partial x}(s))^{2}+\tint_{\mathbb{R}_{0}}[\log (1+%
\tfrac{\partial \tilde{\gamma}}{\partial x}(s,\zeta ))-\tfrac{\partial 
\tilde{\gamma}}{\partial x}(s,\zeta )]\nu (d\zeta )\}ds  \notag \\
& \quad \quad +\tint_{\tau +\epsilon }^{t}\tfrac{\partial \tilde{\sigma}}{%
\partial x}(s)dB(s)+\tint_{\tau +\epsilon }^{t}\tint_{\mathbb{R}_{0}}\log (1+%
\tfrac{\partial \tilde{\gamma}}{\partial x}(s,\zeta ))\tilde{N}(ds,d\zeta
));\quad \tau +\epsilon \leq t\leq T.  \label{z1}
\end{align}%
and%
\begin{equation}
\begin{array}{ll}
Z^{\epsilon }(t) & =\Upsilon (t)^{-1}[\tint_{0}^{t}\Upsilon (s)(\tfrac{%
\partial \tilde{b}}{\partial u}(s)(u^{\epsilon }(s)-\hat{u}(s)) \\ 
& +{\dint_{\mathbb{R}_{0}}}\left( \tfrac{1}{1+\tfrac{\partial \tilde{\gamma}%
}{\partial x}(s,\zeta )}-1\right) \tfrac{\partial \tilde{\gamma}}{\partial u}%
(s,\zeta )(v-\hat{u}(s))\nu (d\zeta ))ds+\tint_{0}^{t}\Upsilon (s)\tfrac{%
\partial \tilde{\sigma}}{\partial u}(s)(v-\hat{u}(s))dB(s) \\ 
& +\tint_{0}^{t}\tint_{\mathbb{R}_{0}}\Upsilon (s)\left( \tfrac{\tfrac{%
\partial \tilde{\gamma}}{\partial u}(s,\zeta )(v-\hat{u}(s))}{1+\tfrac{%
\partial \tilde{\gamma}}{\partial x}(s,\zeta )}-1\right) \tilde{N}(ds,d\zeta
)];\tau \leq t\leq \tau +\epsilon ,%
\end{array}
\label{z2}
\end{equation}%
where%
\begin{equation*}
\left\{ 
\begin{array}{ll}
d\Upsilon (t) & =\Upsilon (t^{-})\left[ -\tfrac{\partial \tilde{b}}{\partial
x}(t)+(\tfrac{\partial \tilde{\sigma}}{\partial x}(t)(u^{\epsilon }(t)-\hat{u%
}(t)))^{2}\right. \\ 
& +\dint_{\mathbb{R}_{0}}\left\{ \tfrac{1}{1+\tfrac{\partial \tilde{\gamma}}{%
\partial x}(t,\zeta )}-1+\tfrac{\partial \tilde{\gamma}}{\partial x}(t,\zeta
)\right\} \nu (d\zeta )dt-\tfrac{\partial \tilde{\sigma}}{\partial x}(t)dB(t)
\\ 
& \left. +\dint_{\mathbb{R}_{0}}\left( \tfrac{1}{1+\tfrac{\partial \tilde{%
\gamma}}{\partial x}(t,\zeta )}-1\right) \tilde{N}(dt,d\zeta )\right] ;\tau
\leq t\leq \tau +\epsilon , \\ 
\Upsilon (0) & =1.%
\end{array}%
\right.
\end{equation*}%
From \eqref{z2} we see that $Z^{\epsilon }(\tau +\epsilon )\rightarrow 0$ as 
$\epsilon \rightarrow 0^{+}$, and then from \eqref{z1} we deduce that $%
Z^{\epsilon }(t)\rightarrow 0$ as $\epsilon \rightarrow 0^{+}$, for all $t$.%
\newline
\linebreak The BSDE $\left( \ref{p-}\right) $ is linear, and we can write
the solution explicitly as follows: 
\begin{equation}
\begin{array}{l}
p^{\epsilon }(t)=\mathbb{E}[\tfrac{\Gamma (T)}{\Gamma (t)}\tfrac{\partial 
\tilde{g}}{\partial x}(\tilde{X}(T))+\tint_{t}^{T}\tfrac{\Gamma (s)}{\Gamma
(t)}\tfrac{\partial \tilde{f}}{\partial x}(s)ds|\mathcal{F}_{t}];\qquad t\in
\lbrack 0,T],%
\end{array}
\label{pfor}
\end{equation}%
where $\Gamma (t)\in \mathcal{S}^{2}$ is the solution of the linear SDE%
\begin{equation*}
\begin{cases}
d\Gamma (t) & =\Gamma (t^{-})[\tfrac{\partial \tilde{b}}{\partial x}(t)dt+%
\tfrac{\partial \tilde{\sigma}}{\partial x}(t)dB(t)+\tint_{\mathbb{R}_{0}}%
\tfrac{\partial \tilde{\gamma}}{\partial x}(t,\zeta )\tilde{N}(dt,d\zeta
)];\qquad t\in \lbrack 0,T], \\ 
\Gamma (0) & =1.%
\end{cases}%
\end{equation*}%
From this, we deduce that $p^{\epsilon }(t)\rightarrow \hat{p}(t),$ $%
q^{\epsilon }(t)\rightarrow \hat{q}(t)$ and $r^{\epsilon }(t,\zeta
)\rightarrow \hat{r}(t,\zeta )$ as $\epsilon \rightarrow 0^{+}.\square 
\newline
$\linebreak We now state and prove the main result of this part.

\begin{theorem}[Necessary maximum principle \protect\cite{AO-}]
\label{Thm.Ness} Suppose $\hat{u}\in\mathcal{A}$ is maximizing the
performance $\left( \ref{perf-}\right) $. Then for all $t\in\lbrack0,T)$ and
all bounded $\mathcal{F}_{t}$-measurable $v\in V$, we have%
\begin{equation*}
\tfrac{\partial H}{\partial u}(t,\hat{X}(t),\hat{u}(t))(v-\hat{u}(t))\leq0.
\end{equation*}
\end{theorem}

\noindent{Proof.} \quad Consider 
\begin{equation}
J(u^{\epsilon})-J(\hat{u})=I_{1}+I_{2},  \label{j}
\end{equation}
where 
\begin{equation}
I_{1}=\mathbb{E}\Big[\tint _{\tau}^{T}\{f(t,X^{\epsilon}(t),u^{\epsilon}(t))-f(t,%
\hat{X}(t),\hat{u}(t))\}dt\Big],  \label{eq3.11}
\end{equation}
and 
\begin{equation}
I_{2}=\mathbb{E}[g(X^{\epsilon}(T))-g(\hat{X}(T))].  \label{eq3.12}
\end{equation}
By the mean value theorem, we can write

\begin{equation}
I_{1}=\mathbb{E}\Big[\tint_{\tau }^{\tau +\epsilon }\{\tfrac{\partial \tilde{f}}{%
\partial x}(t)Z^{\epsilon }(t)+\tfrac{\partial \tilde{f}}{\partial u}%
(t)(u^{\epsilon }(t)-\hat{u}(t))\}dt+\tint_{\tau +\epsilon }^{T}\tfrac{%
\partial \tilde{f}}{\partial x}(t)Z^{\epsilon }(t)dt\Big],  \label{I1}
\end{equation}%

and, applying the Itô formula to $p^{\epsilon }(t)Z^{\epsilon }(t)$ and by $%
\left( \ref{p-}\right) ,\left( \ref{var1}\right) $ and $\left( \ref{var2}%
\right) $, we have 
\begin{align}
I_{2}& =\mathbb{E}\Big[\tfrac{\partial \tilde{g}}{\partial x}(\tilde{X}%
(T))Z^{\epsilon }(T)]=\mathbb{E[}p^{\epsilon }(T)Z^{\epsilon }(T)\Big]  \notag \\
& =\mathbb{E[}p^{\epsilon }(\tau +\epsilon )Z^{\epsilon }(\tau +\epsilon )] 
\notag \\
& +\mathbb{E}\Big[\tint_{\tau +\epsilon }^{T}p^{\epsilon }(t)dZ^{\epsilon
}(t)+\tint_{\tau +\epsilon }^{T}Z^{\epsilon }(t)dp^{\epsilon
}(t)+\tint_{\tau +\epsilon }^{T}d\left\langle p^{\epsilon },Z^{\epsilon
}\right\rangle (t)\Big]  \notag \\
& =\mathbb{E}\Big[p^{\epsilon }(\tau +\epsilon )(\tint_{\tau }^{\tau +\epsilon
}\{\tfrac{\partial \tilde{b}}{\partial x}(t)Z^{\epsilon }(t)+\tfrac{\partial 
\tilde{b}}{\partial u}(t)(u^{\epsilon }(t)-\hat{u}(t))\}dt  \notag \\
& +\tint_{\tau }^{\tau +\epsilon }\{\tfrac{\partial \tilde{\sigma}}{\partial
x}(t)Z^{\epsilon }(t)+\tfrac{\partial \tilde{\sigma}}{\partial u}%
(t)(u^{\epsilon }(t)-\hat{u}(t))\}dB(t)  \notag \\
& +\tint_{\tau }^{\tau +\epsilon }\tint_{\mathbb{R}_{0}}\{\tfrac{\partial 
\tilde{\gamma}}{\partial x}(t,\zeta )Z^{\epsilon }(t)+\tfrac{\partial \tilde{%
\gamma}}{\partial u}(t,\zeta )(u^{\epsilon }(t)-\hat{u}(t))\}\tilde{N}%
(dt,d\zeta ))\Big]  \notag \\
& +\mathbb{E}\Big[\tint_{\tau +\epsilon }^{T}\{p^{\epsilon }(t)\tfrac{\partial 
\tilde{b}}{\partial x}(t)Z^{\epsilon }(t)-\tfrac{\partial \tilde{H}}{%
\partial x}(t)Z^{\epsilon }(t)+q^{\epsilon }(t)\tfrac{\partial \tilde{\sigma}%
}{\partial x}(t)Z^{\epsilon }(t)  \notag \\
& +\tint_{\mathbb{R}_{0}}r^{\epsilon }(t,\zeta )\tfrac{\partial \tilde{\gamma%
}}{\partial x}(t,\zeta )Z^{\epsilon }(t)\nu (d\zeta )\}dt\Big].
\end{align}%
Using the generalized duality formula, we get%
\begin{align}
I_{2}& =\mathbb{E}\Big[\tint_{\tau }^{\tau +\epsilon }\{p^{\epsilon }(\tau
+\epsilon )(\tfrac{\partial \tilde{b}}{\partial x}(t)Z^{\epsilon }(t)+\tfrac{%
\partial \tilde{b}}{\partial u}(t)(u^{\epsilon }(t)-\hat{u}(t)))  \notag \\
& +\mathbb{E[}D_{t}p^{\epsilon }(\tau +\epsilon )|\mathcal{F}_{t}](\tfrac{%
\partial \tilde{\sigma}}{\partial x}(t)Z^{\epsilon }(t)+\tfrac{\partial 
\tilde{\sigma}}{\partial u}(t)(u^{\epsilon }(t)-\hat{u}(t)))  \notag \\
& +\tint_{\mathbb{R}_{0}}\mathbb{E[}D_{t,\zeta }p^{\epsilon }(\tau +\epsilon
)|\mathcal{F}_{t}]\{\tfrac{\partial \tilde{\gamma}}{\partial x}(t,\zeta
)Z^{\epsilon }(t)+\tfrac{\partial \tilde{\gamma}}{\partial u}(t,\zeta
)(u^{\epsilon }(t)-\hat{u}(t))\}\nu (d\zeta )\}dt\Big]  \notag \\
& -\mathbb{E}\Big[\tint_{\tau +\epsilon }^{T}\tfrac{\partial \tilde{f}}{\partial
x}(t)Z^{\epsilon }(t)dt\Big],  \label{I2-}
\end{align}%
where by the definition of $H$ $\left( \ref{h}\right) $ 
\begin{equation*}
\begin{array}{ll}
\tfrac{\partial \tilde{f}}{\partial x}(t) & =\tfrac{\partial \tilde{H}}{%
\partial x}(t)-\tfrac{\partial \tilde{b}}{\partial x}(t)p^{\epsilon }(t)-%
\tfrac{\partial \tilde{\sigma}}{\partial x}(t)q^{\epsilon }(t)-\tint_{%
\mathbb{R}_{0}}\tfrac{\partial \tilde{\gamma}}{\partial x}(t,\zeta
)r^{\epsilon }(t,\zeta )\nu (d\zeta ).%
\end{array}%
\end{equation*}%
Summing $\left( \ref{I1}\right) $ and $\left( \ref{I2-}\right) $, we obtain%
\begin{equation}
\begin{array}{ll}
I_{1}+I_{2} & =\mathbb{E[}\tint_{\tau }^{\tau +\epsilon }\{\tfrac{\partial 
\tilde{f}}{\partial x}(t)+p^{\epsilon }(\tau +\epsilon )\tfrac{\partial 
\tilde{b}}{\partial x}(t)+\mathbb{E[}D_{t}p^{\epsilon }(\tau +\epsilon )|%
\mathcal{F}_{t}]\tfrac{\partial \tilde{\sigma}}{\partial x}(t) \\ 
& +\tint_{\mathbb{R}_{0}}\mathbb{E[}D_{t,\zeta }p^{\epsilon }(\tau +\epsilon
)|\mathcal{F}_{t}]\tfrac{\partial \tilde{\gamma}}{\partial x}(t,\zeta )\nu
(d\zeta )\}Z^{\epsilon }(t)dt] \\ 
& +\mathbb{E[}\tint_{\tau }^{\tau +\epsilon }\{\tfrac{\partial \tilde{f}}{%
\partial u}(t)+p^{\epsilon }(\tau +\epsilon )\tfrac{\partial \tilde{b}}{%
\partial u}(t)+\mathbb{E[}D_{t}p^{\epsilon }(\tau +\epsilon )|\mathcal{F}%
_{t}]\tfrac{\partial \tilde{\sigma}}{\partial u}(t) \\ 
& +\tint_{\mathbb{R}_{0}}\mathbb{E[}D_{t,\zeta }p^{\epsilon }(\tau +\epsilon
)|\mathcal{F}_{t}]\tfrac{\partial \tilde{\gamma}}{\partial u}(t,\zeta )\nu
(d\zeta )\}(u^{\epsilon }(t)-\hat{u}(t))dt].%
\end{array}
\label{sum}
\end{equation}%
By the estimate of $Z^{\epsilon }$ $\left( \ref{esz}\right) $, we get 
\begin{equation}
\underset{\epsilon \rightarrow 0^{+}}{\lim }X^{\epsilon }(t)=\hat{X}(t)\text{%
; for all }t\in \lbrack \tau ,T]\text{,}  \label{estx}
\end{equation}%
and by $\left( \ref{esa}\right) $ we have 
\begin{equation}
p^{\epsilon }(t)\rightarrow \hat{p}(t)\text{, }q^{\epsilon }(t)\rightarrow 
\hat{q}(t)\text{ and }r^{\epsilon }(t,\zeta )\rightarrow \hat{r}(t,\zeta )%
\text{ when }\epsilon \rightarrow 0^{+},  \label{estpq}
\end{equation}%
where $(\hat{p},\hat{q},\hat{r})$ solves the BSDE 
\begin{equation}
\left\{ 
\begin{array}{ll}
d\hat{p}(t) & =-\tfrac{\partial \hat{H}}{\partial x}(t)dt+\hat{q}%
(t)dB(t)+\tint_{\mathbb{R}_{0}}\hat{r}(t,\zeta )\tilde{N}(dt,d\zeta );\tau
\leq t\leq T, \\ 
\hat{p}(T) & =\tfrac{\partial g}{\partial x}(\hat{X}(T)).%
\end{array}%
\right.   \label{adj}
\end{equation}%
Using the above and the assumption that $\hat{u}$ is optimal, we get%
\begin{align*}
0& \geq \underset{\epsilon \rightarrow 0^{+}}{\lim }\tfrac{1}{\epsilon }%
(J(u^{\epsilon })-J(\hat{u})) \\
& =\mathbb{E[}\{\tfrac{\partial f}{\partial u}(\tau ,\hat{X}(\tau ),\hat{u}%
(\tau ))+\hat{p}(\tau )\tfrac{\partial b}{\partial u}(\tau ,\hat{X}(\tau ),%
\hat{u}(\tau ))+\mathbb{E[}D_{\tau }\hat{p}(\tau ^{+})|\mathcal{F}_{t}]%
\tfrac{\partial \sigma }{\partial u}(\tau ,\hat{X}(\tau ),\hat{u}(\tau )) \\
& \text{ \ \ \ \ \ \ \ \ }+\tint_{\mathbb{R}_{0}}\mathbb{E[}D_{\tau ,\zeta }%
\hat{p}(\tau ^{+})|\mathcal{F}_{t}]\tfrac{\partial \gamma }{\partial u}(\tau
,\hat{X}(\tau ),\hat{u}(\tau ),\zeta )\nu (d\zeta )\}(v-\hat{u}(\tau ))],
\end{align*}%
where, by Theorem \ref{th7.3}, 
\begin{equation*}
\begin{array}{lll}
\mathbb{E[}D_{\tau }\hat{p}(\tau ^{+})|\mathcal{F}_{t}] & =\underset{%
\epsilon \rightarrow 0^{+}}{\lim }\mathbb{E[}D_{\tau }\hat{p}(\tau +\epsilon
)|\mathcal{F}_{t}] & =\hat{q}(\tau ), \\ 
\mathbb{E[}D_{\tau ,\zeta }\hat{p}(\tau ^{+})|\mathcal{F}_{t}] & =\underset{%
\epsilon \rightarrow 0^{+}}{\lim }\mathbb{E[}D_{\tau ,\zeta }\hat{p}(\tau
+\epsilon )|\mathcal{F}_{t}] & =\hat{r}(\tau ,\zeta ).%
\end{array}%
\end{equation*}%
Hence 
\begin{equation*}
\mathbb{E}[\tfrac{\partial H}{\partial u}(\tau ,\hat{X}(\tau ),\hat{u}(\tau
))(v-\hat{u}(\tau ))]\leq 0.
\end{equation*}%
Since this holds for all bounded $\mathcal{F}_{\tau }$-measurable $v$, we
conclude that%
\begin{equation*}
\tfrac{\partial H}{\partial u}(\tau ,\hat{X}(\tau ),\hat{u}(\tau ))(v-\hat{u}%
(\tau ))\leq 0\text{ for all }v.
\end{equation*}%
$\qquad \qquad \qquad \qquad \qquad \qquad \qquad \qquad \qquad \qquad
\qquad \qquad \qquad \qquad \qquad \qquad \square $

\subsubsection{Example \protect\cite{AO-}}

\noindent We now illustrate Theorem \ref{Thm.Ness} by applying it to a
linear-quadratic stochastic control problem with a constraint, as follows:%
\newline
Consider a controlled SDE of the form 
\begin{equation*}
\left\{ 
\begin{array}{ll}
dX(t) & =u(t)dt+\sigma dB(t)+\tint_{\mathbb{R}_{0}}\gamma (\zeta )\tilde{N}%
(dt,d\zeta );\quad t\in \lbrack 0,T], \\ 
X(0) & =x_{0}\in \mathbb{R}.%
\end{array}%
\right.
\end{equation*}%
Here $u\in \mathcal{A}$ is our control process (see below) and $\sigma $ and 
$\gamma $ is a given constant in $\mathbb{R}$ and function from $\mathbb{R}%
_{0}$ into $\mathbb{R}$, respectively, with 
\begin{equation*}
\tint_{\mathbb{R}_{0}}\gamma ^{2}(\zeta )\nu (d\zeta )<\infty .
\end{equation*}%
We want to control this system in such a way that we minimize its value at
the terminal time $T$ with a minimal average use of energy, measured by the
integral $\mathbb{E}[\tint_{0}^{T}u^{2}(t)dt]$ and we are only allowed to
use nonnegative controls. Thus we consider the following constrained optimal
control problem:

\begin{problem}
\label{pro} Find $\hat{u}\in\mathcal{A}$ (the set of admissible controls)
such that 
\begin{equation*}
J(\hat{u})=sup_{u\in\mathcal{A}}J(u),
\end{equation*}
where 
\begin{equation*}
J(u)=\mathbb{E}\Big[-\tfrac{1}{2}X^{2}(T)-\tfrac{1}{2}\tint _{0}^{T}u^{2}(t)dt\Big],
\end{equation*}
and $\mathcal{A}$ is the set of predictable processes $u$ such that $%
u(t)\geq0$ for all $t\in\lbrack0,T]$ and 
\begin{equation*}
\mathbb{E}\Big[\tint _{0}^{T}u^{2}(t)dt\Big]<\infty.
\end{equation*}
Thus in this case the set $V$ of admissible control values is given by $%
V=[0,\infty)$ and we can use $U=V$. The Hamiltonian is given by 
\begin{equation*}
H(t,x,u,p,q,r)=-\tfrac{1}{2}u^{2}+up+\sigma q+\tint _{\mathbb{R}%
_{0}}\gamma(\zeta)r(\zeta)\nu(d\zeta),
\end{equation*}
the adjoint BSDE for the optimal adjoint variables $\hat{p},\hat{q},\hat{r}$
is given by 
\begin{equation*}
\left\{ 
\begin{array}{ll}
d\hat{p}(t) & =\hat{q}(t)dB(t)+\tint _{\mathbb{R}_{0}}\hat{r}(t,\zeta)\tilde{%
N}(dt,d\zeta);t\in\lbrack0,T], \\ 
\hat{p}(T) & =-\hat{X}(T).%
\end{array}
\right.
\end{equation*}
Hence 
\begin{equation}
\hat{p}(t)=-\mathbb{E}[\hat{X}(T)|\mathcal{F}_{t}].
\end{equation}
Theorem \ref{Thm.Ness} states that if $\hat{u}$ is optimal, then 
\begin{equation*}
(-\hat{u}(t)+\hat{p}(t))(v-\hat{u}(t))\leq0;\quad\text{ for all }v\geq0.
\end{equation*}
From this we deduce that 
\begin{equation*}
\left\{ 
\begin{array}{cc}
\text{ (i) if }\hat{u}(t)=0, & \text{ then }\hat{u}(t)\geq\hat{p}(t), \\ 
\text{ (ii) if }\hat{u}(t)>0, & \text{ then }\hat{u}(t)=\hat{p}(t).%
\end{array}
\right.
\end{equation*}
Thus we see that we always have $\hat{u}(t)\geq\max\{\hat{p}(t),0\}$. We
claim that in fact we have equality, i.e. that 
\begin{equation*}
\hat{u}(t)=\max\{\hat{p}(t),0\}=\max\{-\mathbb{E}[\hat{X}(T)|\mathcal{F}%
_{t}],0\}.
\end{equation*}
To see this, suppose the opposite, namely that 
\begin{equation*}
\hat{u}(t)>\max\{\hat{p}(t),0\}.
\end{equation*}
Then in particular $\hat{u}(t)>0$, which by (ii) above implies that $\hat {u}%
(t)=\hat{p}(t)$, a contradiction. We summarize what we have proved as
follows:
\end{problem}

\begin{theorem}
\cite{AO-} Suppose there is an optimal control $\hat{u}\in\mathcal{A}$ for
Problem \ref{pro}. Then 
\begin{equation*}
\hat{u}(t)=\max\{\hat{p}(t),0\}=\max\{-\mathbb{E}[\hat{X}(T)|\mathcal{F}%
_{t}],0\},
\end{equation*}
where $(\hat{p},\hat{X})$ is the solution of the coupled forward-backward
SDE system given by 
\begin{align*}
& 
\begin{cases}
d\hat{X}(t) & =\max\{\hat{p}(t),0\}dt+\sigma dB(t)+\tint _{\mathbb{R}%
_{0}}\gamma(\zeta)\tilde{N}(dt,d\zeta);\quad t\in\lbrack0,T], \\ 
\hat{X}(0) & =x_{0}\in\mathbb{R},%
\end{cases}
\\
& 
\begin{cases}
d\hat{p}(t) & =\hat{q}(t)dB(t)+\tint _{\mathbb{R}_{0}}\hat{r}(t,\zeta)\tilde{%
N}(dt,d\zeta);t\in\lbrack0,T], \\ 
\hat{p}(T) & =-\hat{X}(T).%
\end{cases}%
\end{align*}
\end{theorem}

\begin{remark}
For comparison, in the case when there are no constraints on the control $u$%
, we get from the well-known solution of the classical linear-quadratic
control problem (see e.g. Øksendal \cite{O}, Example 11.2.4) that the
optimal control $u^{\ast}$ is given in feedback form by 
\begin{equation*}
u^{\ast}(t)=-\frac{\hat{X}(t)}{T+1-t};\quad t\in\lbrack0,T].
\end{equation*}
\end{remark}

\noindent \emph{EXERCISE}\newline

\begin{enumerate}
\item Let $X(t)$ satisfy the equation%
\begin{equation*}
\left\{ 
\begin{array}{ll}
dX(t) & =(b_{0}(t)+b_{1}(t)X(t))dt+(\sigma_{0}(t)+\sigma_{1}(t)X(t))dB(t) \\ 
& +\tint _{\mathbb{R}_{0}}(\gamma_{0}\left( t,\zeta\right) +\gamma_{1}\left(
t,\zeta\right) X(t))\tilde{N}(dt,d\zeta)];t\in\left[ 0,T\right] , \\ 
X(0) & =x_{0},%
\end{array}
\right.
\end{equation*}
for given $\mathbb{F}$-predictable processes $b_{0}n(t),b_{1}(t),\sigma
_{0}(t),\sigma_{1}(t),\gamma_{0}\left( t,\zeta\right) ,\gamma_{1}\left(
t,\zeta\right) $ with $\gamma_{i}\left( t,\zeta\right) \geq-1$ for $i=0,1$.%
\newline
Suppose%
\begin{align*}
\Upsilon(t) & =\exp\Big[\tint _{0}^{t}(-b_{1}(s)+\tfrac{1}{2}\sigma_{1}^{2}(s)-%
\tint _{\mathbb{R}_{0}}\{\log(1+\gamma_{1}\left( s,\zeta\right)
)-\gamma_{1}\left( s,\zeta\right) \}\nu(d\zeta))ds \\
& \text{ \ \ \ \ \ \ \ \ \ \ \ \ \ \ \ \ }-\tint
_{0}^{t}\sigma_{1}(s)dB(s)+\tint _{0}^{t}\tint _{\mathbb{R}%
_{0}}\log(1+\gamma_{1}\left( s,\zeta\right) )\tilde{N}(ds,d\zeta)\Big];t\in\left[
0,T\right] .
\end{align*}
Then the unique solution $X(t)$ is given by 
\begin{align*}
X(t) & =\Upsilon(t)^{-1}\Big[x_{0}+\tint _{0}^{t}\Upsilon(s)(b_{0}(s)+{\tint _{%
\mathbb{R}_{0}}}(\tfrac{1}{1+\gamma_{1}\left( s,\zeta\right) }%
-1)\gamma_{0}(s,\zeta )\nu(d\zeta))ds \\
& +\tint _{0}^{t}\Upsilon(s)\sigma_{0}(s)dB(s)+\tint _{0}^{t}\tint _{\mathbb{%
R}_{0}}\Upsilon(s)(\tfrac{\gamma_{0}(s,\zeta)}{1+\gamma_{1}\left(
s,\zeta\right) })\tilde{N}(ds,d\zeta)\Big];t\in\left[ 0,T\right] .
\end{align*}

\item Suppose that $Y(t)$ satisfies the linear BSDE%
\begin{equation*}
\left\{ 
\begin{array}{l}
dY(t)=-[\alpha (t)Y(t)+\beta (t)Z(t)+{\tint_{\mathbb{R}_{0}}}\eta (t,\zeta
)K(t,\zeta )\nu (d\zeta )+\gamma (t)]dt \\ 
\text{ \ \ \ \ \ \ \ \ \ \ \ \ \ }+Z(t)dB(t)+{\tint_{\mathbb{R}_{0}}}%
K(t,\zeta )\tilde{N}(dt,d\zeta )\,, \\ 
Y(t)=\xi .%
\end{array}%
\right. 
\end{equation*}%
Prove that the component of the solution $Y(t)$ can be written on its closed
formula as 
\begin{equation*}
Y(t)=\mathbb{E}\Big[(\xi \Gamma (t,T)+\tint_{t}^{T}\Gamma (t,s)\gamma (s)ds)|%
\mathcal{F}_{t}\Big],\quad t\in \left[ 0,T\right] \,,
\end{equation*}%
where $\Gamma (t,s)$ is the solution of the following linear sde 
\begin{equation*}
\left\{ 
\begin{array}{ll}
d\Gamma (t,s) & =\Gamma (t,s^{-})[\alpha (t)dt+\beta (t)dB(t)+\tint_{\mathbb{%
R}_{0}}\eta (t,\zeta )\widetilde{N}(dt,d\zeta )],\quad s\in \left[ t,T\right]
, \\ 
\Gamma (t,t) & =1\,.%
\end{array}%
\right. 
\end{equation*}
\end{enumerate}

\subsection{Stochastic Volterra integral equations (SVIEs)}

In the following we put $\triangle:=\{\left( t,s\right) \in\left[ 0,T\right]
^{2}:t\leq s\}$. We define the following spaces:

\begin{itemize}
\item $L_{y}^{2}$ consists of the $\mathbb{F}$-adapted càdlàg processes $%
Y:[0,T]\times\Omega\rightarrow\mathbb{R}$ equipped with the norm 
\begin{equation*}
\parallel Y\parallel_{L_{y}^{2}}^{2}:=\mathbb{E[}\tint
_{0}^{T}|Y(t)|^{2}dt]<\infty.
\end{equation*}

\item $L_{z}^{2}$ consists of the $\mathbb{F}$-predictable processes 
\begin{equation*}
Z:\triangle\times\Omega\rightarrow\mathbb{R},
\end{equation*}
such that $\mathbb{E[}\tint _{0}^{T}\tint _{t}^{T}\left\vert
Z(t,s)\right\vert ^{2}dsdt]<\infty$ with $s\mapsto Z(t,s)$ being $\mathbb{F}$%
-predictable on $[t,T].$ We equip $L_{z}^{2}$ with the norm 
\begin{equation*}
\parallel Z\parallel_{L_{z}^{2}}^{2}:=\mathbb{E}\Big[\tint _{0}^{T}\tint
_{t}^{T}\left\vert Z(t,s)\right\vert ^{2}dsdt\Big].
\end{equation*}
{}

\item $L_{\nu}^{2}$ consists of all Borel functions $K:%
\mathbb{R}
_{0}\rightarrow\mathbb{R},$ such that 
\begin{equation*}
\parallel K\parallel_{L_{\nu}^{2}}^{2}:=\tint _{\mathbb{R}%
_{0}}K(t,s,\zeta)^{2}\nu(d\zeta)<\infty.
\end{equation*}

\item $H_{\nu}^{2}$ consists of $\mathbb{F}$-predictable processes$\
K:\triangle\times%
\mathbb{R}
_{0}\times\Omega\rightarrow\mathbb{R},$ such that 
\begin{equation*}
\mathbb{E[}\tint _{0}^{T}\tint _{t}^{T}\tint _{\mathbb{R}_{0}}|K(t,s,%
\zeta)|^{2}\nu(d\zeta)dsdt]<\infty
\end{equation*}
and $s\mapsto K\left( t,s,\cdot\right) $ being $\mathbb{F}$-predictable on $%
[t,T].$ We equip $H_{\nu}^{2}$ with the norm 
\begin{equation*}
\parallel K\parallel_{H_{\nu}^{2}}^{2}:=\mathbb{E[}\tint _{0}^{T}\tint
_{t}^{T}\tint _{\mathbb{R}_{0}}|K(t,s,\zeta)|^{2}\nu(d\zeta)dsdt].
\end{equation*}

\item Let $L_{\mathcal{F}_{T}}^{2}[0,T]$ be the space of all processes $%
\psi:[0,T]\times\Omega\rightarrow\mathbb{R}$ and $\psi$ is $\mathcal{F}_{T}$%
-measurable for all $t\in\lbrack0,T],$ such that%
\begin{equation*}
||\psi||_{L_{\mathcal{F}_{T}}^{2}[0,T]}^{2}=\mathbb{E}\Big[\tint
_{0}^{T}|\psi(t)|^{2}dt\Big]<\infty.
\end{equation*}

\item $L_{\mathbb{F}}^{2}[0,T]$ is the space of all $\psi\in L_{\mathcal{F}%
_{T}}^{2}[0,T]$ that are $\mathbb{F}$-adapted.
\end{itemize}

Let us start by motivating what is a forward SVIE and then we will go to the
BSVIE. 

\subsubsection{A motivating example}

Stochastic Volterra integral equations (SVIEs) are a special type of
integral equations. They represent interesting models for stochastic
dynamics with memory, with applications to e.g.

\begin{itemize}
\item engineering,

\item biology (e.g. population dynamics) and

\item finance.
\end{itemize}

Moreover, they are useful tools for studying

\begin{itemize}
\item fractional Brownian motion,

\item stochastic differential equations with delay and

\item stochastic partial differential equations.
\end{itemize}



For example, let $X^{u}(t)=X(t)$ be a given cash flow, modelled by the
following stochastic Volterra integral equation: 
\begin{equation}
\begin{array}{c}
X(t)=x_{0}+\tint _{0}^{t}[b_{0}(t,s)X(s)-u(s)]ds+\tint
_{0}^{t}\sigma_{0}(s)X(s)dB(s) \\ 
+\tint _{0}^{t}\tint _{\mathbb{R}_{0}}\gamma_{0}\left( s,\zeta\right) X(s)%
\tilde{N}(ds,d\zeta);\quad t\geq0,%
\end{array}
\label{eq5.12}
\end{equation}

or, in differential form,

\begin{equation}
\left\{ 
\begin{array}{l}
dX(t)=[b_{0}(t,t)X(t)-u(t)]dt+\sigma_{0}(t)X(t)dB(t) \\ 
+\tint _{\mathbb{R}_{0}}\gamma_{0}\left( t,\zeta\right) X(t)\tilde{N}%
(dt,d\zeta)+(\int_{0}^{t}\frac{\partial b_{0}}{\partial t}%
(t,s)X(s)ds)dt;\quad t\geq0. \\ 
X(0)=x_{0}.%
\end{array}
\right.  \label{eq5.13}
\end{equation}
We see that the dynamics of $X(t)$ contains a history (or memory) term
represented by the $ds$-integral$.$\newline

We assume that $b_{0}(t,s),$ $\sigma_{0}(s)$ and $\gamma_{0}\left(
s,\zeta\right) $ are given deterministic functions of $t$, $s$, and $\zeta$,
with values in $\mathbb{R}$, and that $b_{0}(t,s)$ is continuously
differentiable with respect to $t$ for each $s$. For simplicity we assume
that these functions are bounded, and we assume that there exists $%
\varepsilon>0$ such that $\gamma_{0}(s,\zeta)\geq-1+\varepsilon$ for all $%
s,\zeta$ and the initial value $x_{0}\in \mathbb{R}$. 
Let $\mathcal{A}$ denote the set of asmissible controls $u$. We want to solve the following maximisation problem:

\begin{problem}
Find $\hat{u}\in\mathcal{A},$ such that 
\begin{equation}
\sup_{u}J(u)=J(\hat{u}),  \label{eq6.4}
\end{equation}
where 
\begin{equation}
J(u)=\mathbb{E[}\theta X(T)+\tint _{0}^{T}\log(u(t))dt],  \label{eq5.18}
\end{equation}
$\theta=\theta(\omega)$ being a given $\mathcal{F}_{T}$-measurable random
variable.
\end{problem}

We will return to this example after some general theory on optimal control
of SVIEs.

\subsubsection{Backward stochastic Volterra integral equations (BSVIEs)}

Recall that the BSDE $(Y,Z,K)$%
\begin{equation*}
\begin{array}{c}
-dY(t)=F(t,Y(t),Z(t),K(t,\cdot))dt-Z(t)dB(t) \\ 
\text{ \ \ \ \ \ }-\tint _{\mathbb{R}_{0}}K(t,\zeta)\tilde{N}(dt,d\zeta),%
\text{ \ \ }Y(T)=\varsigma,%
\end{array}%
\end{equation*}
is equivalent to

\begin{align}
Y(t) & =\varsigma+\tint _{t}^{T}F(s,Y(s),Z(s),K(s,\cdot))ds-\tint
_{t}^{T}Z(s)dB(s)  \label{bsde} \\
& \text{ \ \ \ \ \ }-\tint _{t}^{T}\tint _{\mathbb{R}_{0}}K(s,\zeta)\tilde{N}%
(ds,d\zeta).  \notag
\end{align}
The corresponding BSVIE has the form%
\begin{equation}
\begin{array}{cc}
Y(t) & =\xi(t)+\tint _{t}^{T}F(t,s,Y(s),Z(t,s),K(t,s,\cdot))ds \\ 
& -\tint _{t}^{T}Z(t,s)dB(s)-\tint _{t}^{T}\tint _{\mathbb{R}%
_{0}}K(t,s,\zeta)\tilde{N}(ds,d\zeta).%
\end{array}
\label{BSDE5}
\end{equation}

\subsubsection{Representation of solutions of BSVIE}

\begin{theorem}
\cite{AOY} Suppose that $F,Y,Z$ and $K$ are given càdlàg adapted processes
which satisfy a BSVIE of the form 
\begin{equation*}
\begin{array}{cc}
Y(t) & =\xi(t)+\tint _{t}^{T}F(t,s,Y(s),Z(t,s),K(t,s,\cdot))ds \\ 
& -\tint _{t}^{T}Z(t,s)dB(s)-\tint _{t}^{T}\tint _{\mathbb{R}%
_{0}}K(t,s,\zeta)\tilde{N}(ds,d\zeta).%
\end{array}%
\end{equation*}

Then for a.a. $t$ and $\zeta$ the following holds: 
\begin{equation*}
Z(t,s)=D_{t}Y(t^{+}):=\underset{\varepsilon\rightarrow0^{+}}{\lim}%
D_{t}Y(t+\varepsilon)\text{ (limit in }(\mathcal{S})^{\ast}),
\end{equation*}%
\begin{equation*}
Z(t,s)=\mathbb{E}[D_{t}Y(t^{+})|\mathcal{F}_{t}]:=\underset{\varepsilon
\rightarrow0^{+}}{\lim}\mathbb{E}[D_{t}Y(t+\varepsilon)|\mathcal{F}_{t}]%
\text{ (limit in }L^{2}(P)),
\end{equation*}

and 
\begin{equation*}
K(t,s,\zeta)=D_{t^{,}\zeta}Y(t^{+}):=\underset{\varepsilon\rightarrow0^{+}}{%
\lim}D_{t,\zeta}Y(t+\varepsilon)\text{ (limit in }(\mathcal{S})^{\ast}),
\end{equation*}%
\begin{equation*}
K(t,s,\zeta)=\mathbb{E}[D_{t^{,}\zeta}Y(t^{+})|\mathcal{F}_{t}]:=\underset{%
\varepsilon\rightarrow0^{+}}{\lim}\mathbb{E}[D_{t,\zeta}Y(t+\varepsilon )|%
\mathcal{F}_{t}]\text{ (limit in }L^{2}(P)).
\end{equation*}
\end{theorem}

\subsubsection{Closed formula for linear BSVIE}

\noindent Consider now the linear form. Let $\left( \Phi(t,s),0\leq t<s\leq
T\right) $ and $(\xi(s),\beta(s,\zeta);0\leq s\leq T\,,\zeta\in\mathbb{R}%
_{0})$ be given (deterministic) measurable functions of $t$, $s$, and $\zeta$%
, with values in $\mathbb{R}_{0}$. For simplicity we assume that these
functions are bounded, and we assume that there exists $\varepsilon>0$ such
that $\beta(s,\zeta)\geq-1+\varepsilon$ for all $s,\zeta$. We consider the
following linear backward stochastic Volterra integral equations in the
unknown process triplet $(Y(t),Z(t,s),K(t,s,\zeta))$: 
\begin{equation}
\begin{array}{c}
Y(t)=F(t)+\int_{t}^{T}\left[ \Phi(t,s)Y(s)+\xi(s)Z(t,s)+{\tint _{\mathbb{R}%
_{0}}}\beta(s,\zeta)K(t,s,\zeta)\nu(d\zeta)\right] ds \\ 
-\int_{t}^{T}Z(t,s)dB(s)-\int_{t}^{T}{\tint _{\mathbb{R}_{0}}}K(t,s,\zeta)%
\tilde{N}(ds,d\zeta)\,,%
\end{array}
\label{e.1.1}
\end{equation}
where $0\leq t\leq T$ and $\tilde{N}(dt,d\zeta)=N(dt,d\zeta)-\nu(d\zeta)dt$
is the compensated Poisson random measure.

To this end, we define the probability measure $Q$ by 
\begin{equation}
dQ=M(T)dP\text{ on }\mathcal{F}_{T},
\end{equation}
where 
\begin{align}
M(t) & :=\exp\Big(\int_{0}^{t}\xi(s)dB(s)-\frac{1}{2}\int_{0}^{t}\xi
^{2}(s)ds+\int_{0}^{t}\int_{\mathbb{R}_{0}}\ln(1+\beta(s,\zeta))\tilde {N}%
(ds,d\zeta)  \notag \\
& +\int_{0}^{t}\int_{\mathbb{R}_{0}}\{\ln(1+\beta(s,\zeta))-\beta
(s,\zeta)\}\nu(d\zeta)ds\Big);\quad0\leq t\leq T.
\end{align}
Then under the new probability measure $Q$ the process 
\begin{equation}
B_{Q}(t):=B(t)-\int_{0}^{t}\xi(s)ds\,,\quad0\leq t\leq T\,.  \label{e.def_bq}
\end{equation}
is a Brownian motion, and the random measure 
\begin{equation}
\widetilde{N}_{Q}(dt,d\zeta):=\widetilde{N}(dt,d\zeta)-\beta(t,\zeta)\nu(d\zeta)dt
\label{e.def_nq}
\end{equation}
is the $Q$-compensated Poisson random measure of $N(\cdot,\cdot)$, in the
sense that the process

\begin{equation*}
\widetilde{N}_{\gamma}(t):=\int_{0}^{t}\int_{\mathbb{R}_{0}}\gamma(s,\zeta )%
\widetilde{N}_{Q}(ds,d\zeta)
\end{equation*}
is a local $Q$-martingale, for all predictable processes $\gamma(t,\zeta)$
such that 
\begin{equation}
\int_{0}^{T}\int_{\mathbb{R}_{0}}\gamma^{2}(t,\zeta)\beta^{2}(t,\zeta
)\nu(d\zeta)dt<\infty.
\end{equation}
We also introduce, for $0\leq t\leq{r}\leq T$, 
\begin{equation*}
\Phi^{(1)}(t,r)=\Phi(t,r)\,,\quad\Phi^{(2)}(t,r)=\int_{t}^{r}\Phi
(t,s)\Phi(s,r)ds,
\end{equation*}
and inductively 
\begin{equation}
\Phi^{(n)}(t,r)=\int_{t}^{r}\Phi^{(n-1)}(t,s)\Phi(s,r)ds\,,\quad
n=3,4,\cdots\,.
\end{equation}

\begin{remark}
Note that if $|\Phi(t,r)|\leq C$ (constant) for all $t,r$, then by induction 
\begin{equation*}
|\Phi^{(n)}(t,r)|\leq\frac{C^{n}T^{n}}{n!}
\end{equation*}
for all $t,r,n$. Hence, 
\begin{equation*}
\sum_{n=1}^{\infty}|\Phi^{(n)}(t,r)|<\infty
\end{equation*}
for all $t,r$.
\end{remark}

\begin{theorem}
\cite{ho} Put 
\begin{equation}
\Psi(t,r)=\sum_{n=1}^{\infty}\Phi^{(n)}(t,r)\,.  \label{e.def_psi}
\end{equation}
Then we have the following explicit form of the solution triplet:

\begin{itemize}
\item[(i)] \ The $Y$ component of the solution triplet is given by 
\begin{align}
Y(t) & =\mathbb{E}_{Q}\left[ F(t)\Big|\mathcal{F}_{t}\right]
+\int_{t}^{T}\Psi(t,r)\mathbb{E}_{Q}\left[ F(r)\Big|\mathcal{F}_{t}%
\right] dr  \notag  \label{eq1.10} \\
& =\mathbb{E}_{Q}\left[ F(t)+\int_{t}^{T}\Psi(t,r)F(r)dr\Big|\mathcal{F}_{t}%
\right] \,.
\end{align}

\item[(ii)] The $Z$ and $K$ components of the solution triplet are given by
the following:\newline
Define 
\begin{equation}
U(t)=F(t)+\int_{t}^{T}\Phi(t,r)Y(r)dr-Y(t);\quad0\leq t\leq T.
\label{eq1.11}
\end{equation}
Then $Z(t,s)$ and $K(t,s,\zeta)$ can be expressed by the Hida-Malliavin
derivatives $D_{s}$ and $D_{s,\zeta}$ with respect to $B$ and $N$,
respectively, as follows: 
\begin{equation}
Z(t,s)=\mathbb{E}_{Q}[D_{s}U(t)-U(t)\int_{s}^{T}D_{s}\xi(r)dB_{Q}(r)|\mathcal{F}_{s}];\quad 0\leq t\leq s\leq T  \label{eq1.14}
\end{equation}
and 
\begin{equation}
K(t,s,\zeta)=\mathbb{E}_{Q}[U(t)(\tilde{H}_{s}-1)+\tilde{H}_{s}D_{s,\zeta
}U(t)|\mathcal{F}_{s}];\quad{0\leq}t\leq s\leq T,  \label{eq1.15}
\end{equation}
where 
\begin{align}
\tilde{H}_{s} & =\exp\Big[\int_{0}^{s}\int_{\mathbb{R}_{0}}[D_{s,x}%
\beta(r,x)+\log(1-\frac{D_{s,x}\beta(r,x)}{1-\beta(r,x)})(1-\beta
(r,x)]\nu(dx)dr  \notag \\
& +\int_{0}^{s}\int_{\mathbb{R}_{0}}\log(1-\frac{D_{s,x}\beta(r,x)}{%
1-\beta(r,x)}\tilde{N}_{Q}(dr,dx)\Big].
\end{align}
\end{itemize}
\end{theorem}

\noindent{Proof.} \quad With the processes $B_{Q}$ and $\tilde{N}_{Q}$
defined in \eqref{e.def_bq}-\eqref{e.def_nq} we can eliminate the unknowns $%
Z(t,s)$ and $K(t,s,\zeta)$ inside the first integral in \eqref{e.1.1}. More
precisely, we can rewrite equation \eqref{e.1.1} as 
\begin{equation}
Y(t)=F(t)+\int_{t}^{T}\Phi(t,s)Y(s)ds-\int_{t}^{T}Z(t,s)dB_{Q}(s)-%
\int_{t}^{T}\int_{\mathbb{R}_{0}}K(t,s,\zeta)\tilde{N}_{Q}(ds,d\zeta)\,,
\label{e.1.2}
\end{equation}
where $0\leq t\leq T$. Taking the conditional $Q$-expectation on $\mathcal{F}%
_{t}$, we get 
\begin{align}
Y(t) & =\mathbb{E}_{Q}\left[ F(t)+\int_{t}^{T}\Phi (t,s)Y(s)ds\big|\mathcal{F%
}_{t}\right]  \notag \\
& =\tilde{F}(t,t)+\int_{t}^{T}\Phi(t,s)\mathbb{E}_{Q}\left[ Y(s)\big|%
\mathcal{F}_{t}\right] ds\,,\quad0\leq t\leq T\,.  \label{1.1.5}
\end{align}
Here, and in what follows, we denote 
\begin{equation}
\tilde{F}(t,s)=\mathbb{E}_{Q}\left[ F(t)\big|\mathcal{F}_{s}\right] \,.
\end{equation}
Fix $r\in\lbrack0,t]$. Taking the conditional $Q$-expectation on $\mathcal{F}%
_{r}$ of \eqref{1.1.5}, we get 
\begin{equation*}
\mathbb{E}_{Q}\left[ Y(t)|\mathcal{F}_{r}\right] =\tilde{F}(t,r)+\int
_{t}^{T}\Phi(t,s)\mathbb{E}_{Q}\left[ Y(s)\big|\mathcal{F}_{r}\right]
ds\,,\quad r\leq t\leq T\,
\end{equation*}
Denote 
\begin{equation*}
\tilde{Y}(s)=\mathbb{E}_{Q}\left[ Y(s)|\mathcal{F}_{r}\right] \,,\quad r\leq
s\leq T\,.
\end{equation*}
Then the above equation can be written as 
\begin{equation*}
\tilde{Y}(t)=\tilde{F}(t,r)+\int_{t}^{T}\Phi(t,s)\tilde{Y}(s)ds\,,\quad
r\leq t\leq T\,.
\end{equation*}
Substituting $\tilde{Y}(s)=\tilde{F}(s,r)+\int_{s}^{T}\Phi(s,u)\tilde{Y}%
(u)du $ in the above equation, we obtain 
\begin{align}
\tilde{Y}(t) & =\tilde{F}(t,r)+\int_{t}^{T}\Phi(t,s)\left\{ \tilde {F}%
(s,r)+\int_{s}^{T}\Phi(s,u)\tilde{Y}(u)du\right\} ds  \notag \\
& =\tilde{F}(t,r)+\int_{t}^{T}\Phi(t,s)\tilde{F}(s,r)ds+\int_{t}^{T}%
\Phi^{(2)}(t,u)\tilde{Y}(u)du\,,\quad r\leq t\leq T\,,  \notag
\end{align}
By repeatedly using the above argument, we get 
\begin{align}
\tilde{Y}(t) & =\tilde{F}(t,r)+\sum_{n=1}^{\infty}\int_{t}^{T}\Phi
^{(n)}(t,u)\tilde{F}(u,r)du  \notag \\
& =\tilde{F}(t,r)+\int_{t}^{T}\Psi(t,u)\tilde{F}(u,r)du\,,
\end{align}
where $\Psi$ is defined by \eqref{e.def_psi}. 
Now substituting $\mathbb{E}_{Q}(Y(s)|\mathcal{F}_{t})=\tilde{Y}(s)$ (with $%
r=t$) into \eqref{1.1.5} we obtain part (i) of the theorem. It remains to
prove \eqref{eq1.14}-\eqref{eq1.15}. By \eqref{e.1.2} we have 
\begin{equation}
U(t)=\int_{t}^{T}Z(t,s)dB_{Q}(s)+\int_{t}^{T}\int_{\mathbb{R}%
_{0}}K(t,s,\zeta)\tilde{N}_{Q}(ds,d\zeta);\quad0\leq t\leq s\leq T.
\end{equation}
Note that by the Clark-Ocone formula under change of measure (see \cite{H})
, extended to $L^{2}(\mathcal{F}_{T},P)$ as in \cite{AaOPU}, we
get 
\begin{equation}
Z(t,s)=\mathbb{E}_{Q}[D_{s}U(t)-U(t)\int_{s}^{T}D_{s}\xi(r)dB_{Q}(r)|%
\mathcal{F}_{s}];\quad t\leq s\leq T
\end{equation}
and 
\begin{equation}
K(t,s,\zeta)=\mathbb{E}_{Q}[U(t)(\tilde{H}_{s}-1)+\tilde{H}_{s}D_{s,\zeta
}U(t)|\mathcal{F}_{s}];\quad t\leq s\leq T
\end{equation}
where 
\begin{align}
\tilde{H}_{s} & =\exp\Big[\int_{0}^{s}\int_{\mathbb{R}_{0}}[D_{s,x}%
\beta(r,x)+\log\left( 1-\frac{D_{s,x}\beta(r,x)}{1-\beta(r,x)}\right)
(1-\beta(r,x))]\nu(dx)dr  \notag \\
& +\int_{0}^{s}\int_{\mathbb{R}_{0}}\log(1-\frac{D_{s,x}\beta(r,x)}{%
1-\beta(r,x)}\tilde{N}_{Q}(dr,dx)\Big],
\end{align}
as claimed.\qquad\qquad\qquad\qquad\qquad\qquad\qquad$\square$

We illustrate our result by a specific example:

\begin{example}
\cite{ho} Let $\Phi(t, r)=\rho(r-t)$ for some bounded function $\rho$
defined on the positive half line and let $\mathcal{L} \rho(s)=
\int_{0}^{\infty}e^{-st}\rho(t) dt$ be the Laplace transform of $\rho$. Then 
$\Phi^{(n)}(t,r)=\rho _{n}(r-t)$, where $\rho_{n}=\overbrace{\rho*\cdots*\rho%
}^{n}$ is the $n$ fold convolution of $\rho$. The Laplace transform $%
\mathcal{L}\rho_{n}(s)=\left( \mathcal{L} \rho(s)\right) ^{n}$. Thus if $%
\Phi(t, r)=\rho(r-t)$, then 
\begin{equation*}
\Psi(t,r)=\bar\Psi(r-t)\,,
\end{equation*}
where the Laplace transform of $\bar\Psi$ is 
\begin{equation*}
\mathcal{L} \bar\Psi(s) =\sum_{n=1}^{\infty}\left( \mathcal{L} \rho
(s)\right) ^{n} =\frac{\mathcal{L} \rho(s)}{1-\mathcal{L} \rho(s)}\,.
\end{equation*}
In particular if $\rho(x)=e^{-x}, x>0$, then $\mathcal{L} \rho(s)=\frac {1}{%
1+s}$, which implies that $\mathcal{L} \bar\Psi(s)=\frac{1}{s}$. Thus $%
\Psi(x)=1$.
\end{example}

\subsubsection{Smoothness of the solution triplet}

\noindent It is of interest to study when the solution components $%
Z(t,s),K(t,s,\zeta)$ are smooth ($C^{1}$) with respect to $t$. Such
smoothness properties are important in the study of optimal control. It is
also important in the numerical solutions. Using the explicit form of the
solution triplet, we can give sufficient conditions for such smoothness in
the linear case.

\begin{theorem}
\cite{ho} Assume that $\xi,\beta$ are deterministic and that $F(t)$ and $%
\Phi(t,s)$ are $C^{1}$ with respect to $t$ satisfying 
\begin{equation}
\mathbb{E}_{Q}\Big[\int_{0}^{T}\left\{ \int_{t}^{T}\Big\{F^{2}(t)+\Phi
^{2}(t,s)+\left( \frac{dF(t)}{dt}\right) ^{2}+\left( \frac{\partial \Phi(t,s)%
}{\partial t}\right) ^{2}\Big\}ds\right\} dt\Big]<\infty.
\end{equation}
Then, for $t<s\leq T,$ 
\begin{align}
Z(t,s) & =\mathbb{E}_{Q}[D_{s}F(t)+\int_{t}^{T}\Phi(t,r)D_{s}Y(r)dr|\mathcal{%
F}_{s}], \\
K(t,s,\zeta) & =\mathbb{E}_{Q}[D_{s,\zeta}F(t)+\int_{t}^{T}\Phi
(t,r)D_{s,\zeta}Y(r)dr|\mathcal{F}_{s}].
\end{align}
In particular, we have 
\begin{equation}
\mathbb{E}_{Q}\left[ \int_{0}^{T}\left\{ \int_{t}^{T}\left( \frac{\partial Z%
}{\partial t}\left( t,s\right) \right) ^{2}ds\right\} dt+\int_{0}^{T}\left\{
\int_{t}^{T}\int_{\mathbb{R}_{0}}\left( \frac{\partial K}{\partial t}\left(
t,s,\zeta\right) \right) ^{2}\nu\left( d\zeta\right) ds\right\} dt\right]
<\infty.
\end{equation}
\end{theorem}

\noindent{Proof.} \quad Since $Y(t)$ is $\mathcal{F}_{t}$-measurable, we get
that $D_{s}Y(t)=D_{s,\zeta}Y(t)=0$ for all $s>t$. Hence by \eqref{eq1.11} 
\begin{equation}
\mathbb{E}_{Q}[D_{s}U(t)|\mathcal{F}_{s}]=\mathbb{E}_{Q}[D_{s}F(t)+\int
_{t}^{T}\Phi(t,r)D_{s}Y(r)dr|\mathcal{F}_{s}]
\end{equation}
and 
\begin{equation}
\mathbb{E}_{Q}[D_{s,\zeta}U(t)|\mathcal{F}_{s}]=\mathbb{E}_{Q}[D_{s,\zeta
}F(t)+\int_{t}^{T}\Phi(t,r)D_{s,\zeta}Y(r)dr|\mathcal{F}_{s}].
\end{equation}
Then the result follows from \eqref{eq1.14} and \eqref{eq1.15}.\qquad
\qquad\qquad\qquad$\square$

\subsection{Stochastic maximum principles for SVIEs}

In this section, we study stochastic maximum principles of stochastic
Volterra integral systems under partial information, i.e., the information
available to the controller is given by a sub-filtration $\mathbb{G}=\{%
\mathcal{G}_{t}\}_{t\geq0}$ such that $\mathcal{G}_{t}\subseteq\mathcal{F}%
_{t}$ for all $t\geq0.$ The set $U\subset\mathbb{R}$ is assumed to be
convex. The set of admissible controls, i.e. the strategies available to the
controller is given by a subset $\mathcal{A}_{\mathbb{G}}$ of the càdlàg, $U$%
-valued and $\mathbb{G}$-adapted processes.\newline
The state of our system $X^{u}(t)=X(t)$ satisfies the following SVIE%
\begin{equation}
\begin{array}{cc}
X(t) & =\xi(t)+\tint _{0}^{t}b\left( t,s,X(s),u(s)\right) ds+\tint
_{0}^{t}\sigma\left( t,s,X(s),u(s)\right) dB(s) \\ 
& +\tint _{0}^{t}\tint _{\mathbb{R}_{0}}\gamma\left(
t,s,X(s),u(s),\zeta\right) \tilde{N}(ds,d\zeta);t\in \lbrack0,T],%
\end{array}
\label{sde}
\end{equation}
where $b(t,s,x,u)=b(t,s,x,u,\omega):\left[ 0,T\right] ^{2}\times \mathbb{R}%
\times U\times\Omega\rightarrow\mathbb{R}$, $\sigma(t,s,x,u)=\sigma
(t,s,x,u,\omega):\left[ 0,T\right] ^{2}\times\mathbb{R}\times U\times
\Omega\rightarrow%
\mathbb{R}
$ and $\gamma(t,s,x,u,\zeta)=\gamma(t,s,x,u,\zeta,\omega):\left[ 0,T\right]
^{2}\times\mathbb{R}\times U\times%
\mathbb{R}
_{0}\times\Omega\rightarrow%
\mathbb{R}
$.\newline
\newline
The \emph{performance functional} has the form%
\begin{equation}
\begin{array}{ll}
J(u) & =\mathbb{E[}{\tint _{0}^{T}}\text{ }f(t,X(t),u(t))dt+g(X(T))],\quad
u\in\mathcal{A}_{\mathbb{G}},%
\end{array}
\label{perf}
\end{equation}
with given functions $f(t,x,u)=f(t,x,u,\omega):\left[ 0,T\right] \times%
\mathbb{R}\times U\times\Omega\rightarrow\mathbb{R}$ and $g(x)=g(x,\omega):%
\mathbb{R\times}\Omega\rightarrow\mathbb{R}$.\newline
We impose the following assumption:\newline
\textbf{Assumption A1} \newline
\emph{The processes $b,\sigma,f$ and $\gamma$ are $\mathcal{F}_{s}$-adapted
for all $s\leq t,$ and twice continuously differentiable ($C^{2}$) with
respect to $t$, $x$ and continuously differentiable ($C^{1}$) with respect
to $u$ for each $s$. The driver $g$ is assumed to be $\mathcal{F}_{T}$%
-measurable and $C^{1}$ in $x$. Moreover, all the partial derivatives are
supposed to be bounded.\newline
}Note that the performance functional (\ref{perf}) is not of Volterra type.%
\emph{\newline
\newline
}Define the \emph{Hamiltonian functional} associated to our control problem (%
\ref{sde}) and (\ref{perf}), as 
\begin{equation}
\begin{array}{l}
\mathcal{H}(t,x,v,p,q,r(\cdot)) \\ 
:=H^{0}(t,x,v,p,q,r(\cdot))+H^{1}(t,x,v,p,q,r(\cdot)),%
\end{array}
\label{eq3.3}
\end{equation}
where

\begin{equation*}
H^{0}:[0,T]\times\mathbb{R}\times U\times\mathbb{R}\times\mathbb{R}\times
L_{\nu}^{2}\rightarrow\mathbb{R}
\end{equation*}
and 
\begin{equation*}
H^{1}:[0,T]\times\mathbb{R}\times U\times\mathbb{R}\times\mathbb{R}\times
L_{\nu}^{2}\rightarrow\mathbb{R}
\end{equation*}
by 
\begin{align}
H^{0}(t,x,v,p,q,r(\cdot)) & :=f(t,x,v)+p(t)b(t,t,x,v)+q(t,t)\sigma (t,t,x,v)
\label{eq2.1} \\
& +\tint _{\mathbb{R}_{0}}r(t,t,\zeta)\gamma(t,t,x,v,\zeta)\nu(d\zeta), 
\notag
\end{align}%
\begin{align*}
H^{1}(t,x,v,p,q,r(\cdot)) & :=\tint _{t}^{T}p(s)\tfrac{\partial b}{\partial s%
}(s,t,x,v)ds+\tint _{t}^{T}q(s,t)\tfrac{\partial\sigma}{\partial s}%
(s,t,x,v)ds \\
& +\tint _{t}^{T}\tint _{\mathbb{R}_{0}}r(s,t,\zeta)\tfrac{\partial\gamma}{%
\partial s}(s,t,x,v,\zeta)\nu(d\zeta)ds.
\end{align*}
We may regard $x,p,q,r=r(\cdot)$ as generic values for the processes $%
X(\cdot),$ $p(\cdot),$ $q(\cdot),$ $r(\cdot)$, respectively.\newline
The BSVIE for the adjoint processes $p(t),q(t,s),r(t,s,\cdot)$ is defined by%
\begin{equation}
\begin{array}{cc}
p(t) & =\tfrac{\partial g}{\partial x}(X(T))+\tint _{t}^{T}\tfrac{\partial%
\mathcal{H}}{\partial x}(s)ds-\tint _{t}^{T}q(t,s)dB(s) \\ 
& -\tint _{t}^{T}\tint _{\mathbb{R}_{0}}r(t,s,\zeta)\tilde{N}%
(ds,d\zeta);t\in\lbrack0,T],%
\end{array}
\label{p}
\end{equation}
where we have used the simplified notation 
\begin{equation*}
\tfrac{\partial\mathcal{H}}{\partial x}(t)=\tfrac{\partial\mathcal{H}}{%
\partial x}(t,X(t),u(t),p(t),q(t,t),r(t,t,\cdot)).
\end{equation*}

\begin{remark}
Using the definition of $\mathcal{H}$ and the Fubini theorem, we see that
the driver in the BSVIE \eqref{p} can be explicitly written 
\begin{align}
& \int_{t}^{T}\frac{\partial\mathcal{H}}{\partial x}(s)ds=\int_{t}^{T}\Big\{%
\frac{\partial f}{\partial x}(s,x,v)+p(s)\frac{\partial b}{\partial x}%
(s,s,x,v)+\int_{s}^{T}p(z)\frac{\partial^{2}b}{\partial z\partial x}%
(z,t,x,v)dz  \notag \\
& +q(s,s)\frac{\partial\sigma}{\partial x}(s,s,x,v)+\int_{s}^{T}q(z,t)\frac{%
\partial^{2}\sigma}{\partial z\partial x}(z,t,x,v)dz  \notag \\
& +\int_{\mathbb{R}_{0}}r(s,s,\zeta)\frac{\partial\gamma}{\partial x}%
(s,s,x,v,\zeta)\nu(d\zeta)+\int_{s}^{T}\int_{\mathbb{R}_{0}}r(z,t,\zeta )%
\frac{\partial^{2}\gamma}{\partial z\partial x}(z,t,x,v,\zeta)\nu (d\zeta)dz%
\Big\}ds  \notag \\
& =\int_{t}^{T}\Big\{\frac{\partial f}{\partial x}(s,x,v)+p(s)\Big[\frac {%
\partial b}{\partial x}(s,s,x,v)+(s-t)\frac{\partial^{2}b}{\partial
s\partial x}(s,t,x,v)\Big]  \notag \\
& +q(s,s)\frac{\partial\sigma}{\partial x}(s,s,x,v)+(s-t)q(s,t)\frac {%
\partial^{2}\sigma}{\partial s\partial x}(s,t,x,v)  \notag \\
& +\int_{\mathbb{R}_{0}}\Big[r(s,s,\zeta)\frac{\partial\gamma}{\partial x}%
(s,s,x,v,\zeta)+(s-t)r(s,t,x,v,\zeta)\frac{\partial^{2}\gamma}{\partial
s\partial x}(s,t,x,v,\zeta)\Big]\nu(d\zeta)\Big\}ds.
\end{align}
From this it follows by Theorem 3.1 in Agram \textit{et al }\cite{AOY}, that
we have existence and uniqueness of the solution of equation \eqref{p}. 
\newline
From now on we also make the following assumption:\newline
Note that from equation (\ref{sde}), we get the following equivalent
formulation, for each $(t,s)\in\lbrack0,T]^{2},$ 
\begin{equation}
\begin{array}{ll}
dX(t) & =\xi^{\prime}(t)dt+b\left( t,t,X(t),u(t)\right) dt+(\tint _{0}^{t}%
\tfrac{\partial b}{\partial t}\left( t,s,X(s),u(s)\right) ds)dt \\ 
& +\sigma\left( t,t,X(t),u(t)\right) dB(t)+(\tint _{0}^{t}\tfrac{%
\partial\sigma}{\partial t}\left( t,s,X(s),u(s)\right) dB(s))dt \\ 
& +\tint _{\mathbb{R}_{0}}\gamma\left( t,t,X(t),u(t),\zeta\right) \tilde{N}%
(dt,d\zeta)+(\tint _{0}^{t}\tint _{\mathbb{R}_{0}}\tfrac{\partial\gamma}{%
\partial t}\left( t,s,X(s),u(s),\zeta\right) \tilde{N}(ds,d\zeta))dt,%
\end{array}
\label{eq3.6}
\end{equation}
and from equation (\ref{p}) under assumption $\emph{A2}$, we have the
following differential form 
\begin{equation}
\left\{ 
\begin{array}{ll}
dp(t) & =-[\tfrac{\partial\mathcal{H}}{\partial x}(t)+\tint _{t}^{T}\tfrac{%
\partial q}{\partial t}(t,s)dB(s)+\tint _{t}^{T}\tint _{\mathbb{R}_{0}}%
\tfrac{\partial r}{\partial t}(t,s,\zeta)\tilde{N}(ds,d\zeta)]dt \\ 
& +q(t,t)dB(t)+\tint _{\mathbb{R}_{0}}r(t,t,\zeta)\tilde{N}(dt,d\zeta), \\ 
p(T) & =\tfrac{\partial g}{\partial x}(X(T)).%
\end{array}
\right.  \label{eq3.7}
\end{equation}
\end{remark}

\begin{remark}
Assumption $\emph{A2}$ is verified in a subclass of linear BSVIE with jumps,
as we will see in section 5. For more details, we refer to Hu and Øksendal 
\cite{ho}.
\end{remark}

\subsubsection{A sufficient maximum principle}

\noindent We now state and prove a sufficient version of the maximum
principle approach (a verification theorem).

\begin{theorem}[Sufficient maximum principle \protect\cite{AOY}]
Let $\hat{u}\in\mathcal{A}_{\mathbb{G}},$ with corresponding solutions $\hat{%
X}(t),$ $\left( \hat{p}(t),\hat{q}(t,s),\hat {r}(t,s,\cdot)\right) $ of (\ref%
{sde}) and (\ref{p}) respectively. Assume that

\begin{itemize}
\item The functions%
\begin{equation*}
x\mapsto g(x),
\end{equation*}
and 
\begin{equation*}
(x,u)\mapsto\mathcal{H}(t,x,u,\hat{p},\hat{q},\hat{r}(\cdot))
\end{equation*}
are concave.

\item (The maximum condition) 
\begin{align}
& \underset{v\in\mathbb{U}}{\sup}\text{ }\mathbb{E[}\mathcal{H}(t,\hat {X}%
(t),v,\hat{p}(t),\hat{q}(t,t),\hat{r}(t,t,\cdot))|{\mathcal{G}}_{t}]  \notag
\\
& =\mathbb{E[}\mathcal{H}(t,\hat{X}(t),\hat{u}(t),\hat{p}(t),\hat {q}(t,t),%
\hat{r}(t,t,\cdot))|{\mathcal{G}}_{t}]\text{ }\forall t\text{ }P\text{-a.s.}
\label{eq3.9}
\end{align}
\newline
Then $\hat{u}$ is an optimal control for our problem.
\end{itemize}
\end{theorem}

\noindent{Proof.} \quad By considering a sequence of stopping times
converging upwards to $T$, we see that we may assume that all the $dB$- and $%
\tilde{N}$- integrals in the following are martingales and hence have
expectation $0$. \newline
\newline
Choose $u\in\mathcal{A}_{\mathbb{G}},$ we want to prove that $J(u)\leq J(%
\hat{u})$.\newline
By the definition of the cost functional (\ref{perf})$,$ we have 
\begin{equation}
J(u)-J(\hat{u})=I_{1}+I_{2},  \label{eq2.8}
\end{equation}
where we have used the shorthand notations \newline
\begin{equation*}
I_{1}=\mathbb{E}\Big[\tint _{0}^{T}\tilde{f}\left( t\right) dt\Big],\quad I_{2}=%
\mathbb{E[}\tilde{g}(T)],
\end{equation*}
\newline
and 
\begin{equation*}
\tilde{f}\left( t\right) =f(t)-\hat{f}(t),
\end{equation*}
with 
\begin{equation*}
\begin{array}{ll}
f(t) & =f\left( t,X(t),u(t)\right) , \\ 
\hat{f}\left( t\right) & =f(t,\hat{X}(t),\hat{u}(t)),%
\end{array}%
\end{equation*}
and similarly for $b(t,t)=b\left( t,t,X(t),u(t)\right) ,$ and the other
coefficients. By the definition of the Hamiltonian $(\ref{eq2.1})$, we get 
\begin{equation}
I_{1}=\mathbb{E}\Big[\tint _{0}^{T}\{\tilde{H}^{0}(t)-\hat{p}(t)\tilde{b}(t,t)-%
\hat{q}(t,t)\tilde{\sigma}(t,t)-\tint _{\mathbb{R}_{0}}\hat{r}(t,t,\zeta)%
\tilde{\gamma}(t,t,\zeta)\nu(d\zeta)\}dt\Big],  \label{i1}
\end{equation}
where $\tilde{H}^{0}(t)=H^{0}(t)-\hat{H}^{0}(t)$ with%
\begin{equation*}
\begin{array}{ll}
H^{0}(t) & =H^{0}(t,X(t),u(t),\hat{p}(t),\hat{q}(t,t),\hat{r}(t,t,\cdot)),
\\ 
\hat{H}^{0}(t) & =\hat{H}^{0}(t,\hat{X}(t),\hat{u}(t),\hat{p}(t),\hat {q}%
(t,t),\hat{r}(t,t,\cdot)).%
\end{array}%
\end{equation*}
By the concavity of $g$ and the terminal value of the BSVIE $\left( \ref{p}%
\right) $, we obtain%
\begin{equation*}
\begin{array}{lll}
I_{2} & \leq\mathbb{E}[\tfrac{\partial\hat{g}}{\partial x}(T)\tilde{X}(T)] & 
=\mathbb{E}[\hat{p}(T)\tilde{X}(T)].%
\end{array}%
\end{equation*}
Applying the Itô formula to $\hat{p}(t)\tilde{X}(t)$, we get 
\begin{align}
I_{2} & \leq\mathbb{E}\Big[\hat{p}(T)\tilde{X}(T)\Big]  \notag \\
& =\mathbb{E}\Big[\tint _{0}^{T}\hat{p}(t)\{\tilde{b}(t,t)+\tint _{0}^{t}\tfrac{%
\partial\tilde{b}}{\partial t}(t,s)ds+\tint _{0}^{t}\tfrac{\partial\tilde{%
\sigma}}{\partial t}\left( t,s\right) dB(s)  \notag \\
& +\tint _{0}^{t}\tint _{\mathbb{R}_{0}}\tfrac{\partial\tilde{\gamma}}{%
\partial t}\left( t,s,\zeta\right) \tilde {N}(ds,d\zeta)\}dt+\tint _{0}^{T}%
\tilde{X}(t)\{-\tfrac{\partial\widehat{\mathcal{H}}}{\partial x}(t)+\tint
_{t}^{T}\tfrac{\partial\hat{q}}{\partial t}(t,s)dB(s)  \notag \\
& +\tint _{t}^{T}\tint _{\mathbb{R}_{0}}\tfrac{\partial\hat{r}}{\partial t}%
(t,s,\zeta)\tilde{N}(ds,d\zeta)\}dt+\tint _{0}^{T}\hat{q}(t,t)\tilde{\sigma}%
(t,t)dt+\tint _{0}^{T}\tint _{\mathbb{R}_{0}}\hat{r}(t,t,\zeta)\tilde{\gamma}%
(t,t,\zeta)\nu(d\zeta)dt\Big].  \label{I2}
\end{align}
By the Fubini theorem, we get 
\begin{equation}
\tint _{0}^{T}\hat{p}(t)\Big(\tint _{0}^{t}\tfrac{\partial\tilde{b}}{\partial t}%
(t,s)ds\Big)dt=\tint _{0}^{T}\Big(\tint _{s}^{T}\hat{p}(t)\tfrac{\partial\tilde{b}}{%
\partial t}(t,s)dt\Big)ds=\tint _{0}^{T}\Big(\tint _{t}^{T}\hat{p}(s)\tfrac{\partial%
\tilde{b}}{\partial s}(s,t)ds\Big)dt.  \label{eq2.11}
\end{equation}
The generalized duality formula for the Brownian motion, yields 
\begin{align*}
\mathbb{E}\Big[\tint _{0}^{T}\hat{p}(t)(\tint _{0}^{t}\tfrac{\partial\tilde{%
\sigma}}{\partial t}(t,s)dB(s))dt\Big] & =\tint _{0}^{T}\mathbb{E}\Big[\tint _{0}^{t}%
\hat{p}(t)\tfrac{\partial\tilde{\sigma}}{\partial t}(t,s)dB(s)\Big]dt \\
& =\tint _{0}^{T}\mathbb{E}\Big[\tint _{0}^{t}\mathbb{E}[D_{s}\hat{p}(t)|%
\mathcal{F}_{s}]\tfrac{\partial\tilde{\sigma}}{\partial t}(t,s)ds\Big]dt.
\end{align*}
Fubini's theorem gives%
\begin{align*}
\mathbb{E}\Big[\tint _{0}^{T}\hat{p}(t)(\tint _{0}^{t}\tfrac{\partial\tilde{%
\sigma}}{\partial t}(t,s)dB(s))dt\Big] & =\tint _{0}^{T}\mathbb{E}\Big[\tint _{s}^{T}%
\mathbb{E}[D_{s}\hat{p}(t)|\mathcal{F}_{s}]\tfrac{\partial\tilde{\sigma}}{%
\partial t}(t,s)dt\Big]ds \\
& =\mathbb{E}\Big[\tint _{0}^{T}\tint _{t}^{T}\mathbb{E}[D_{t}\hat{p}(s)|%
\mathcal{F}_{t}]\tfrac{\partial\tilde{\sigma}}{\partial s}(s,t)dsdt\Big],
\end{align*}
and by equality $\left( \ref{eq2.25}\right) $, we end up with%
\begin{equation}
\mathbb{E}\Big[\tint _{0}^{T}\hat{p}(t)(\tint _{0}^{t}\tfrac{\partial\tilde{%
\sigma}}{\partial t}(t,s)dB(s))dt\Big]=\mathbb{E}\Big[\tint _{0}^{T}\tint _{t}^{T}%
\hat{q}(s,t)\tfrac{\partial\tilde{\sigma}}{\partial s}(s,t)dsdt\Big].
\label{eq2.12}
\end{equation}
Doing similar considerations as for the Brownian setting for the jumps, such
as the Fubini theorem, the generalized duality formula for jumps, we obtain

\begin{align}
\mathbb{E}\Big[\tint _{0}^{T}(\tint _{0}^{t}\tint _{\mathbb{R}_{0}}\hat{p}(t)%
\tfrac{\partial\tilde{\gamma}}{\partial t}(t,s,\zeta)\tilde {N}%
(ds,d\zeta))dt\Big] & =\tint _{0}^{T}\mathbb{E}\Big[\tint _{0}^{t}\tint _{\mathbb{R}%
_{0}}\hat{p}(t)\tfrac{\partial\tilde{\gamma}}{\partial t}(t,s,\zeta)\tilde {N%
}(ds,d\zeta))\Big]dt  \notag \\
& =\tint _{0}^{T}\mathbb{E}\Big[\tint _{0}^{t}\tint _{\mathbb{R}_{0}}\mathbb{E}%
[D_{s,\zeta}\hat{p}(t)|\mathcal{F}_{s}]\tfrac{\partial\tilde{\gamma }}{%
\partial t}(t,s,\zeta)\nu(d\zeta)ds\Big]dt  \notag \\
& =\tint _{0}^{T}\mathbb{E}\Big[\tint _{s}^{T}\tint _{\mathbb{R}_{0}}\mathbb{E}%
\Big[D_{s,\zeta}\hat{p}(t)|\mathcal{F}_{s}]\tfrac{\partial\tilde{\gamma }}{%
\partial t}(t,s,\zeta)\nu(d\zeta)dt\Big]ds  \notag \\
& =\mathbb{E}\Big[\tint _{0}^{T}\tint _{t}^{T}\tint _{\mathbb{R}_{0}}\mathbb{E}%
\Big[D_{t,\zeta}\hat{p}(s)|\mathcal{F}_{t}]\tfrac{\partial\tilde{\gamma }}{%
\partial s}(s,t,\zeta)\nu(d\zeta)dsdt\Big]  \notag \\
& =\mathbb{E}\Big[\tint _{0}^{T}\tint _{t}^{T}\tint _{\mathbb{R}_{0}}\hat{r}%
(s,t,\zeta)\tfrac{\partial\tilde{\gamma}}{\partial s}(s,t,\zeta
)\nu(d\zeta)dsdt\Big].  \label{eq2.13}
\end{align}
Substituting $\left( \ref{eq2.11}\right) ,\left( \ref{eq2.12}\right) $ and $%
\left( \ref{eq2.13}\right) $ combined with $\left( \ref{eq3.3}\right) $ in $%
\left( \ref{eq2.8}\right) $, yields%
\begin{equation}
J(u)-J(\hat{u})\leq\mathbb{E}\Big[\tint _{0}^{T}\{\mathcal{H}(t)-\widehat{%
\mathcal{H}}(t)-\tfrac{\partial\widehat{\mathcal{H}}}{\partial x}(t)\tilde{X}%
(t)\}dt\Big].  \notag
\end{equation}
\newline
By the concavity of $\mathcal{H},$ we have%
\begin{equation*}
\mathcal{H}(t)-\widehat{\mathcal{H}}(t)\leq\tfrac{\partial\widehat {\mathcal{%
H}}}{\partial x}(t)\tilde{X}(t)+\tfrac{\partial\widehat{\mathcal{H}}}{%
\partial u}(t)\tilde{u}(t).
\end{equation*}
Hence, since $u=\hat{u}$ is $\mathbb{G}$-adapted and maximizes the
conditional Hamiltonian, 
\begin{align}
& J(u)-J(\hat{u})\leq\mathbb{E}\Big[\tint _{0}^{T}\tfrac{\partial\mathcal{H}}{%
\partial u}(t)(u(t)-\hat{u}(t))dt\Big]  \notag \\
& =\mathbb{E}\Big[\tint _{0}^{T}\mathbb{E}[\tfrac{\partial\mathcal{H}}{\partial u%
}(t)|\mathcal{G}_{t}](u(t)-\hat{u}(t))dt\Big]\leq0,
\end{align}
which means that $\hat{u}$ is an optimal control. \hfill$\square$ \bigskip

\subsubsection{A necessary maximum principle}

\noindent Suppose that a control $u\in\mathcal{A}_{\mathbb{G}}$ is optimal
and that $\beta${\normalsize $\in\mathcal{A}_{\mathbb{G}}.$ If the function }%
$\lambda${\normalsize $\longmapsto J(u+\lambda\beta)$ is well-defined and
differentiable on a neighbourhood of $0$, then 
\begin{equation*}
\tfrac{d}{d\lambda}J(u+\lambda\beta)\mid_{\lambda=0}=0.
\end{equation*}
Under a set of suitable assumptions on the coefficients, we will show that%
\begin{equation*}
\tfrac{d}{d\lambda}J(u+\lambda\beta)\mid_{\lambda=0}=0
\end{equation*}
is equivalent to%
\begin{equation*}
\mathbb{E}[\tfrac{\partial\mathcal{H}}{\partial u}\left( t\right) \mathcal{%
\mid G}_{t}]=0\text{ \ }P-\text{a.s. for each }t\in\lbrack0,T].
\end{equation*}
}\newline
The details are as follows:\newline
For each given $t\in \lbrack0,T],$ let $\eta=\eta(t)$ be a bounded $\mathcal{%
G}_{t}$-measurable random variable, let $h\in\lbrack T-t,T]$ and define

\begin{equation}
\beta(s):=\eta1_{\left[ t,t+h\right] }(s);s\in\left[ 0,T\right] .
\label{eq4.1}
\end{equation}
Assume that 
\begin{equation}
u+\lambda\beta\in\mathcal{A}_{\mathbb{G}},
\end{equation}
for all $\beta$ and all $u\in\mathcal{A}_{\mathbb{G}}$, and all non-zero $%
\lambda$ sufficiently small. Assume that the \emph{derivative process} $Y(t)$%
, defined by 
\begin{equation}
Y(t)=\tfrac{d}{d\lambda}X^{(u+\lambda\beta)}(t)|_{\lambda=0},  \label{1.13}
\end{equation}
exists.\newline
Then we see that%
\begin{align*}
Y(t) & =\tint _{0}^{t}\Big(\tfrac{\partial b}{\partial x}(t,s)Y(s)+\tfrac{%
\partial b}{\partial u}(t,s)\beta(s)\Big)ds \\
& +\tint _{0}^{t}\Big(\tfrac{\partial\sigma}{\partial x}(t,s)Y(s)+\tfrac{%
\partial\sigma}{\partial u}(t,s)\beta(s)\Big)dB(s) \\
& +\tint _{0}^{t}\tint _{\mathbb{R}_{0}}\Big(\tfrac{\partial\gamma}{\partial x}%
(t,s,\zeta)Y(s)+\tfrac{\partial\gamma }{\partial u}(t,s,\zeta)\beta(s)\Big)%
\tilde{N}(ds,d\zeta),
\end{align*}
and hence

\begin{align}
dY(t) & =\Big[\tfrac{\partial b}{\partial x}(t,t)Y(t)+\tfrac{\partial b}{%
\partial u}(t,t)\beta(t)+\tint _{0}^{t}(\tfrac{\partial^{2}b}{\partial
t\partial x}(t,s)Y(s)+\tfrac{\partial^{2}b}{\partial t\partial u}%
(t,s)\beta(s))ds  \notag \\
& +\tint _{0}^{t}\Big(\tfrac{\partial^{2}\sigma}{\partial t\partial x}(t,s)Y(s)+%
\tfrac{\partial ^{2}\sigma}{\partial t\partial u}(t,s)\beta(s)\Big)dB(s)  \notag
\\
& +\tint _{0}^{t}\tint _{\mathbb{R}_{0}}(\tfrac{\partial^{2}\gamma}{\partial
t\partial x}(t,s,\zeta)Y(s)+\tfrac {\partial^{2}\gamma}{\partial t\partial u}%
(t,s,\zeta)\beta(s))\tilde {N}(ds,d\zeta)\Big]dt  \notag \\
& +\Big(\tfrac{\partial\sigma}{\partial x}(t,t)Y(t)+\tfrac{\partial\sigma }{%
\partial u}(t,t)\beta(t)\Big)dB(t)  \notag \\
& +\tint _{\mathbb{R}_{0}}\Big(\tfrac{\partial\gamma}{\partial x}(t,t,\zeta)Y(t)+%
\tfrac{\partial\gamma }{\partial u}(t,t,\zeta)\beta(t)\Big)\tilde{N}(dt,d\zeta).
\label{1.14}
\end{align}
We are now ready to formulate the result:

\begin{theorem}[Necessary maximum principle \protect\cite{AOY}]
Suppose that $\hat{u}\in$ $\mathcal{A}_{\mathbb{G}}$ is such that, for all $%
\beta$ as in \eqref{eq4.1}, 
\begin{equation}
\tfrac{d}{d\lambda}J(\hat{u}+\lambda\beta)|_{\lambda=0}=0  \label{eq4.5}
\end{equation}
and the corresponding solution $\hat{X}(t),(\hat{p}(t),\hat{q}(t,t),\hat {r}%
(t,t,\cdot))$ of (\ref{sde}) and (\ref{p}) exists. Then, 
\begin{equation}
\mathbb{E[}\tfrac{\partial\mathcal{H}}{\partial u}(t)|\mathcal{G}_{t}]_{u=%
\hat{u}(t)}=0.  \label{eq4.6}
\end{equation}
Conversely, if \eqref{eq4.6} holds, then \eqref{eq4.5} holds.
\end{theorem}

\noindent{Proof.} \quad By considering a suitable increasing family of
stopping times converging to $T$, we may assume that all the local
martingales ($dB$- and $\tilde{N}$- integrals) appearing in the proof below
are martingales. For simplicity of notation we drop the "hat" everywhere and write $%
u $ in stead of $\hat{u}$, $X$ in stead of $\hat{X}$ etc in the following.
Consider 
\begin{equation}
\begin{array}{l}
\tfrac{d}{d\lambda}J(u+\lambda\beta)|_{\lambda=0} \\ 
=\mathbb{E[}\tint _{0}^{T}\{\tfrac{\partial f}{\partial x}(t)Y(t)+\tfrac{%
\partial f}{\partial u}(t)\beta(t)\}dt+\tfrac{\partial g}{\partial x}%
(X(T))Y(T)].%
\end{array}
\label{1.15}
\end{equation}
Applying the Itô formula, we get 
\begin{align*}
\begin{array}{l}
\mathbb{E[}\tfrac{\partial g}{\partial x}(X(T))Y(T)]=\mathbb{E[}p(T)Y(T)] \\ 
=\mathbb{E[}\tint _{0}^{T}p(t)(\tfrac{\partial b}{\partial x}(t,t)Y(t)+%
\tfrac{\partial b}{\partial u}(t,t)\beta(t))dt \\ 
+\tint _{0}^{T}p(t)\{\tint _{0}^{t}(\tfrac{\partial^{2}b}{\partial t\partial
x}(t,s)Y(s)+\tfrac{\partial^{2}b}{\partial t\partial u}(t,s)\beta(s))ds\}dt
\\ 
+\tint _{0}^{T}p(t)\{\tint _{0}^{t}(\frac{\partial^{2}\sigma}{\partial
t\partial x}(t,s)Y(s)+\frac{\partial ^{2}\sigma}{\partial t\partial u}%
(t,s)\beta(s))dB(s)\}dt \\ 
+\tint _{0}^{T}p(t)\{\tint _{0}^{t}\tint _{\mathbb{R}_{0}}(\tfrac{%
\partial^{2}\gamma}{\partial t\partial x}(t,s,\zeta)Y(s)+\tfrac {%
\partial^{2}\gamma}{\partial t\partial u}(t,s,\zeta)\beta(s))\tilde {N}%
(ds,d\zeta)\}dt \\ 
-\tint _{0}^{T}Y(t)\tfrac{\partial\mathcal{H}}{\partial x}(t)dt+\tint
_{0}^{T}q(t,s)(\tfrac{\partial\sigma}{\partial x}(t,t)Y(t)+\tfrac{%
\partial\sigma }{\partial u}(t,t)\beta(t))dt\\
+\tint _{0}^{T}\tint _{\mathbb{R}_{0}}r(t,s,\zeta)(\tfrac{\partial\gamma}{%
\partial x}(t,t,\zeta)Y(t)+\tfrac {\partial\gamma}{\partial u}%
(t,t,\zeta)\beta(t))\nu(d\zeta)dt].
\end{array}{|}
\end{align*}
From $\left( \ref{eq2.12}\right) $ and $\left( \ref{eq2.13}\right) $, we have%
\begin{align*}
\begin{array}{l}
\mathbb{E}\left[ p(T)Y(T)\right] \\ 
=\mathbb{E}[\tint _{0}^{T}\{\tfrac{\partial b}{\partial x}(t,t)p(t)+\tint
_{t}^{T}(\tfrac{\partial^{2}b}{\partial s\partial x}(s,t)p(s)+\tfrac{%
\partial ^{2}\sigma}{\partial s\partial x}(s,t)q(s,t) \\ 
+\tint _{\mathbb{R}_{0}}\frac{\partial^{2}\gamma}{\partial s\partial x}%
(s,t,\zeta)r(s,t,\zeta )\nu(d\zeta))ds\}Y(t)dt \\ 
\mathbb{+}\tint _{0}^{T}\{\tfrac{\partial b}{\partial u}(t,t)p(t)+\tint
_{t}^{T}(\frac{\partial^{2}b}{\partial s\partial u}(s,t)p(s)+\frac{%
\partial^{2}\sigma }{\partial s\partial u}(s,t)q(s,t) \\ 
+\tint _{\mathbb{R}_{0}}\frac{\partial^{2}\gamma}{\partial s\partial u}%
(s,t,\zeta)r(s,t,\zeta )\nu(d\zeta))ds\}\beta(t)dt \\ 
-\tint _{0}^{T}\frac{\partial\mathcal{H}}{\partial x}(t)Y(t)dt+\tint
_{0}^{T}(\frac{\partial\sigma}{\partial x}(t,t)Y(t)+\frac{\partial\sigma}{%
\partial u}(t,t)\beta(t))q(t,t)dt%
\\
+\tint _{0}^{T}\tint _{\mathbb{R}_{0}}(\tfrac{\partial\gamma}{\partial x}%
(t,t,\zeta)Y(t)+\tfrac{\partial\gamma }{\partial u}(t,t,\zeta)%
\beta(t))r(t,t,\zeta)\nu(d\zeta)dt].
\end{array}
\end{align*}

Using the definition of $\mathcal{H}$ in (\ref{eq3.3}) and the definition of 
$\beta$, we obtain

\begin{equation}
\tfrac{d}{d\lambda}J(u+\lambda\beta)|_{\lambda=0}=\mathbb{E}\Big[\tint _{0}^{T}%
\tfrac{\partial\mathcal{H}}{\partial u}(s)\beta(s)ds\Big]=\mathbb{E}\Big[\tint
_{t}^{t+h}\tfrac{\partial\mathcal{H}}{\partial u}(s)ds\alpha \Big].  \label{eq4.8}
\end{equation}
Now suppose that 
\begin{equation}
\tfrac{d}{d\lambda}J(u+\lambda\beta)|_{\lambda=0}=0.  \label{eq4.9}
\end{equation}
Differentiating the right-hand side of \eqref{eq4.8} at $h=0$, we get%
\begin{equation*}
\mathbb{E}[\tfrac{\partial\mathcal{H}}{\partial u}(t)\eta]=0.
\end{equation*}
Since this holds for all bounded $\mathcal{G}_{t}$-measurable $\eta$, we have%
\begin{equation}
\mathbb{E[}\tfrac{\partial\mathcal{H}}{\partial u}(t)|\mathcal{G}_{t}]=0.
\label{eq4.10}
\end{equation}
Conversely, if we assume that (\ref{eq4.10}) holds, then we obtain (\ref%
{eq4.9}) by reversing the argument we used to obtain (\ref{eq4.8}).

\hfill$\square$ \bigskip

\noindent \emph{EXERCISE}\newline
Let $X^{u}(t)=X(t)$ be a given cash flow, modelled by the following
stochastic Volterra equation: 
\begin{equation}
\begin{array}{c}
X(t)=x_{0}+\tint _{0}^{t}[b_{0}(t,s)X(s)-u(s)]ds+\tint
_{0}^{t}\sigma_{0}(s)X(s)dB(s) \\ 
+\tint _{0}^{t}\tint _{\mathbb{R}_{0}}\gamma_{0}\left( s,\zeta\right) X(s)%
\tilde{N}(ds,d\zeta);\quad t\geq0,%
\end{array}
\label{eq5.12}
\end{equation}
or, in differential form,%
\begin{equation}
\left\{ 
\begin{array}{l}
dX(t)=[b_{0}(t,t)X(t)-u(t)]dt+\sigma_{0}(t)X(t)dB(t) \\ 
+\tint _{\mathbb{R}_{0}}\gamma_{0}\left( t,\zeta\right) X(t)\tilde{N}%
(dt,d\zeta)+[\int_{0}^{t}\frac{\partial b_{0}}{\partial t}%
(t,s)X(s)ds]dt;\quad t\geq0. \\ 
X(0)=x_{0}.%
\end{array}
\right.  \label{eq5.13}
\end{equation}
We see that the dynamics of $X(t)$ contains a history or memory term
represented by the $ds$-integral$.$\newline
We assume that $b_{0}(t,s),$ $\sigma_{0}(s)$ and $\gamma_{0}\left(
s,\zeta\right) $ are given 
deterministic functions of $t$, $s$, and $\zeta$, with values in $\mathbb{R}$%
, and that $b_{0}(t,s)$ is continuously differentiable with respect to $t$
for each $s$. For simplicity we assume that these functions are bounded, and
we assume that there exists $\varepsilon>0$ such that $\gamma_{0}(s,\zeta
)\geq-1+\varepsilon$ for all $s,\zeta$ and the initial value $x_{0}\in%
\mathbb{R}
$. We want to solve the following maximisation problem: \newline
Find $\hat {u}\in\mathcal{A_{\mathbb{G}}},$ such that 
\begin{equation}
\sup_{u}J(u)=J(\hat{u}),  \label{eq6.4}
\end{equation}
where 
\begin{equation}
J(u)=\mathbb{E}\Big[\theta X(T)+\tint _{0}^{T}\log(u(t))dt\Big].  \label{eq5.18}
\end{equation}
Here $\theta=\theta(\omega)$ is a given $\mathcal{F}_{T}$-measurable random
variable.\\

$\bold{Acknowledgments}.$\\
We are grateful to José Luis da Silva for his valuable comments.

\end{document}